\documentclass[
        parskip = half, 
        twoside,
specimen_mek        footinclude = false, 
       a4paper ]{scrartcl}

\usepackage{ifthen} 

\newboolean{boldeqnitalic}
\setboolean{boldeqnitalic}{true}
\usepackage{stylefiles/defdb1} 

\usepackage{stylefiles/rublkm_paper} 

\newcommand{\figpath}{figures}       
                                                     
\usepackage[linesnumbered]{algorithm2e}
\usepackage{tikz}
\usepackage{pgfplots}

\usetikzlibrary{shapes,arrows}

\begin{document}


\setcounter{page}{1}

\begin{center}
{\LARGE Simulation of crack propagation through voxel-based, heterogeneous structures based on eigenerosion and finite cells}

\vspace{5mm}

Dennis Wingender$^{1}$, Daniel Balzani$^{1\star}$

\vspace{3mm}

{\small $^1$Chair of Continuum Mechanics, Ruhr-Universit\"at Bochum,\\ 
Universit\"atsstra{\ss}e~150, 44801~Bochum, Germany}\\[3mm]

\vspace{3mm}

{\small ${}^{\star}$E-mail address of corresponding author: daniel.balzani@ruhr-uni-bochum.de}

\vspace{10mm}

\begin{minipage}{15.0cm}
\textbf{Abstract}\hspace{3mm}
This paper presents an algorithm for the efficient simulation of ductile crack propagation through heterogeneous structures, as e.g. metallic microstructures, which are given as voxel data. 
These kinds of simulations are required for e.g., the numerical investigation of wear mechanisms at small length scales, which is still a challenging task in engineering. 
The basic idea of the proposed algorithm is to combine the advantages of the Finite Cell Method allowing for a convenient integration of heterogeneous finite element problems with the eigenerosion approach to still enable the mesh-independent simulation of crack propagation. 
The major component is to switch from finite subcells to finite elements wherever the crack progresses, thereby automatically adaptively refining at the crack tip by managing the newly appearing nodes as hanging nodes. 
Technically relevant problems of crack propagation at the microscale are mostly linked with sub-critical crack growth where the crack moves fast and stepwise with subsequent load cycles. 
Therefore, inertia may become important which is why dynamics are taken into account by spreading the mass of the eroded elements to the nodes to avoid a loss in mass resulting from the erosion procedure. 
Furthermore, a certain treatment for the finite cell decomposition is considered in order to ensure efficiency and accuracy. 
The numerical framework as well as the voxel decomposition techniques are analyzed in detail in different three-dimensional numerical examples to show the performance of the proposed approach. 

\end{minipage}
\end{center}

\medskip{}

\medskip

\section{Introduction}
In many mechanical engineering applications, the analysis of cracks propagating through heterogeneous structures is subject to investigation. 
One prominent example is wear of hardened tool materials, which undergo cyclic loading leading to crack propagation on the microscopic level resulting in wear in form of surface spalling and even structural failure.
For instance in the context of mechanized tunneling, mining tools connected to the cutting wheel penetrate the soil and are thus subjected to high, cyclically applied forces. 
These mining tools appear in form of cutting disks and chisels consisting of ductile steel armored by wear-resistant weldings made out of hard metals or metal matrix composites (MMC). 
In their microstructures, the brittle inclusions, which supply a high hardness of the composite, are surrounded by a ductile metal matrix, leading to a high crack resistance and thus a high resistance regarding wear in different forms.
Firstly, overload breakage and abrasion appear on the macroscale.
Secondly, surface spalling governed by sub-critical crack propagation on the microscale under cyclic loading occurs. 
In order to investigate the latter case, an efficient, robust and mesh-independent simulation framework for calculations including brittle as well as ductile crack propagation at finite strains is required.
Additionally for microstructure morphologies obtained from micro-CT scans, a special numerical treatment directly taking advantage from the given voxel data is helpful to decrease computational costs. \\
The simulation of ductile crack propagation on the microscale of metallic structures requires computational methods capable of handling geometrical nonlinearities, finite strain elasto-plasticity and crack propagation along arbitrary crack paths through complex three-dimensional structures. 
Additionally, computational efficiency is desired.
Altogether, these requirements still pose a challenging task in computational engineering. 
In \citet{ShaNavMunBerBou:2019:Cmf} multiple approaches for simulating crack propagation on the microscale based on the Finite Element Method (FEM) for tackling these problems are examined and compared. 
As a continuous approach, the phase field introduced in \citet{MieWelHof:2010:ijnm} assumes that the sharp cracks are smoothed by a continuous damage field. 
This method is applied on microstructures, for instance in \citet{NguYvoZhuBor:2016:Apf} considering small deformations and brittle material behavior, in \citet{ShaGho:2019:Cpf} for elastic polycrystals and in \citet{CheTuGho:2020:Wea} assuming crystal plasticity.
However, an additional degree of freedom, namely the damage variable, occurs additionally to the displacements in the phase field method, which increases the computational costs. 
As another continuous approach, a gradient-enhanced damage model as for instance introduced in~\citet{JunSchJanHac:2019:afa} for small strains and in~\citet{JunRieBal:2021:ear} extended to finite strains, also assumes a smooth damage field for the sharp crack.
Furthermore, cohesive-zone models as introduced in \citet{Bar:1962:Tmt} for brittle and \citet{Dug:1960:Yos} for ductile materials which were also applied in an FE framework in e.g., \citet{HilModPet:1976:Aoc} make use of interface elements underlying decohesion due to cracking. 
Examples of this technique for the debonding of the inclusion from the matrix is shown in \citet{LiaSof:2003:Man,MenWan:2015:Poi}. 
These methods may suffer from erroneous crack patterns as shown in \citet{SchDeB:1993:Otn}. 
Additionally, if the crack path is unknown, the interface elements have to be applied between all elements as e.g. in \citet{XuNee:1994:Nso} leading to an increased number of elements and thus computational costs. 
On the other hand, discontinuous approaches, allowing jumps in the mechanical fields for example in the displacement fields, are applied as well to simulate the crack propagation on microscopic level. 
For example, the Extended FEM (XFEM) from \cite{BelBla:1999:ecg} handles those by application of enriched shape functions. 
For instance \citet{SukSroBakPre:2003:Bfi} shows microscopic simulations with brittle crack propagation and \citet{BeeLoeWri:2018:3dc} ductile crack propagation at finite strains. 
Other approaches, for example early element erosion techniques applied at microscale as e.g. in \citet{WulFis:1996:Fso}, suffer from mesh-dependency.\\
Therefore, a new efficient method for the simulation of ductile crack propagation through heterogeneous, metallic structures is presented in this work. 
To this end, the eigenerosion approach for ductile crack propagation as presented in \citet{WinBal:2022:SoC} is combined with the Finite Cell Method (FCM) introduced in \citet{ParDueRan:2007:fcm}.
The considered eigenerosion approach is based on the original framework presented in \cite{PanOrt:2012:aea}, which extended the idea of eigenfracture~\cite{SchFraOrt:2009:eae}. 
It has been shown to enable robust, efficient and mesh-independent simulations. 
The basic idea behind the eigenerosion strategy is to erode finite elements whenever a regularized Griffith-type criterion formulated in terms of the energy release rate is fulfilled. 
In this case the element is able to undergo {\itshape eigen}deformations for which no additional external mechanical work is needed which motivates the method's name. 
\cite{SchFraOrt:2009:eae} proved the $\Gamma$-convergence of the associated regularized energy-dissipation functional to the unregularized one as the neighborhood radius~$\epsilon$ goes to zero. 
This enables the desired mesh-independent crack propagation. 
It was firstly implemented in \cite{PanOrt:2012:aea} for brittle crack propagation at small strains. 
Extensions of it, for example for high impact loading and fragmentation, are presented in \cite{PanLiOrt:2013:ijf,LiPanOrt:2015:mom,NavRenYuLiRui:2018:mtd}. 
In~\citet{QinPanKal:2019:vef}, the eigenerosion has been firstly extended to ductile crack propagation for calculations on concrete by the application of small strain Drucker-Prager elasto-plasticity. 
In our work, the implementation of~\citet{WinBal:2022:SoC} considering finite strain $J_2$-elasto-viscoplasticity is used, where mesh-independent simulations have been numerically shown.
It is combined with the FCM, which extends the basic FEM. 
Similarly, this combination with the phase-field instead of the eigenerosion have been exploited. 
For example, a framework deriving crack initiation is given in~\citet{RanMasParDusRan:2014:Utf}. 
Furthermore,~\citet{NagElhAmbKolDeLRan:2019:Pfm} applied hp-refinement within the FCM for brittle crack propagation at small strains. 
In the FCM, the elements are allowed to contain multiple subdomains with different material properties. 
This enables the use of meshes which do not conform with the material boundaries. 
Hence, structures based on voxel data from micro-CT scans can be discretized with a regular hexahedral mesh containing summarized voxels as subdomains. 
This enables the automated discretization of varying material inhomogeneities which makes it quite advantageous in problems where numerous different material realizations have to be evaluated, e.g., in the context of uncertainty quantification, cf.~\cite{MisBal:2019:qou}. 
As shown in \citet{YanRuesKollDueRan:2012:Aei}, it also circumvents the problem of discretizing complex structures obtained by interpolation and smoothening of the voxel-data with e.g., tetrahedral elements, where the element faces are aligned with the material interfaces. 
There, complicated and computationally costly techniques as for example shown in \citet{SchKluBarg:2016:atd} have to be applied which often precludes an automatized discretization. 
The algorithm proposed here, is mainly based on switching the subcells to finite elements in those cells where the crack evolves. 
This leads to an automatic mesh refinement at the crack tip enabling a high accuracy and computational efficiency at the same time. 
For the newly appearing hanging nodes at the interfaces between the refined cells and the neighboring elements, the concept of Lagrange multipliers in \citet{DemOdeRacHar:1989:tau,OdeDemRacWes:1989:tau,RacOdeDem:1989:tau} is applied. 
Different algorithms for decomposing the elements into subcells based on voxel data are proposed and compared in~\citet{FanMisBal:2020:aso}. 
For the combination with the eigenerosion, additional restrictions have to be considered regarding the decomposition because of the change from subcells to finite elements. 
Whereas the aspect ratio of the individual subcell's dimensions in different directions can be rather arbitrary, it is bounded in case of finite elements. \\ 
In this paper, a new algorithm for the voxel-based simulation of crack propagation through metallic microstructures based on a combination of the eigenerosion approach extended to finite strain elasto-plasticity and the FCM is presented. 
The paper is organized as follows. 
In section 2, the basic concept of the eigenerosion approach for ductile crack propagation at finite strains is recapitulated as it represents a key component. 
Afterwards, the new algorithm combining the FCM with eigenerosion is presented in section 3. 
Furthermore, different subcell decomposition techniques for voxel-based microstructure data and their implications to the proposed algorithm are discussed. 
Section 4 presents different, three-dimensional numerical examples based on artificially generated voxel data and on real voxel data which show the performance of the framework. 
A final conclusion is given in Section 5.

\section{Eigenerosion for Ductile Fracture at Finite Strains}

As the eigenerosion framework is considered here as major component for the description of crack propagation, the concept  for ductile fracture at finite strains from \citet{WinBal:2022:SoC} extending the original ideas of \citet{SchFraOrt:2009:eae,PanOrt:2012:aea} is briefly recapitulated.
Furthermore, the specific material description considered in this paper is explained. 

\subsection{Eigenerosion Algorithm}

\begin{figure}[t]
\unitlength1cm
\begin{picture}(16.0,8.0)
\put( 0.2, 0.0){\includegraphics[width=0.55\linewidth]{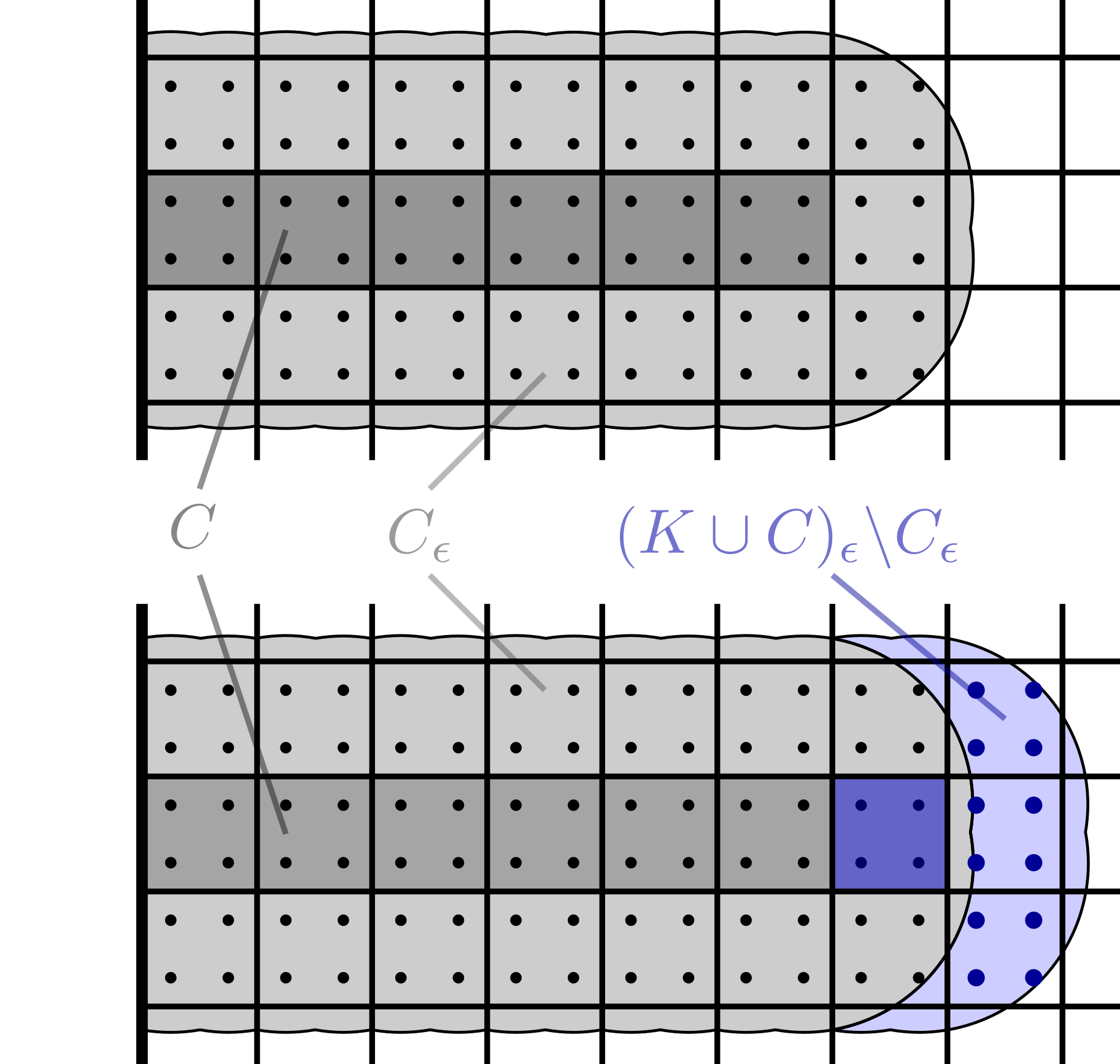}}
\put( 0.0, 0.0){(a)}
\put( 9.0, 0.0){\scalebox{0.75}{
\tikzstyle{decision} = [diamond, draw, fill=blue!40, 
    text width=2.0cm, text badly centered, node distance=2.5cm, inner sep=0pt]
\tikzstyle{block} = [rectangle, draw, fill=lightgray!60, 
    text width=2.5cm, text centered,  minimum height=1.0cm, node distance=1.8cm]
\tikzstyle{line} = [draw,ultra thick, -latex']
\tikzstyle{cloud} = [draw, ellipse,fill=red!20, node distance=2.0cm,
    minimum height=2em]  

\begin{tikzpicture}[node distance = 1.8cm, font=\small,auto]
    \node [block] (init) { initialize FE model};
    \node [block, below of=init, node distance=1.8cm] (crack_init) {compute list with Gau\ss points for each element $K$};
    \node [block, below of=crack_init] (timeupdate) {update time step $n\rightarrow n+1$};
    \node [block, below of=timeupdate, node distance=1.4cm] (solve) {solve nonlinear FE problem};
    \node [block, below of=solve,fill=blue!15, node distance=1.6cm ] (netgain) {compute $-\Delta F_K$ for each element $K$};
    \node [decision, below of=netgain, node distance=2.4cm] (decide) {any $-\Delta F_K>0$};
    \node [block, right of=decide, fill=blue!40, node distance=3.5cm] (update) {erode element with largest net energy gain $-\Delta F_K$};
    \node [block, right of=solve, fill=blue!15, node distance=3.5cm] (netupdate) {update list with Gau\ss points};
    
    \path [line] (init) -- (crack_init);
    \path [line] (crack_init) -- (timeupdate);
    \path [line] (timeupdate) -- (solve);
    \path [line] (solve) -- (netgain);
    \path [line] (netgain) -- (decide);
    \path [line] (decide) -- node  {yes} (update);
    \path [line] (netupdate) -- (solve);
    \path [line] (decide) -- ++ (-2.0cm,0) |-  node [near start] {no} (timeupdate);
    \path [line] (update) -- (netupdate);;
\end{tikzpicture}}}
\put( 9.6, 0.0){(b)}
\end{picture}
\caption{(a)~Mesh with eroded elements $C$ (dark grey), the Gau{\ss} points in their $\epsilon$-neighborhood $C_{\epsilon}$ (light grey) and additional crack area $\left( C \cup K\right)_{\epsilon} \backslash C_{\epsilon}$ (light blue) for element domain $K$. (b)~Schematic illustration of the considered eigenerosion algorithm within a finite element framework.}
\label{fig:eigenerosion}
\end{figure}
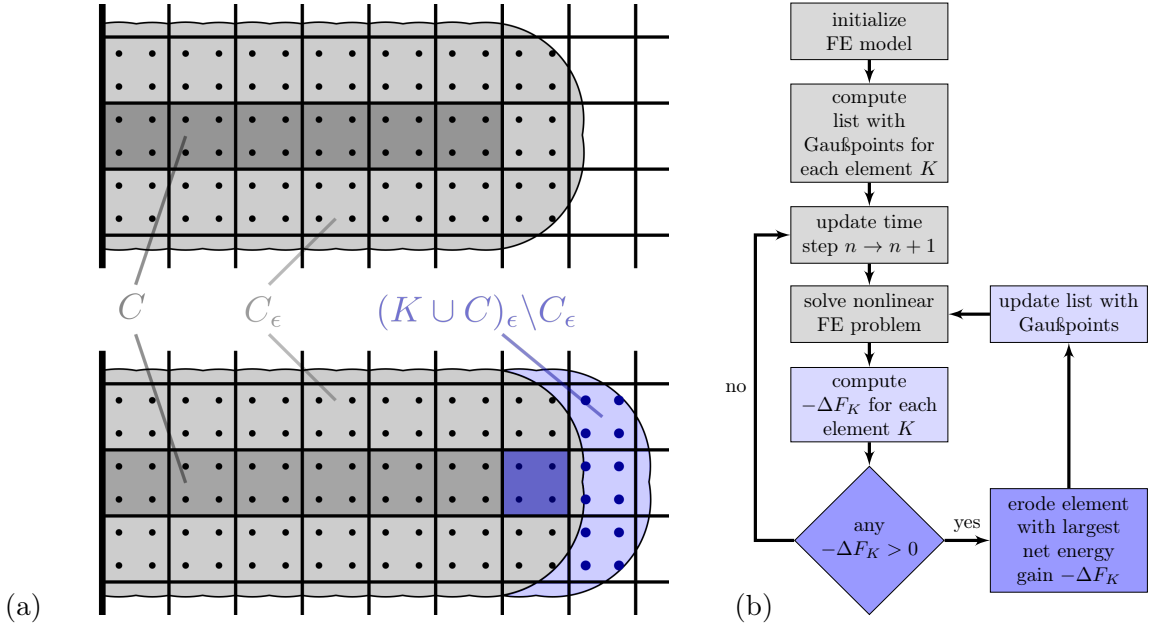

The eigenerosion approach is implemented in connection with the finite cell method, which is in turn based on the standard FE framework for solid mechanics.
Because of that, the basic equations of the FEM are recapitulated briefly to define notation.
Therein, the displacement~$\bu$ which minimizes the total potential energy~$\Pi(\grad\bu,\bu)$ is saught under given boundary constraints.
To this end, standard terms of variational calculus require the 1st variation of $\Pi$ to vanish and thus, the weak form reads
\eb
\delta\Pi=\int\limits_{\Omega} \nabla^\text{S}_x\delta\bu:\Btau\, \mathrm{d}V
- \int\limits_{\partial \Omega^{\mathrm{N}}} \delta\bu\cdot \bt \, \mathrm{d}A
+\int\limits_{\Omega} \rho_0\delta\bu\cdot\ddot{\bu} \, \mathrm{d}V = 0
\ee
including inertia and neglecting body forces is sought.
Herein, the Kirchhoff stress tensor~$\Btau$ is double-contracted with the symmetric part of the spatial gradient of displacement variations~$\nabla^{\text{S}}_x\delta\bu=\half(\grad( \delta \bu) + \grad(\delta \bu)^{\mathrm{T}})$ and integrated over the volume~$\Omega$ in the reference configuration.
The external traction forces acting on the Neumann surface~$\delta \Omega^{\mathrm{N}}$ and the density in the reference configuration are denoted by~$\bt$ and~$\rho_0$, respectively.
Applying the FEM, this nonlinear partial differential equation can be solved by spatial decomposition of the body~$\Omega$ into finite elements~$K$ with their domains~$\Omega_K$.
In the individual elements a standard polynomial approximation for the displacements and the displacement variations is introduced as $\bu \approx \bN\bd_K$ and $\delta\bu \approx \bN \delta\bd_K$, and following therefrom one obtains $\nabla^{\text{S}}_x \delta\bu \approx \bB\delta\bd_K$.
Herein, the vectors~$\bd_K$ and $\delta\bd_K$ contain all nodal displacements and nodal displacement variations per element~$K$.
Furthermore, matrix notation has been used together with standard matrices~$\bN$ and $\bB$ including the ansatz functions and their spatial derivatives, respectively.
Linearizing the resulting approximated version of the weak form in terms of the Newton-Raphson scheme, and applying the Newmark approach for the time integration of the inertia term, results in the equation
\eb
\sum_K\delta\bd_K^{\text{T}}\left[
\left(\bk_K + \bbm_K\right)\,\Delta\bd_K
- \bbr_K
- \bbr^{\text{m}}_K
+ \bq_K
\right] = \bzero
\label{eq:eqsystem}
\ee
with the linear increment of the displacements~$\Delta\bd_K$.
Herein, the element tangent stiffness matrix~$\bk_K:=\int_{\Omega_K}\bB^{\text{T}}\mathbb{C}\bB\,\text{d}V$ with the material tangent modulus matrix~$\mathbb{C}$ in Voigt notation, the internal element residual vector~$\bbr_K:=-\int_{\Omega_{K}} \bB^{\text{T}}\Btau\text{d}V$, and the element vector of external forces~$\bq_K:= -\int_{\partial \Omega^{\text{N}}_K}  \bN^{\text{T}} \bt \, \text{d}A$.
The mass matrices $\bbm_K$ and inertia residual vectors~$\bbr^{\text{m}}_K$ are chosen as the consistent mass matrix and residual vector resulting from the Newmark scheme where no erosion takes place, and the lumped mass matrix and residual vector are considered otherwise.
Note that in principle any other time integration scheme may be used instead of the Newmark method.
The equation~\eqref{eq:eqsystem} is rewritten in terms of the classical global matrices, the vectors of unified nodal displacement variations and increments $\delta\bD=\Unifargone{K}\delta\bd_K$ and $\Delta\bD=\Unifargone{K}\Delta\bd_K$,
the tangent stiffness matrix $\bK:=\Assemargone{K}[\bk_K]$,
the mass matrix~$\bM:=\Assemargone{K}[\bbm^{\text{e}}_K]$,
the internal residual vector $\bR:=\Assemargone{K}[\bbr_K]$,
the inertia residual vector~$\bR^{\text{m}}:=\Assemargone{K}[\bbr^{\text{m}}_K]$,
and the vector of external forces at the Neumann boundary~$\bQ:=\Assemargone{K}[\bq_K]$.
Together with the incorporation of Dirichlet boundary conditions the Newton-Raphson scheme yields the linearized system of equations $(\bK+\bM)\Delta\bD=\bR+\bR^{\text{m}}-\bQ$, which is solved repeatedly for the incremental displacements~$\Delta\bD$ within the Newton-Raphson iteration and updated as $\bD\Leftarrow\bD+\Delta\bD$ to compute the global vector of nodal displacements~$\bD$ which fulfills mechanical equilibrium.\\
The original eigenerosion approach is based on the theory of brittle fracture, firstly published in \citet{Gri:1921:tpo}. 
Therein, the existence of a Griffith-type energy release rate $G:=-\pp{U(\bu)}{|C|}$ is assumed. 
This rate relates the potential energy~$U$ stored by imposing mechanical work depending on the displacement vector~$\bu$ to the area $|C|$ of the crack set~$C$. 
Additionally, irreversibility of crack propagation and no healing of the material are considered. 
Hence, the monotonicity constraint $C(t)\subset C(t+\Delta t)$
%
results for the crack set $C(t)$ at time $t$ and the crack set $C(t+\Delta t)$ at a later time step $t+\Delta t$ assuming discrete time steps. 
In terms of FE discretizations, the crack set $C$ consists of elements that are assumed to be eroded, which enables them to undergo eigendeformations, for which no additional energy is required. 
It increases with the crack front velocity~$v$ if the Griffith-type energy release rate~$G$ reaches a critical material-dependent value~$G_c$ and thus, the criterion $G-G_c \le 0$ has to hold. 
If not, the crack rests so that the crack front velocity~$v$ becomes zero. 
Combining these conditions yields the expression $(G-G_c)\,v=0$ to hold. 
Together with the concept presented in \citet{MieOrt:2008:aco}, the energy dissipation functional $F=U(\bu)+G_c|C|$ is formulated. 
This functional has to be minimized at every time $t$ with respect to the crack set~$C$ and the displacement field~$\bu$. 
In this naive form, the eigenerosion would lead to mesh-dependent results. 
Thus, in order to avoid this, the energy-dissipation functional is regularized by $F_{\epsilon}=U+G_c|C_{\epsilon}|/(2\,\epsilon)$. 
Herein, the $\epsilon$-neighborhood $C_{\epsilon}$ of the crack set $C$ within the influence radius $\epsilon>0$ is considered instead of only the crack set $C$. 
Hence, the mesh independence of this algorithm can be recovered. 
\citet{SchFraOrt:2009:eae} has proven that the regularized energy-dissipation functional $F_{\epsilon}$ $\Gamma$-converges to the unregularized energy-dissipation functional $F$ as $\epsilon \rightarrow 0$. 
Transferring this concept to element erosion in finite element simulations, the net energy gain of each element~$K$ 
\begin{equation}
-\Delta F_{K}=-\Delta U_K-G_c\,\Delta A_K
\end{equation}
is evaluated for each element $K$ containing the effective crack area $\Delta A_K$ in each time step $t_n$ after solving the mechanical fields. 
Here, the difference in its potential energy $-\Delta U_K$ before and after erosion becomes the energy $U_K$ which is stored in the element~$K$. 
If $-\Delta F_K$ becomes larger than zero, crack propagation occurs, otherwise the crack rests. 
The regularized crack area 
\begin{equation}
\Delta A_K =\frac{|\left( C \cup K\right)_{\epsilon} \backslash C_{\epsilon}|}{2\,\epsilon}
\end{equation}
represents the volume of the $\epsilon$-neighborhood $|C\cup K|_{\epsilon}$ of the whole crack $C$ including the one of element $K$ without the one of the previous crack $C_{\epsilon}$, cf. figure \ref{fig:eigenerosion}a, divided by the influence radius $\epsilon$. 
This represents the regularized crack area that would additionally occur if element $K$ was eroded. 
In our case, this term is evaluated based on the points of the Gau\ss\, point quadrature and their volumes that are also used for deriving the element residual and tangent stiffness.
If the distance of any Gau\ss\, point to a Gau\ss\, point of element $K$ is lower than the influence radius $\epsilon$, it is assumed to be part of the $\epsilon$-neighborhood of element $K$ and it is saved in a list. 
Gau\ss\, points lying in the $\epsilon$-neighborhood of eroded elements are removed from this list. 
For the evaluation of $\Delta A_K$, the volumes of all Gau\ss points in the list of element $K$ are summed up and normalized by dividing by $2\,\epsilon$. 
Furthermore, the difference in the potential energy~$-\Delta U_K$ is evaluated for each element~$K$ for the net energy gain~$\Delta F_K$. 
Assuming small strains and linear elasticity, it can be calculated in terms of the element stiffness matrix~$\bk_K$ and nodal degrees of freedom $\bd_K$ of element~$K$ by $0.5\,\bd_K^{\mathrm{T}}\,\bk_K\bd_K$, as proposed in~\citet{PanOrt:2012:aea}. 
However, here we focus on more general cases of nonlinear material formulations at finite strains including plasticity. 
Following~\cite{WinBal:2022:SoC}, the total strain energy density~$\psi$ integrated over the element's domain~$\Omega_K$ can be considered, i.e. 
\eb
-\Delta U_K
=\int\limits_{\Omega_K}\left(\psi\,+\int_t {\cal D}^{\mathrm{vis}} \,\mathrm{d}t\right)\mathrm{d}V
\quad\mbox{with}\quad
\psi = \psi^{\text{e}}(\Bvarepsilon^{\text{e}}) + \psi^{\text{p}}(\alpha)
.
\ee
%
Herein, the elastic part of the strain energy density~$\psi^{\text{e}}$ depends on the elastic logarithmic strain tensor~$\Bvarepsilon^{\text{e}}$ and the plastic part~$\psi^{\text{p}}$ depends on a variable~$\alpha$ associated with microstructure changes due to hardening. 
In~\cite{WinBal:2022:SoC} it has been shown that regions of localized plastic strains such as shear bands induce issues with mesh-dependence already before a crack evolves, and thus, a viscous part is additionally considered. 
The associated viscous dissipation~$\int_t\mathcal{D}^{\text{vis}}\text{d}t$ is therefore added to the difference of the potential energy. 
This approach coincides with the extension of \citet{Irw:1957:aos} additionally considering the energy that dissipates in the crack tip due to plasticity since it is part of the imposed energy as well. 
The individual energetic terms~$\psi^{\text{e}}$, $\psi^{\text{p}}$ and~${\cal D}^{\text{vis}}$ depend on the chosen material law which will be specified later.\\
The resulting eigenerosion algorithm is illustrated in figure~\ref{fig:eigenerosion}b. 
Therein, the mechanical equilibrium equations are firstly solved for the displacements~$\bu$ in every time step while not assuming any crack propagation.
Thus, this step represents a trial step. 
Afterwards, the net energy gain $-\Delta F_K$ is evaluated.
If the net energy gain of any element $K$ becomes larger than zero, the element with the largest net energy gain is eroded, which means that its residual and static tangent stiffness matrices are set to zero and the Gau\ss\, point list is updated. 
The algorithm considers inertia effects by applying the Newmark-scheme~\cite{New:1959:amc}. 
This recovers numerical stability even if the body completely breaks in multiple parts. 
Furthermore, the incorporation of inertia may allow for more accurate simulations in case of fast crack evolution. 
In order to not loose mass as a result from the erosion procedure and still enable the detached nodes to move independently, the lumped mass matrix only containing entries on the diagonal is considered for the eroded elements. 
For all intact elements, the consistent mass matrix is used. 
Then, all previous steps are repeated in a staggered scheme until no element erodes anymore in this time step. 
After that, the algorithm proceeds with the next time step.

\subsection{Material Model for Ductile Fracture
\label{sec:matmod}}
For the simulation of ductile fracture, the finite $J_2$-elasto-plasticity material model with exponential isotropic hardening as introduced in \cite{Sim:1992:afs, SimMie:1992:cmam, SimTay:1984:cmam} is applied. 
The numerical implementation is given in \citet{MieSteWag:1994:ame} and \citet{Kli:2000:tun} and its modification for the eigenerosion approach in \citet{WinBal:2022:SoC}. 
As demonstrated in the latter reference, the elasto-plastic material model leads to numerical localization of the plastic fields and therefore to mesh dependence of the framework of ductile fracture. 
To avoid this, viscous regularization is applied.
Furthermore, the viscosity, here the Perzyna-type, cf. e.g., \citet{Per:1966:Fpi,JunSchMakHac:2017:cmt}, leads to numerically more stable calculations especially if large plastic deformations occur, because its internal variables are derived explicitly. 
The material models considered in the numerical analyses presented later are briefly recapitulated in this section. 
Those are based on the multiplicative decomposition of the deformation gradient $\bF=\bF^{\text{e}} \cdot \bF^{\text{p}}$
into an elastic part~$\bF^{\text{e}}$ and a plastic part~$\bF^{\text{p}}$ as e.g. shown in~\citet{Kro:1959:akd,Lee:1969:epd}. 
Based on this, the elastic Cauchy-Green tensor $\bb^{\mathrm{e}}=\bF^{\mathrm{e}} \cdot \bF^{\mathrm{eT}}=\textstyle\sum_{i=1}^3(\lambda_i^{\mathrm{e}})^2\,\bn_i \otimes \bn_i$
allows the spectral decomposition into the eigenvalues resulting into the principal elastic stretches~$\lambda_i^{\mathrm{e}}$ and their corresponding eigenvectors~$\bn_i$ in form of the principal stretch directions. 
Following this, the principal logarithmic elastic strains $\epsilon_i^\mathrm{e}=\log(\lambda_i^\mathrm{e})$ and their corresponding eigenvectors $\bn_i$ are composed to the elastic strain tensor 
\eb
\boldsymbol{\varepsilon}^{\mathrm{e}}=\sum\limits_{i=1}^3 \epsilon_i^{\mathrm{e}}\,\bn_i \otimes \bn_i
\ee
enabling the additive split of the strain tensor $\boldsymbol{\varepsilon}=\boldsymbol{\varepsilon}^{\mathrm{e}}+\boldsymbol{\varepsilon}^{\mathrm{p}}$ into an elastic and plastic part $\boldsymbol{\varepsilon}^{\mathrm{e}}$ and $\boldsymbol{\varepsilon}^{\mathrm{p}}$, respectively.
Based thereon, a quadratic elastic part of the strain energy density is considered as 
\eb
\psi^\mathrm{e}=\frac{\kappa}{2}\, \mathrm{tr}(\boldsymbol{\varepsilon}^{\mathrm{e}})^2+\mu\, \mathrm{dev}(\boldsymbol{\varepsilon}^{\mathrm{e}}) : \mathrm{dev}(\boldsymbol{\varepsilon}^{\mathrm{e}})
\ee
containing the compression modulus $\kappa$, the shear modulus $\mu$ and the deviator operation $\mathrm{dev}(\boldsymbol{\varepsilon}) := \boldsymbol{\varepsilon} - 1/3\, \mathrm{tr}(\boldsymbol{\varepsilon})\bI$ with the identity tensor $\bI$. 
Furthermore, the convex plastic dissipation is assumed as the superposition of linear hardening controlled by its slope~$h^{\mathrm{lin}}$ and exponential hardening with the degree of exponential hardening $h^{\mathrm{exp}}$, i.e. 
\eb
\psi^{\mathrm{p}}=y_0\,\alpha +  (y_{\infty}-y_0) \,\left[\alpha+\frac{\exp(-h^{\mathrm{exp}}\,\alpha)-1}{ h^{\mathrm{exp}}}\right]+\frac{1}{2}\,h^{\mathrm{lin}}\,\alpha^2
\ee
as used in~\citet{Voc:1955:aps}. 
The parameter $y_0$ describes the initial yield stress and $y_{\infty}$ the stress where the exponential hardening approaches almost purely linear hardening. 
By derivation of the two parts of the strain energy density function, the Kirchhoff stress tensor 
\eb
\boldsymbol{\tau}=\pp{\psi^\mathrm{e}(\boldsymbol{\varepsilon}^{\mathrm{e}})}{\boldsymbol{\varepsilon}^{\mathrm{e}}}=\kappa\,\mathrm{tr}(\boldsymbol{\varepsilon}^\mathrm{e})\,\bI+2\,\mu\, \mathrm{dev}(\boldsymbol{\varepsilon}^\mathrm{e}).
\ee
and the hardening function
\eb
\beta=\frac{\partial \, \psi^{\mathrm{p}}(\alpha)}{\partial \, \alpha}=y_0+(y_{\infty}-y_0)\,\left[1-\exp(-h^{\mathrm{exp}}\,\alpha)\right]+h^{\mathrm{lin}}\,\alpha.
\ee
follow.
The von Mises type flow condition $\phi=\tau^{\mathrm{vM}}-\beta \leq 0$ 
is applied with the von Mises stress $\tau^{\mathrm{vM}}=\sqrt{3/2}||\mathrm{dev}\,\Btau||$ assuming plastic incompressibility, which is typical for metal plasticity. 
The equation $\phi = 0$ is solved for the plastic parameter $\lambda^{\mathrm{p}}$ describing the plastic
evolution of the norm of plastic strains $||\dot{\boldsymbol{\varepsilon}}^{\mathrm{p}}||=\lambda^{\mathrm{p}}$
and with respect to the evolution of the equivalent plastic strains $\dot{\alpha}=\sqrt{\frac{2}{3}}\,\lambda^{\mathrm{p}}$ 
%
using a local Newton iteration. 
The plastic parameter is obtained by 
\eb
\lambda^{\mathrm{p}}=\frac{1}{\eta}\,\left\langle \phi\right\rangle_{+}
\label{eq:evo_evp}
\ee
with the Macaulay bracket $\left\langle (\bullet) \right\rangle_{\pm}=((\bullet) \pm |(\bullet)|)/2$ and the viscosity $\eta$. 
Note that for the description of crack propagation through brittle materials, simply the plastic part of the strain energy density as well as the viscous dissipation are set to zero and the elastic strains become the total strains.

\section{Proposed Combined Algorithm}
%
\begin{figure}[t]
\begin{picture}(26.60,170)
\put(-75 ,-50){
\includegraphics[width=0.82\textwidth]{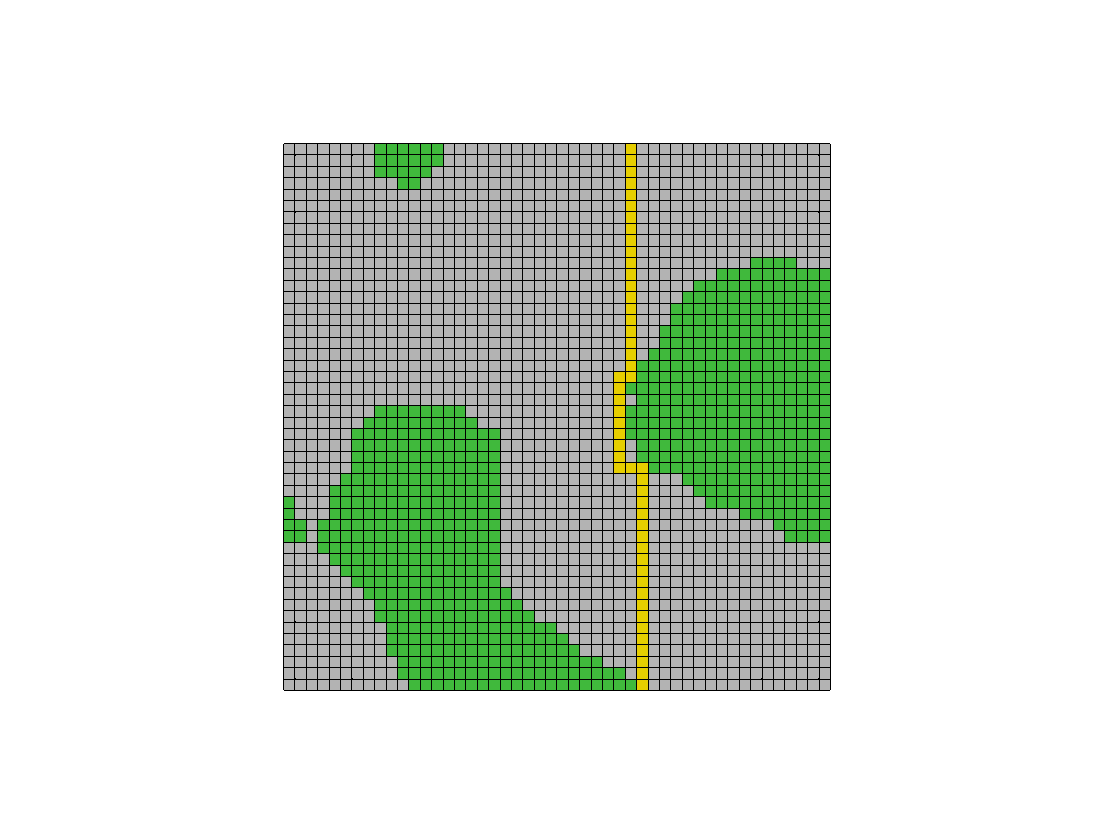}}
\put(190 ,-18){\includegraphics[width=0.64\textwidth]{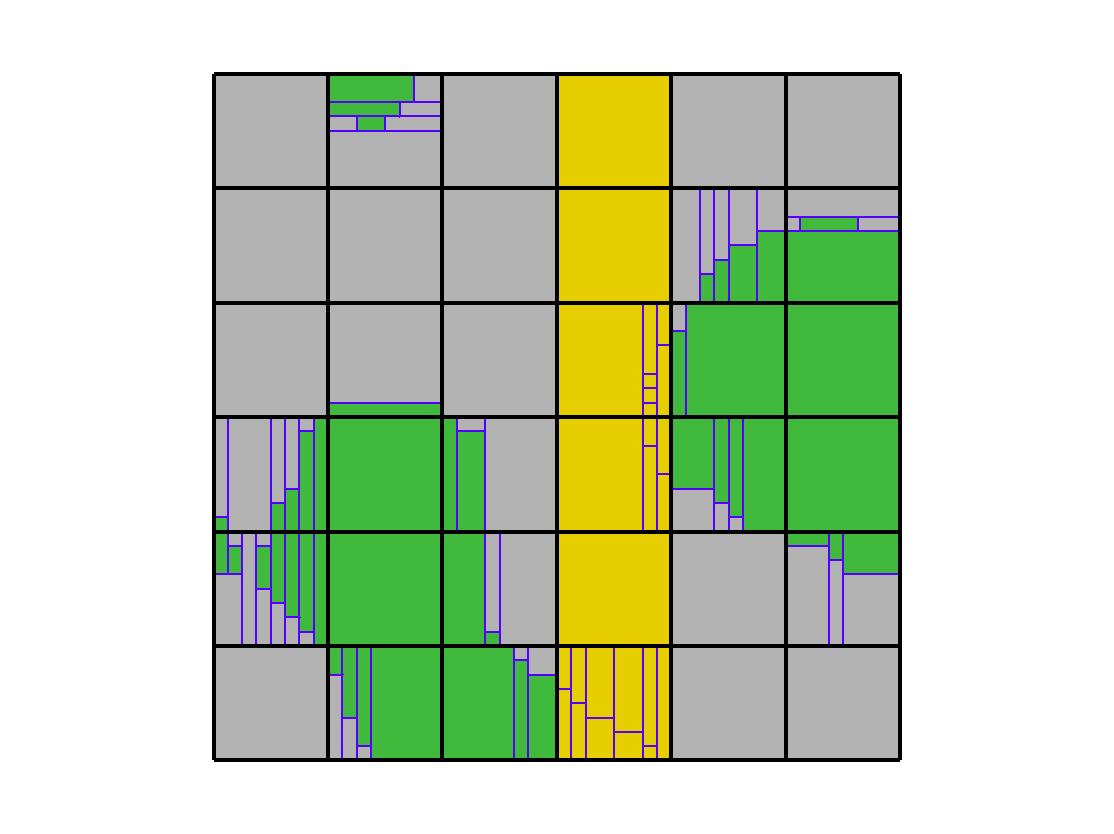}}
\put(  0 ,0){(a)}
\put(220 ,0){(b)}
\end{picture}
\caption{Schematic illustration of a possible crack path (yellow) through a heterogeneous structure of $48\times 48$ voxels naively using the eigenerosion scheme: (a) discretization of each voxel with one finite element and (b) decomposition of the structure into $6\times6$ finite subcells using the optimal decomposition approach from \cite{FanMisBal:2020:aso}.
\label{fig:decomposition_motivation}}
\end{figure}
Here, the focus is on the analysis of heterogeneous structures given as voxel data. 
In principle, conforming, irregular finite element discretizations can be constructed to capture a smoothened interface morphology which has been previously generated from the voxel data. 
However, these discretizations are often difficult to construct for complex morphologies and an automated meshing is quite difficult if at all possible, cf.~\citet{SchKluBarg:2016:atd}. 
Furthermore, such discretizations usually require a large number of finite elements even for coarse representations. 
Additionally, ill-shaped elements might occur, which lead to inaccuracies in the evaluation of the volume integral in the elements and numerical instabilities. 
Therefore, in this paper we are interested in numerical schemes which do not require such complicated, conforming discretizations. \\
Crack propagation through heterogeneous structures which are given as voxel data can then be directly described using the eigenerosion approach by considering one finite element for each voxel. 
An example is demonstrated in figure~\ref{fig:decomposition_motivation}a where a possible crack through a heterogeneous structure of $48\times48$ voxels/elements is depicted. 
Of course this straightforward approach is not efficient since a large number of finite elements is necessitated by the fact that for conforming discretizations the element faces need to coincide with the material interfaces. 
Thus, severe computational effort in terms of computing time and memory is required. \\
An alternative is to directly summarize voxels to larger elements. 
Considering this, simulations with a semi-regular, conforming mesh including hanging nodes reduce the computational costs in the assembling process but still lead to a high number of FE equations which have to be solved, especially due to many additional constraint equations which appear everywhere at the material interfaces. \\
A more efficient alternative is the Finite Cell Method (FCM) \cite{ParDueRan:2007:fcm} which enables the consideration of material interfaces inside the elements by exchanging the finite element by a finite cell, where the integration is split into subcells. 
Thereby, a significantly reduced number of finite cells is achieved, especially if higher order polynomials are used for the approximation of the displacements. 
However, the naive application of eigenerosion in combination with the FCM leads to unsuitably large finite cells which are to be eroded, if the net energy gain is compared at the finite cell (not subcell) level. 
Figure \ref{fig:decomposition_motivation}b demonstrates the scenario for a possible crack through a heterogeneous structure discretized with $6\times 6$ finite cells. 
As can be seen, especially close to the material interfaces where a high resolution would be required, the crack would only be captured with poor accuracy. 
Even if the erosion procedure was transformed such that individual subcells would be eroded instead of complete cells, the erosion of single subcells would in fact result in large deformations of the complete cell due to the loss of stiffness resulting from the eroded subcells. 
This would in turn lead to subsequent erosion of all subcells within the particular finite cell rendering the procedure inaccurate again. 

Therefore, we propose to switch from finite subcells to finite elements of relatively low polynomial order wherever subcells would be eroded and then proceed with the regular (extended) eigenerosion approach. 
Thereby, an automatic adaptive mesh refinement is realized at the crack tip whenever the finite cells are decomposed into subcells. 
This will always be the case at the material interfaces, which reflects the necessity of an increased resolution of the crack at the material interfaces. 
We consider the subcells to be converted to finite elements which use the same integration order as the subcells and thus, all material history is already known at the according Gau{\ss} points, which avoids expensive projection computations. 
Conceptually, the numerical procedure in each time step is thus divided into the following steps 

\begin{itemize}
 \item[1.] Solve mechanical equilibrium equations using the Newton-Raphson scheme to obtain trial state
 \item[2.] Compute net energy gain $-\Delta F^{\text{trial}}$ for every subcell in every finite cell
 \item[3.] If $-\Delta F^{\text{trial}} > 0$ anywhere, identify the particular finite cell $K$ where the particular subcell $S$ with the largest value of $-\Delta F^{\text{trial}}$ appears
 \item[4.] Switch all finite subcells contained in finite cell $K$ to finite elements and erode the element which has been transformed from subcell $S$
 \item[5.] Repeat steps 1-4 until no finite cell contains a subcell where $-\Delta F^{\text{trial}} > 0$
\end{itemize}
Further details regarding the FCM, the switch from finite subcells to finite elements, and its algorithmic implementation are given in the following subsections. 

\subsection{Finite Cell Method}
\begin{figure}[t]
\begin{picture}(26.60,150.0)
\put(-50 ,00){\includegraphics[width=0.55\textwidth]{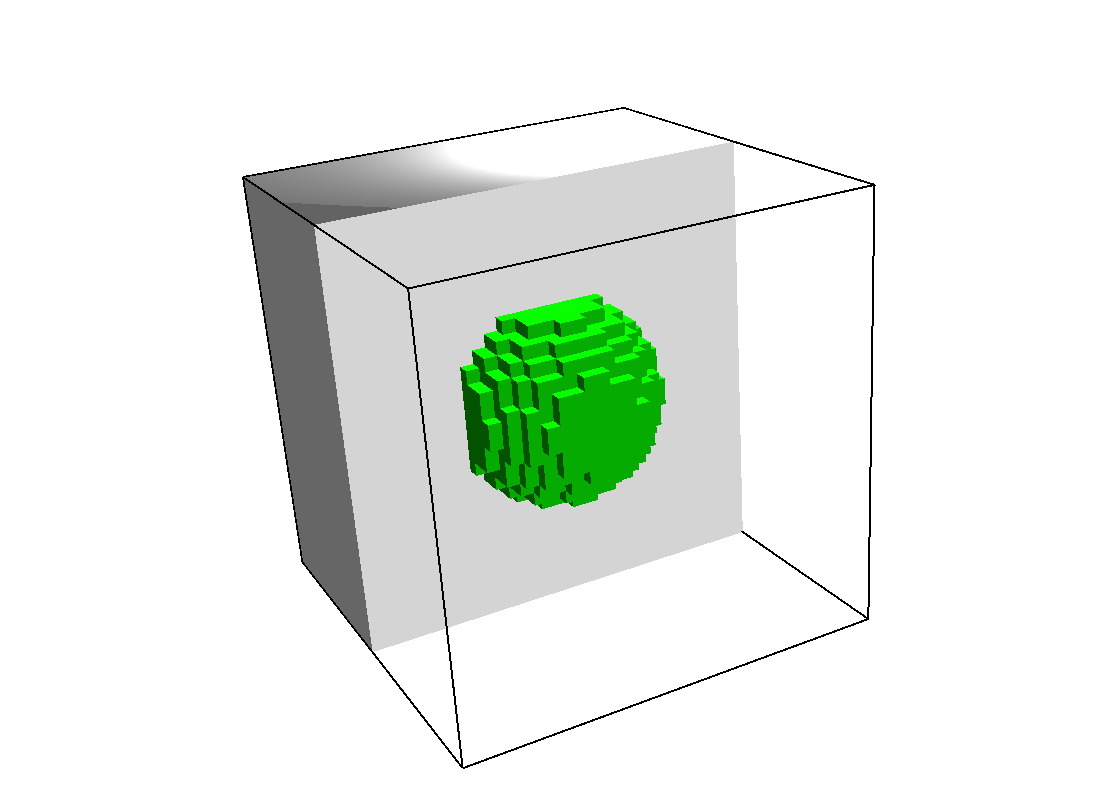}}
\put(0 , 00){(a)}
\put(100 ,00){\includegraphics[width=0.55\textwidth]{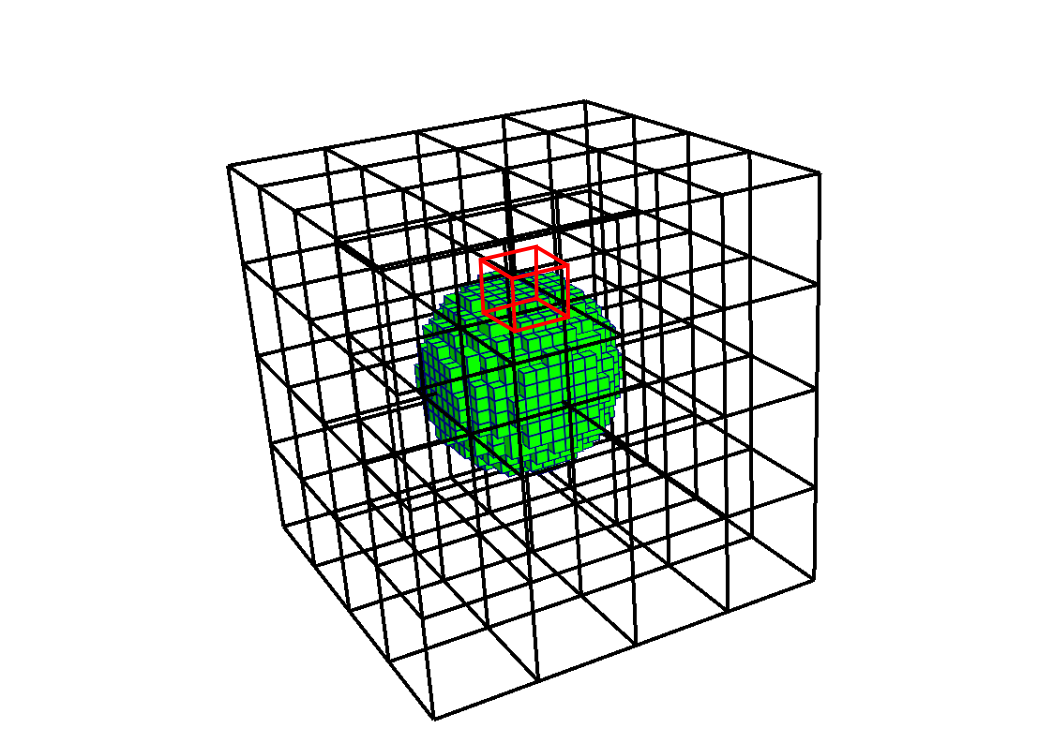}}
\put(155 ,00){(b)}
\put(260 ,10){\includegraphics[width=0.5\textwidth]{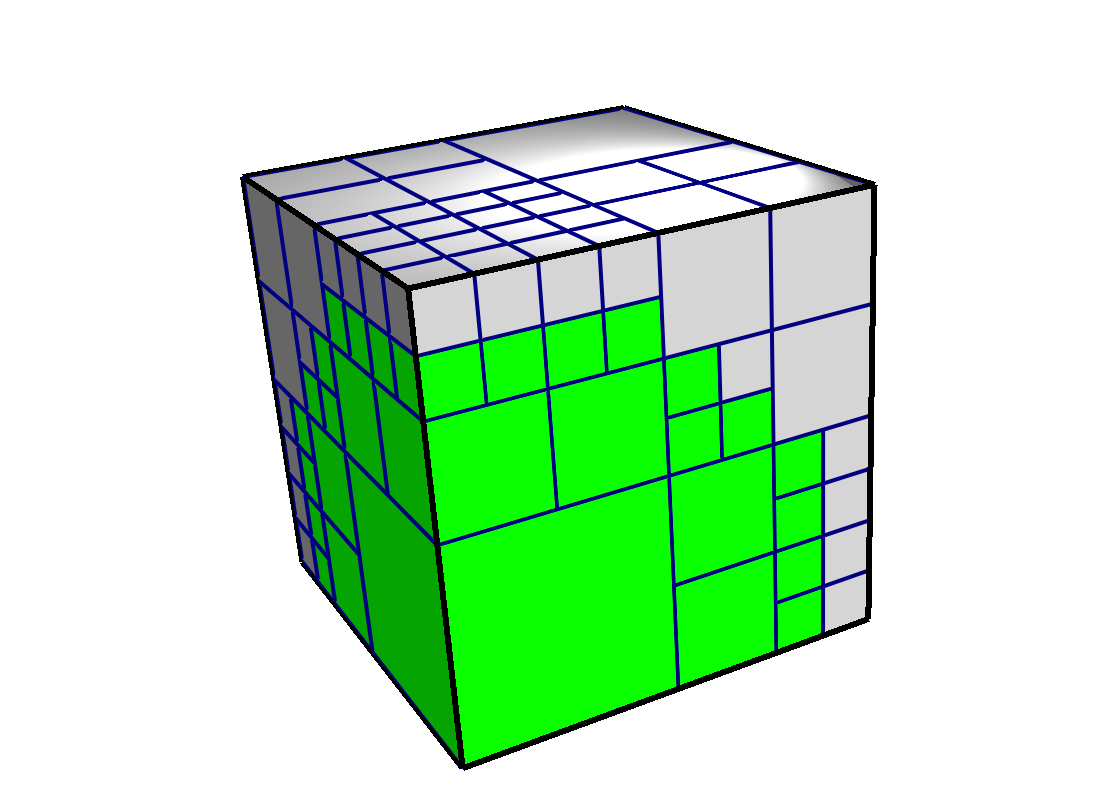}}
\put(315 ,00){(c)}
\put(230 ,112){\tikz\draw[red,line width=1] (0,0)--(2.8,0.75);}
\put(230 ,15){\tikz\draw[red,line width=1] (0,0)--(4.3,-3.0);}
\end{picture}
\caption{(a) Biphasic body $\mathcal{B}$ consisting of spherical inclusion (green) in cubic matrix (grey) given as voxel data, (b) structured, regular finite cell mesh and (c) one amplified finite cell consisting of multiple subcells with different material properties.}
\label{fig:FCM_concept}
\end{figure}
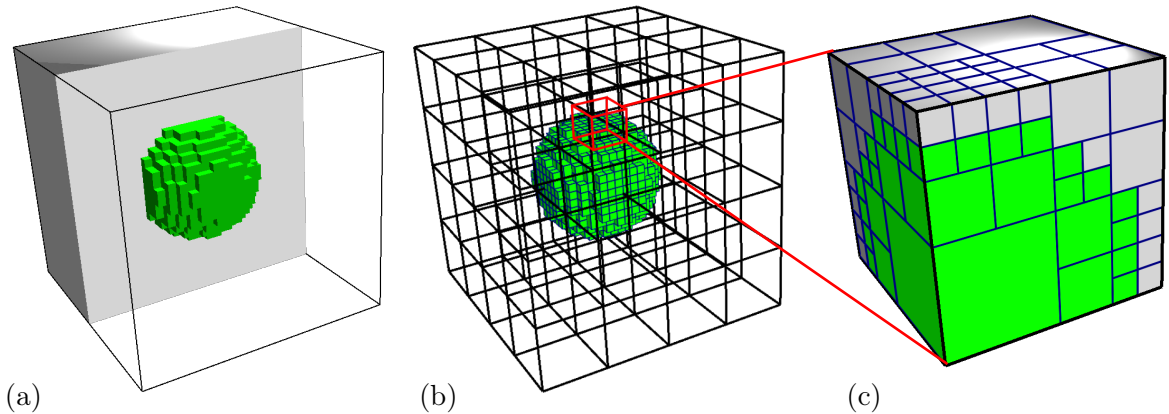

Since the FCM as introduced in \cite{ParDueRan:2007:fcm} is major component of the proposed algorithm, it is briefly recapitulated in this section and its adaption to its use within the proposed algorithm is explained. 
In the FCM the finite elements are replaced by finite cells which allow the decomposition of the cell domain $\Omega^{\text{fc}}$ into $n^\text{sc}$ subdomains $\Omega_1^{\text{sc}}\cup \Omega_2^{\text{sc}} \cup \cdots\cup \Omega_{n^\text{sc}}^{\text{sc}}=\Omega^{\text{fc}}$ referred to as subcells. 
The individual subcells are assumed to be homogeneous, but each subcell is allowed to have different material properties. 
The subcell decomposition is only considered whenever there are material inhomogeneities within one finite cell. 
This permits meshing of the domain of interest with larger domains of simple shape, for instance regular hexahedrals, wherein material interfaces are allowed to appear, see figure \ref{fig:FCM_concept} for an illustration. 
This leads to the extension of the tangent stiffness matrix of each finite cell $K$ in Voigt notation as 
\eb
\bk^{\text{fc}}_K
:=\sum\limits_{S=1}^{n^{\text{sc}}}\bk^{\text{sc}}_S
=\sum\limits_{S=1}^{n^\text{sc}}\, \int\limits_{\Omega_S^{\text{sc}}}\bB^{\text{T}}\mathbb{C}_S\bB\,\text{d}V_S
= \sum\limits_{S=1}^{n^\text{sc}}\sum\limits^{l_\text{int}}_{l=1}(\bB^{\text{T}}\mathbb{C}_{S}\bB)|_{\boldsymbol{\xi}_l} \det \bJ^{\text{sc}}_{S,\,l}\, w_{S,\,l}\, \det \bJ^{\text{fc}}
\ee 
representing the summation of the partial stiffness matrices $\bk^{\text{sc}}_S$ of the subcell domains $\Omega_S^{\text{sc}}$.
Herein, $\mathbb{C}_S $ is the material tangent moduli matrix in Voigt notation of the material in subcell~$S$.
In order to perform volume integration analogously to the FE framework, the material response is evaluated in $l_{\text{int}}$ Gau{\ss} points $l$ of each subcell~$S$ (also in form of hexahedrals) at the Gau{\ss} point with parametric coordinates ${\boldsymbol{\xi}_l}$ in the isoparametric space with the Gau{\ss} weight $w_{S,l}$ by mapping the domain of the subcell to the domain of the finite cell with the Jacobian $\bJ^{\text{sc}}_{S,\,l}$.
The Jacobian $\bJ^{\text{fc}}$ maps the finite cell to the isoparametric space. 
Analogously, the residual vector $\bbr^{\text{fc}}_K := - \sum_{S=1}^{n^{\text{sc}}}\int_{\Omega^{\text{sc}}_S} \bB^{\text{T}}\Btau\,\text{d}V_S$ in Voigt matrix notation is integrated. 
These residual vectors~$\bbr^{\text{fc}}_K$ and stiffness matrices~$\bk^{\text{fc}}_K$ are assembled to the global residual vector $\bR:=\Assemargone{K}[\bbr_K^{\text{fc}}]$ and the global stiffness matrix~$\bK:=\Assemargone{K}[\bk_K^{\text{fc}}]$.
Furthermore, the global mass matrix~$\bM:=\Assemargone{K}[\bbm_K^{\text{fc}}]$ and inertia residual vector~$\bR^\text{m}:=\Assemargone{K}[\bbr_K^{\text{fc},\text{m}}]$ are derived by assembling the mass matrices~$ \bbm_K^{\text{fc}}$ and residual vectors~$ \bbr_K^{\text{fc},\text{m}}$.
Using these matrices and vectors instead of the ones of the standard FE, the nodal displacements~$\bD$ are derived by the Newton-Raphson scheme analogously to the FEM.
%
\\
The computational advantage in the preprocessing step using voxel data is the simple generation of the regular hexahedral mesh in the domain of the scan and the decomposition of the finite cells into subcells. 
For this decomposition, multiple algorithms have been investigated in \citet{FanMisBal:2020:aso}. 
In the FC calculation, an equation system of a reduced number of degrees of freedom has to be solved, even for a complex heterogeneous structure, because the total number of degrees of freedom only depend\textcolor{blue}{s} on the number of finite cells~$n^{\text{fc}}$ and not on the number of subcells~$n^{\text{sc}}$. 
On the other hand, the assembling effort and required memory depend on the number of subcells $n^{\text{sc}}$. 
If hexahedral voxel data sets are considered, as they are usually obtained for microstructure measurements, the finite cell boundaries conform with the structural boundaries. 
Then this circumvents the problem of applying Dirichlet as well as Neumann boundary conditions onto non-conforming FCM discretizations as shown in \citet{DueParYanRan:2008:Tfc,SchRueZanBazDueRan:2012:Sal}.

\subsection{Switching Finite Elements to Subcells at the Crack Tip}
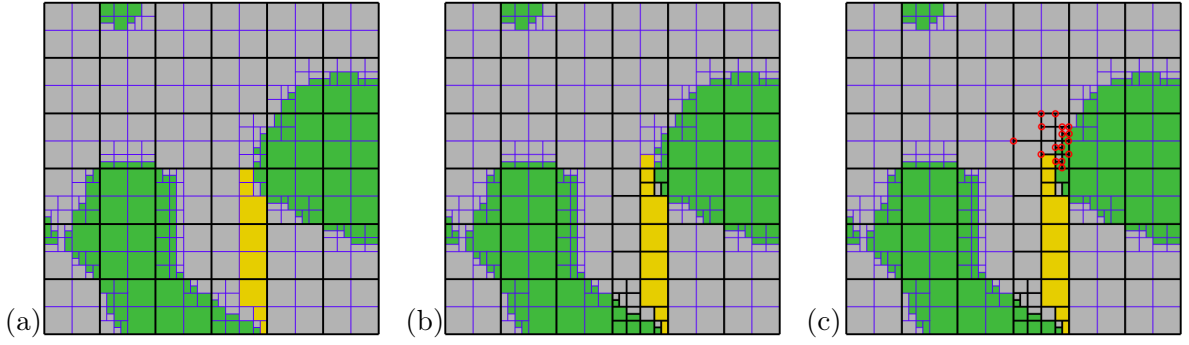
\begin{figure}[t]
\begin{picture}(26.60,130.0)
\put(-24 ,-15){\includegraphics[width=0.46\textwidth]{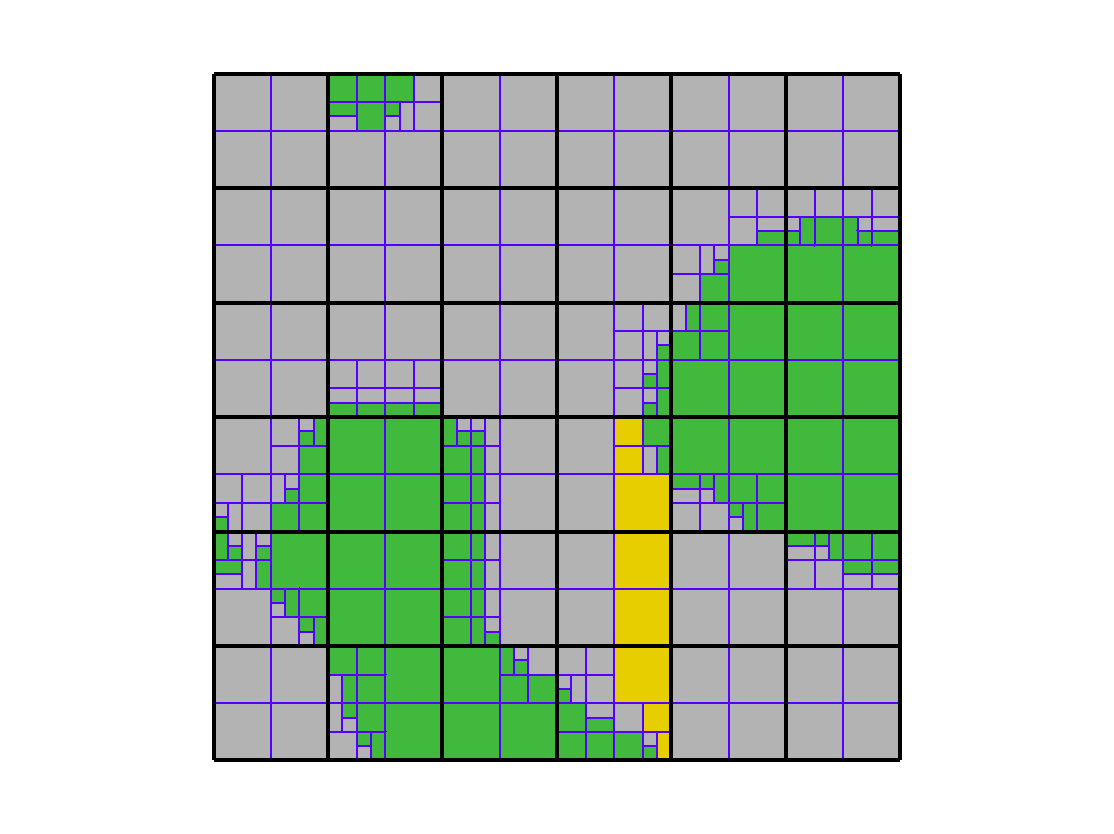}}
\put(126 ,-15){\includegraphics[width=0.46\textwidth]{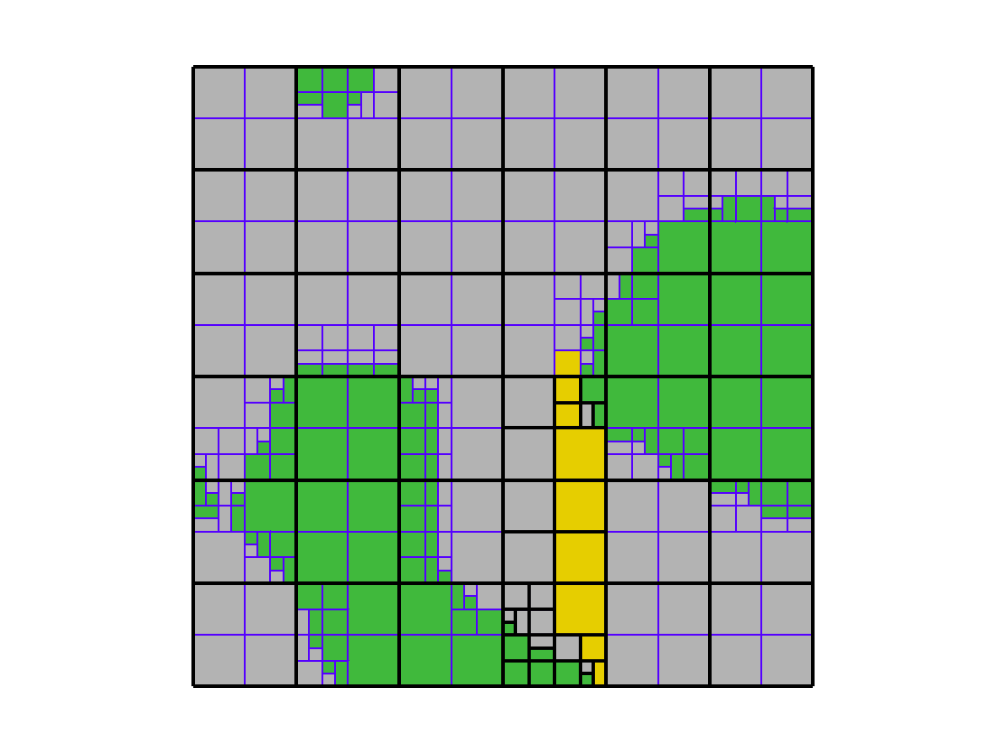}}
\put(276 ,-15){\includegraphics[width=0.46\textwidth]{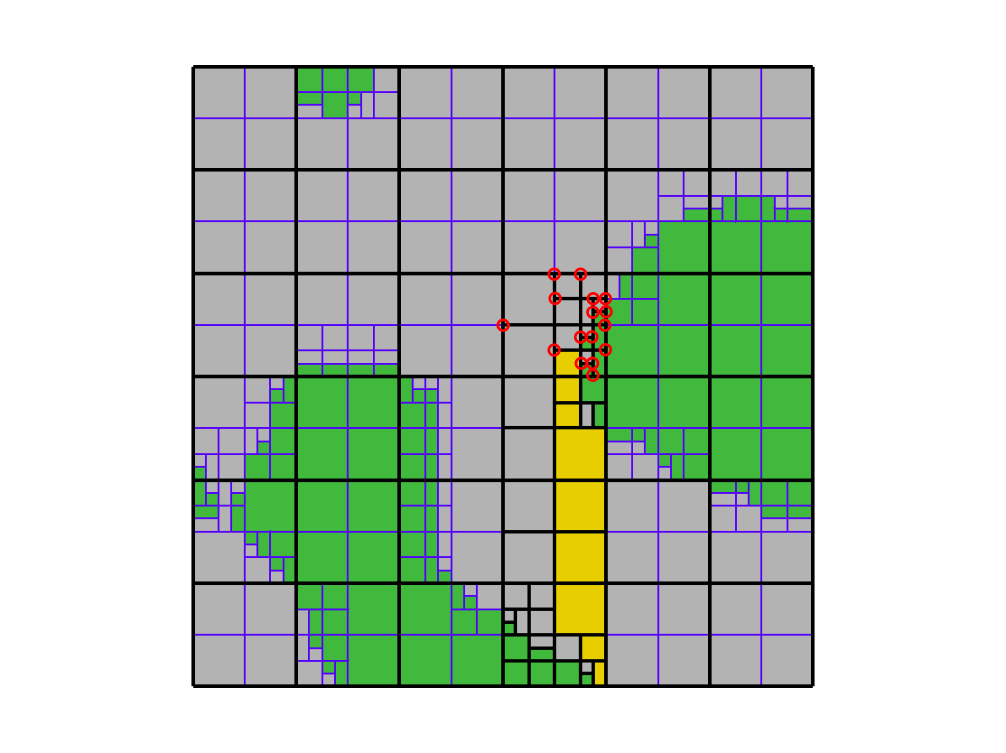}}
\put(0,0){(a)}
\put(150,0){(b)}
\put(300,0){(c)}
\end{picture}
\caption{(a) Crack (yellow) through biphasic, heterogeneous structure with finite cell/element faces (black lines) and subcell faces (blue lines) before crack propagates into the next finite cell, (b) crack enters next finite cell, i.e., the subcells of this finite cell are switched to finite elements which are separately eroded, and (c) illustration of newly added hanging nodes (red circles) resulting from switching to finite elements. Note that in this simple example, linear shape functions are considered. 
\label{fig:FCM_transformation}}
\end{figure}

Major component of the proposed algorithm is to switch all subcells of those finite cells, where erosion is detected, to finite elements and thus, to a separated approximation where also individual subcells can be eroded as elements. 
A schematic illustration of this process is given in figure~\ref{fig:FCM_transformation} where the according crack through the exemplary heterogeneous structure from figure~\ref{fig:decomposition_motivation} is depicted. 
By insertion of new finite elements replacing the subcells, hanging nodes occur as introduced in \citet{DemOdeRacHar:1989:tau,OdeDemRacWes:1989:tau,RacOdeDem:1989:tau}. 
These new nodes either hang at the side of a neighboring finite cell or at the side of a finite element which has previously been a subcell. 
The associated degrees of freedom $\bd^{H}$ at the position $\bX^{H}$ of hanging node $H$ are additionally incorporated into the global vector of nodal displacements~$\bD$. 
However, their values will not automatically match with the ones obtained in the neighboring cell, where these new degrees of freedom are not part of the interpolation. 
In order to ensure continuity of the displacement field across the hanging node, the values of the nodal displacements $\bu^{H}$ in the neighboring element/cell at $\bX^{H}$ can be computed using the standard interpolation 
\eb
\bu^{H} (\bX^{H}) = \sum\limits_I N^I (\bX^{H})\, \bd^I
\ee
with the nodal displacements $\bd_I$ in the neighboring element/cell. 
Then, the constraint condition 
\eb
\bd^H - \bu^H(\bX^H) 
= \bd^H - \sum\limits_I N^I (\bX^{H})\, \bd^I
= \bzero
\label{eq:hanging_node}
\ee
has to be fulfilled to ensure continuity in the hanging nodes~$H$. 
%
For the implementation into the FE framework, these constraints are included using the concept of Lagrange multipliers. 
To this end, the constraint terms $\Pi^H=\Blambda^H\cdot(\bd^H - \sum_I N^I (\bX^{H})\, \bd^I)$ for every hanging node~$H$ are added to the total potential energy~$\Pi$.
Herein, the Lagrange multipliers~$\Blambda^H$ can be interpreted as the interaction force vector between the hanging node and the corresponding point in the neighboring element.  
After deriving the weak form of equilibrium and consistent linearization with respect to the nodal displacements and Lagrange multipliers, the global linearized system of equations becomes
%
\eb
\underbrace{\left[
\begin{array}{ll}
(\bK+\bM) & \bC^{\text{T}} \\
\bC & \mathbf{0} \\
\end{array}
\right]}_{\bar{\bK}} 
\underbrace{\left[\begin{array}{ll} \Delta \bD \\ \Delta\boldsymbol{\lambda} \end{array} \right]}_{\Delta \bar{\bD}}=
\underbrace{\left[\begin{array}{ll}  \bR + \bR^{\text{m}} - \bQ -\bC^{\text{T}}\boldsymbol{\lambda} \\ -\bC\bD\end{array} \right]}_{\bar{\bR}}
\ee
where the global vector of increments of the Lagrange multipliers is denoted by $\Delta\Blambda$. 
The matrix~$\bC$ contains the coefficients resulting from the constraint equations and is obtained by derivation of the constraint equations~\eqref{eq:hanging_node} with respect to the nodal displacements~$\bd$.
It contains a row for each degree of freedom of each hanging node~$H$ which consists of a ``1'' in the column of the corresponding degree of freedom of the hanging node~$H$ in the original system of equations of the FE problem and the shape function~$-N^I (\bX^{H})$ in the columns of the corresponding degrees of freedom of the constraining nodes~$I$.
All other values in the matrix~$\bC$ are zero.
Within the Newton iteration, the modified linearized system of equations $\bar\bK\Delta\bar\bD = \bar\bR$ is solved with regard to~$\Delta\bar\bD$ and the degrees of freedom are updated by~$\bar\bD\Leftarrow\bar\bD+\Delta\bar\bD$, analogously to the standard FEM.
Note that the extended global matrix of the new equation system $\bar{\bK}$ is not positive definite anymore which restricts the choice of the linear equation solver. 
The Newton iteration is stopped if the stopping criterion
\eb
e:=|\bar{\bR}_{\bu}| + c_{\text{norm}}\,| \bar{\bR}_{\Blambda}| < tol
\ee
is fulfilled.
Herein, the vector norms $|\bullet|=\sqrt{\bullet\cdot \bullet}$ of the upper part of the extended global residual vector $\bar{\bR}_{\bu}:=\bR + \bR^{\text{m}} - \bQ -\bC^{\text{T}}\boldsymbol{\lambda}$ and of the lower part $\bar{\bR}_{\Blambda}=-\bC\bD$ are considered.
The numerical weighting constant $c_{\text{norm}}$ has to be chosen in such a way, that the norms of both parts of the residual lie in the same order of magnitude in the converged state in order to ensure that the Newton iteration is converged in both, the displacements $\bD$ as well as the Lagrange multipliers $\Blambda$.

\subsection{Algorithmic Implementation}
The full eigenerosion algorithm with adaptive refinement is demonstrated in figure~\ref{alg:FCM}.
We discretize the time with $n_{\text{t}}$ time steps $\Delta t = t_{k+1} - t_k$ such that the time for the whole process is $t_{\text{end}} = n_{\text{t}} \Delta t$.
In each time step, the algorithm is mainly decomposed into two major parts.
In the first part, the equilibrium equations are solved.
Therefore, the matrices and vectors~$\bk^{\text{fc}}_K$, ~$\bbr^{\text{fc}}_K$,  ~$\bbm^{\text{fc}}_K$ and~$\bbr^{\text{fc,m}}_K$\ of the finite cells are calculated as the sum of the corresponding matrices and vectors~$\bk^{\text{sc}}_S$, ~$\bbr^{\text{sc}}_S$,  ~$\bbm^{\text{sc}}_S$ and~$\bbr^{\text{sc,m}}_S$\ of the subcells.
Only the vector~$\bq_K$ of the external forces is evaluated directly because it is independent of the subcells.
Furthermore, the element stiffness matrix~$\bk_K$, residual~$\bbr_K$, mass matrix~$\bbm_K$, residual of inertia~$\bbr^{\text{m}}$ and vector of external forces~$\bq_K$ of each element~$K$, which has been a subcell previously, is evaluated.
Furthermore, the constraint equation matrix~$\bC$ is computed and, based on these matrices, the global system of equations~$\bar{\bK}\,\Delta\bar{\bD}=\bar{\bR}$ is assembled and solved for the increment of displacements~$\Delta\bar{\bD}$.
Then, the displacement vector $\bar{\bD}\Leftarrow\bar{\bD}+\Delta \bar{\bD}$ is updated.
This procedure is repeated until the residual norm~$e$ falls below the tolerance~$tol$.
\\
The second part of the algorithm in each time step deals with the eigenerosion.
Therein, the net energy gain~$-\Delta F_K$ and~$-\Delta F_S$ is evaluated for each subcell~$S$ and element~$K$, respectively.
If the net energy gain becomes larger than zero in at least one subcell or element, the one with the largest net energy gain is eroded.
Additionally, elements/subcells whose net energy gain~$-\Delta F_{K/S}$ lie within a certain tolerance to the maximum~$\max(-\Delta F_{K/S})$  are also eroded to allow for simultaneous development of crack branches.
If a subcell is eroded, the corresponding finite cell is split into multiple finite elements, each representing a former subcell.
From now on, the former finite cell does not contribute global matrices and vectors anymore.
Instead, the matrices and vectors of the new elements are taken into account.
Note, that for the eroded element, now the lumped mass matrix is considered instead of the consistent one.
The new elements are connected to new nodes.
Their displacements are initially interpolated by using the shape functions and nodal displacements of the previous finite cell.
Because some of those nodes are hanging nodes, their constraint equations are added to the matrix~$\bC$.
Thereby, the size of the global system of equations is increased by the number of degrees of freedom of the new nodes plus the number of additional constraint equations.
Furthermore, the constraint equations of hanging nodes that are constrained by the eroded element are removed.
Because the Gau\ss\, points of the elements are located at the same position as the Gau\ss\, points of the former subcells, the history variables remain the same.
Furthermore, the list of Gau\ss\,points for the calculation of the incremental crack area~$\Delta A_K$ is updated.
Afterwards, mechanical equilibrium is iteratively solved using the Newton-Raphson scheme and all intact elements and subcells are checked for erosion again.
This process is repeated until no erosion occurs anymore.
Then the  whole procedure is applied on the next time step~$t_{k+1}$ until the final time step~$t_{\text{end}}$ is reached.
\\
For imposing an initial crack, the crack propagation procedure and transformation procedure is applied on the corresponding subcells before the simulation of the first time step starts. 
Another aspect, that additionally has to be considered in the proposed algorithm, is the choice of the length scale parameter $\epsilon$ which influences the crack propagation.
The relation $\epsilon=c\,h$ with constant~$c$ and characteristic element size $h$ as proposed in \citet{PanOrt:2012:aea} leads to mesh converging results.
Empirically, it has been found that this relation also holds for this approach by considering the maximum edge length of the subcells instead of the element size $h$.
A validation is shown in the numerical examples. 

\begin{figure}[p]
\hrulefill\\
\begin{algorithm}[H]
\DontPrintSemicolon
Read input\;
Load lists of Gau\ss\, points in the neighborhood domain of each element\;
\While{\upshape time $t_k<t_{\text{end}}$}{
\While{\upshape $e>tol$ ( FE/FCM not converged)}{
\For{$K=1,$ (\upshape number of elements + number of finite cells) }{
     \uIf{$n^{\text{sc}}>1$ }{
     \For{$S=1,$ \upshape $n^{\text{sc}}$ }{
     calculation of~$\bk^{\text{sc}}_S$, ~$\bbr^{\text{sc}}_S$,  ~$\bbm^{\text{sc}}_S$ and~$\bbr^{\text{sc,m}}_S$\;
     }
     \upshape calculation of~$\bk^{\text{fc}}_K
=\sum\limits_{S=1}^{n^{\text{sc}}}\bk^{\text{sc}}_S$,
$\bbr^{\text{fc}}_K
=\sum\limits_{S=1}^{n^{\text{sc}}}\bbr^{\text{sc}}_S$,  $\bbm^{\text{fc}}_K
=\sum\limits_{S=1}^{n^{\text{sc}}}\bbm^{\text{sc}}_S$, $\bbr^{\text{fc,m}}_K
=\sum\limits_{S=1}^{n^{\text{sc}}}\bbr^{\text{sc,m}}_S$ and~$\bq^{\text{fc}}_K$\;
     }
     \Else{ calculation of $\bk_K$, $\bbr_K$, $\bq_K$, $\bbm_K$, $\bbr^{\text{m}}_K$
     }
    }
    compute constraint equation matrix~$\bC$ and assemble global system of equations~$\bar{\bK}\,\Delta \bar{\bD}=\bar{\bR}$\;
    solve system of equations and update vector of solution variables~$\bar\bD\Leftarrow\bar\bD+\Delta\bar\bD$\;}
compute net energy gain $-\Delta F_K$ and $-\Delta F_S$ for every element $K$ and subcell $S$\;
\uIf{\upshape any $-\Delta F_K>0$ or $-\Delta F_S>0$}{
   get element/subcell with $\max(-\Delta F_{K/S})$ and with $|-\Delta F_{K/S}-\max(-\Delta F_{K/S})|<tol$\;
    \For{ \upshape elements to erode}{
   \uIf{ $n^{\text{sc}}>1$}
   {transform finite cell into $n^{\text{sc}}$ finite elements\;
    change hanging node constraint matrix~$\bC$ \;
   }
   change list of Gau\ss\, points\;
   switch element property to eroded\;
   }
   go to $4$ (solve for mechanical equilibrium)\;}
   \Else{
   $ k \Leftarrow {k+1}$
   }
}
\end{algorithm}
\hrulefill
\caption{Algorithm of the proposed strategy for the simulation of crack propagation.
\label{alg:FCM}}
\end{figure}

\clearpage
\subsection{Voxel-Based Discretization and Subcell Decomposition}

Generally, the number of subcells should be as small as possible while keeping a similar accuracy to not unnecessarily increase computational effort. 
For the classical FCM, a low number of subcells is desirable in the classical FCM because the effort of volume integration within the finite cells as major part of the assembling process increases proportionally with the number of subcells. 
For strategy proposed here, where the FCM is combined with the eigenerosion, additional aspects have to be considered for the choice of the decomposition. 
Then, a small number of subcells is specifically desired because it will also increase the number of degrees of freedom in the boundary value problem whenever a crack propagates close to the material interfaces since the subcells will be transformed to finite elements. 
Therefore, the way how the cells are decomposed is quite important for the efficiency of the final strategy. 
Different decomposition schemes and their advantages and disadvantages for the FCM are proposed and discussed in \citet{FanMisBal:2020:aso}. \\
The classical method used for the decomposition is based on octree structures, where each cell is decomposed into 8 subcells of usually equal size. 
In order to simplify the notation for different decomposition techniques we use the abbreviation ``T'' in the name of the method, whenever octree structures are considered. 
Usually a threshold value is introduced to reasonably decide if the finite cell is split or not in order to neglect very small inhomogeneities. 
Then a finite cell is not split but completely assigned to the material properties of the dominant phase if a cell contains a volume fraction of the other phases falling below the given threshold value. 
This splitting procedure is continued on the subcell level until a prescribed amount of possible octree levels is reached or no cells are split anymore. 
As a drawback, a high number of subcells might be generated increasing computational effort.  
In order to arrive at a reduced number of subcells, individual neighboring subcells with equal properties can be merged to larger subcells. 
This approach will be abbreviated by ``M'' in the associated name of the specific method.\\ 
A further approach which enables the reduction of subcells to a minimum is referred to as optimal decomposition (abbreviated by ``OD''), which has been introduced in \cite{FanMisBal:2020:aso}. 
This approach works as follows: 
one direction in the voxel data is chosen, in which neighboring voxels with the same material properties in this direction are combined to connected subcells. 
This step is now repeated in the second and third perpendicular direction of the voxel data set. 
This whole procedure is executed in all possible permutations of directions and the one with the lowest number of subcells is chosen. 
As numerically shown in \cite{FanMisBal:2020:aso}, the mechanical response of FC calculations based on octree decompositions converge to the mechanical response of calculations using the optimal decomposition. 
Furthermore, this convergence is not monotonic in the octree-based computations which makes it difficult to estimate which level of octree structure is indeed sufficient. 
Hence, the optimal decomposition can be considered preferable in FCM simulations due its efficiency. 
However, as a potential drawback, the aspect ratio of dimensions of the individual subcells may become large, i.e. very thin subcells may be obtained. 
This may not be problematic for the classical FCM where large aspect ratios do not matter just for the purpose of integration. 
In the strategy proposed here, however, subcells are switched to finite elements in the crack and thus, large aspect ratios will become an issue as associated finite elements will degenerate. 
Therefore, a modified scheme will then be required.\\ 
In addition to these three decomposition techniques (T, M, OD), also combinations of two or more of these may lead to decompositions combining the advantages of the considered algorithms. 
For example, the combination ``T$[\bullet]$-OD'' first decomposes the cells using the classical octree up to the specified level $[\bullet]$ and then the optimal decomposition is performed in the resulting subcells independently. 
Thereby, a low number of subcells can be obtained which still have a reasonable aspect ratio. 
As a specialized approach, the octree method may be applied up to the level where each resulting subcell only consists of~$2\times 2 \times 2$ voxels.
Application of the optimal decomposition in these subcells will then correspond to a pure merge of the voxels, provided that voxels of equal phase are present.
Therefore, we refer to this specific approach as ``MT''.
Note that each combined decomposition approach will be as accurate as possible as long as the optimal decomposition is performed as last step. 

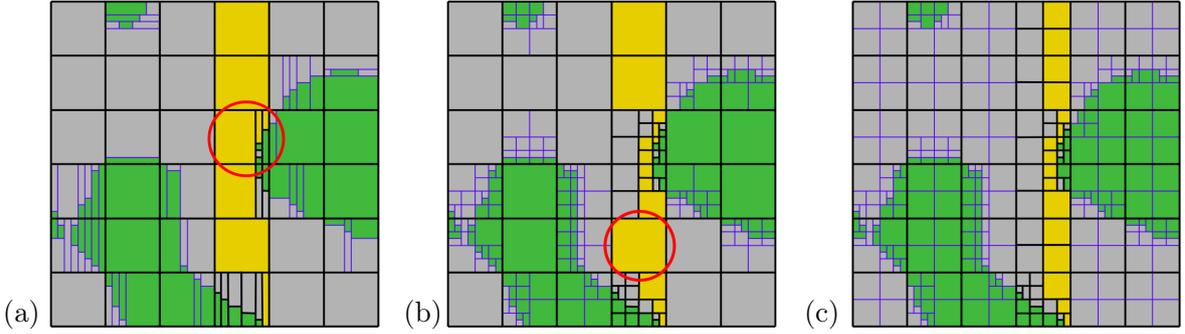
\begin{figure}[tbp]
\begin{picture}(26.60,120)
\put(-20 ,-15){\includegraphics[width=0.45\textwidth]{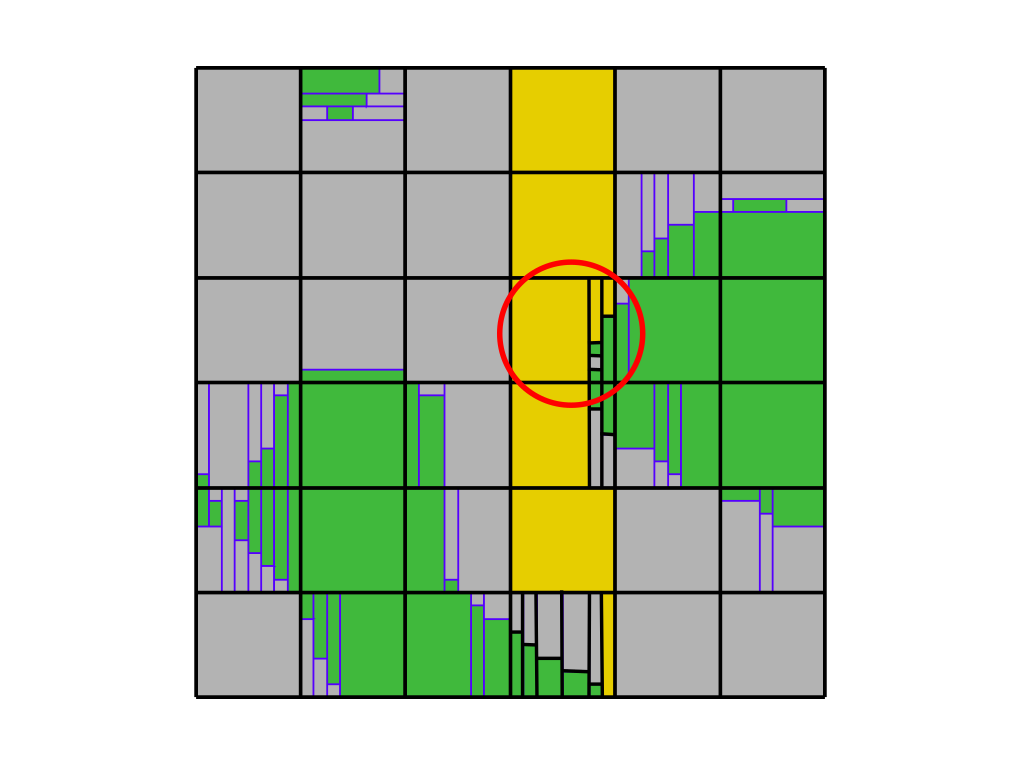}}
\put( 125 ,-15,){
\includegraphics[width=0.45\textwidth]{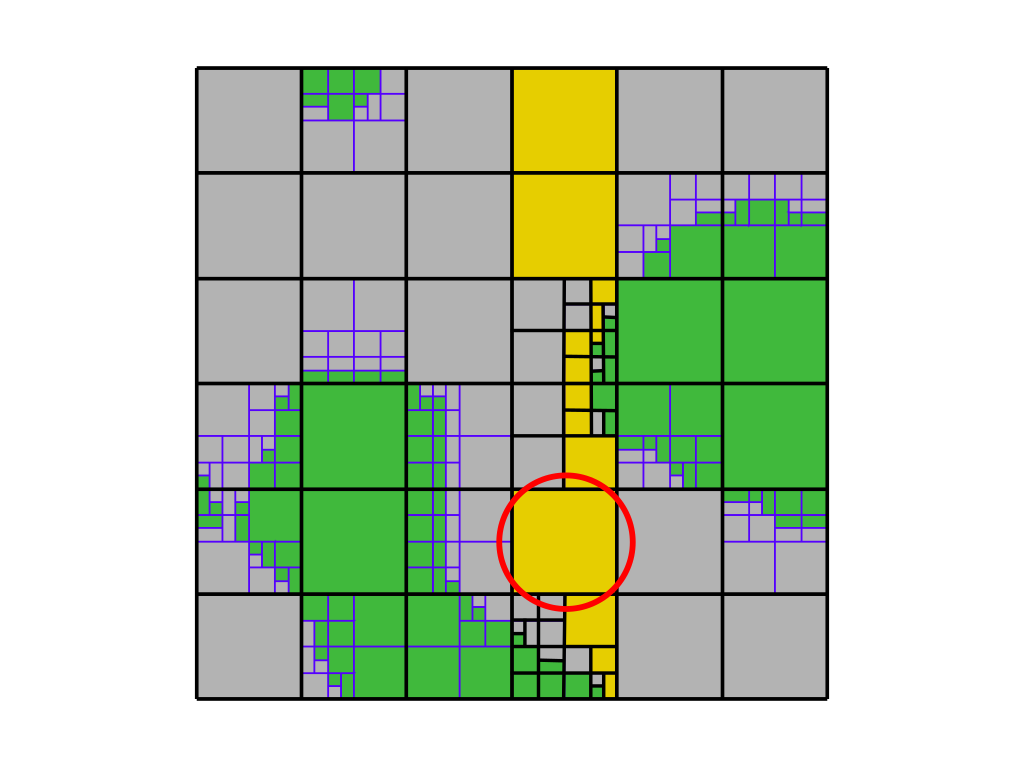}}
\put(280 ,-15){\includegraphics[width=0.45\textwidth]{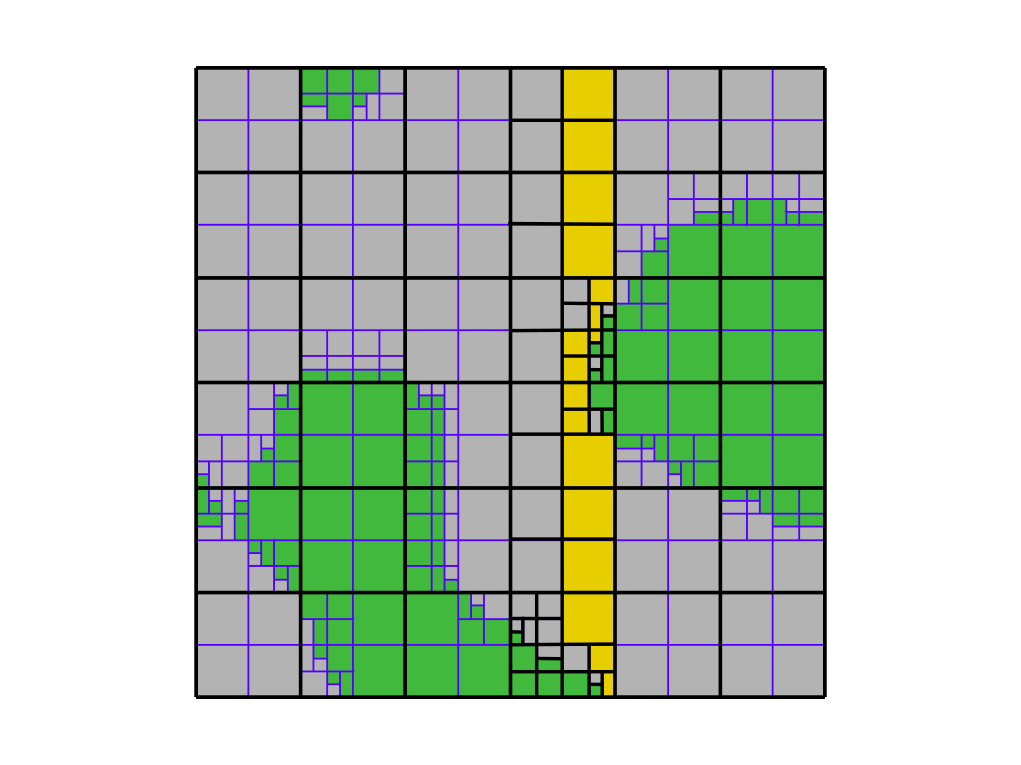}}

\put(0 ,0){(a)}
\put( 150 ,0){(b)}
\put(300 ,0){(c)}
\end{picture}
\caption{Discretization of the binarized scan from figure~\ref{fig:decomposition_motivation} with possible crack path (yellow) consisting of $48\times48$ voxels into $6\times6$ finite cells considering (a) optimal decomposition OD, (b) two levels of octree with optimal decomposition on the lowest level T2-OD = T2-MT, and (c) additional minimum split T2min1-MT.
\label{fig:decomp}}
\end{figure}

For the extended FCM combined with the eigenerosion approach, additional criteria for the voxel decomposition have to be considered due to the potential split of the finite cells into multiple finite elements.
The hanging node constraints might lead to numerical instabilities as a consequence of difficult constraint equations resulting from certain geometric distributions.
For instance, if a hanging node is connected to a point in a neighboring element where the displacements, in turn, are also interpolated based on hanging nodes themselves, then numerical instabilities occur.
This happens if both touching faces between the two neighboring elements overlap each with at least one face of another element.
Even more complex cases, in which cyclic dependencies occur, are possible.
To avoid this case, the voxel data has to be decomposed in such a way that the two faces match completely or one of the two touching faces is a subset of the other one.
If the second case occurs, only the hanging nodes of the smaller face are constrained by the nodes of the larger one.
If a mesh fulfills this requirement everywhere, we refer to it as ``consistent''.
Another aspect of the discretization deals with the crack width.
The subcells have to be arranged in such a way that the crack width may not become too large.
Therefore, the subcell sizes may not be too large near the interfaces of materials and in areas where no material boundaries occur.\\
Exemplarily, different strategies for the decomposition are presented in figure \ref{fig:decomp}.
In figure~\ref{fig:decomposition_motivation}a, a microstructure with $48\times 48$ voxels inheriting a potential crack path is considered which is discretized with $6\times 6$ finite cells.
In figure \ref{fig:decomp}a, the decomposition resulting from using optimal decomposition is presented.
Here, the given voxel set is split into the least number of subcells.
However, the microstructure is not incorporated properly in some parts near the crack due to the lengthy subcells, which may represent distorted finite elements when switching from subcells to elements.
Furthermore, the mesh is ``inconsistent'' in terms of the above-mentioned mismatch regarding the hanging nodes, which leads to numerical problems.
To circumvent these problems, the combination of octree with a merge on the lowest level (T-MT) is concluded favorable.
This also decreases the aspect ratios of many elements as demonstrated in figure~\ref{fig:decomp}b.
Note, that here the octree level has been chosen in such a way, that the merge is only applied in sets of $2\times 2$ voxels so that the mesh becomes consistent.
Therefore, the choice of pixels per structural edge and number of finite subcells per edge enables the octree decomposition and the combination with the merge $T-MT$ to be exact, because of the $8=2^3$ voxels per subcell edge which can be accurately represented by a level 3 octree.
In the generated discretization, the crack near the material interfaces is resolved sufficiently.
Nevertheless, the crack thickness becomes large far away from the material interfaces because some finite cells only consist of one single subcell due to the absence of a second material phase.
Hence, a minimum octree split of the finite cells is considered to decrease the maximum size of the subcells leading to a voxel decomposition capable of representing the crack properly in every part of the structure, cf. figure~\ref{fig:decomp}c.
However, this minimum split increases the computational effort in the assembling procedure of the global equation systems due to its increased number of subcells.
Additionally, the global equations may be increased due to additional hanging nodes if subcells are transformed into finite elements.
Because of that, the number of minimum splits is supposed to be chosen as small as possible in order to keep computational efficiency.
This inefficiency can be reduced by an appropriate management in the software implementation if only those finite cells are minimally decomposed into subcells which contain newly developing cracks.

\section{Numerical examples}

In order to show the performance of the proposed approach, two different three-dimensional numerical examples are analyzed. 
The examples are designed in the context of simulations of crack propagation through metallic microstructures. 
More specifically, two metal matrix composites are considered which are commonly used as protective layers against abrasive wear in e.g., drilling or mining tools.
These materials, however, undergo predominantly another wear mechanism, namely surface-spalling, which is mainly governed by the propagation of microscopic cracks. 
The analysis of the microscopic crack propagation is considered promising for the development of optimized materials ensuring higher protection properties and thus, longer lifespan and reduced costs. 
This is particularly important in the context of tunnel boring machines where an (unexpected) exchange of mining tools leads to a halt of the complete tunneling process representing significant costs and resources. \\
Since three-dimensional microstructure measurements are usually available as voxel data, the algorithm proposed in this paper is particularly beneficial. 
Whereas the first problem considers a more academic example of a simplified metal matrix composite microstructure serving as benchmark problem to enable a more generally meaningful analysis, the second problem addresses a real metal-matrix microstructure obtained from micro-CT.

\subsection{Benchmark Problem}

\begin{figure}[t]
\begin{picture}(26.60,178.0)
\put(00 ,00){\includegraphics[width=0.48\textwidth]{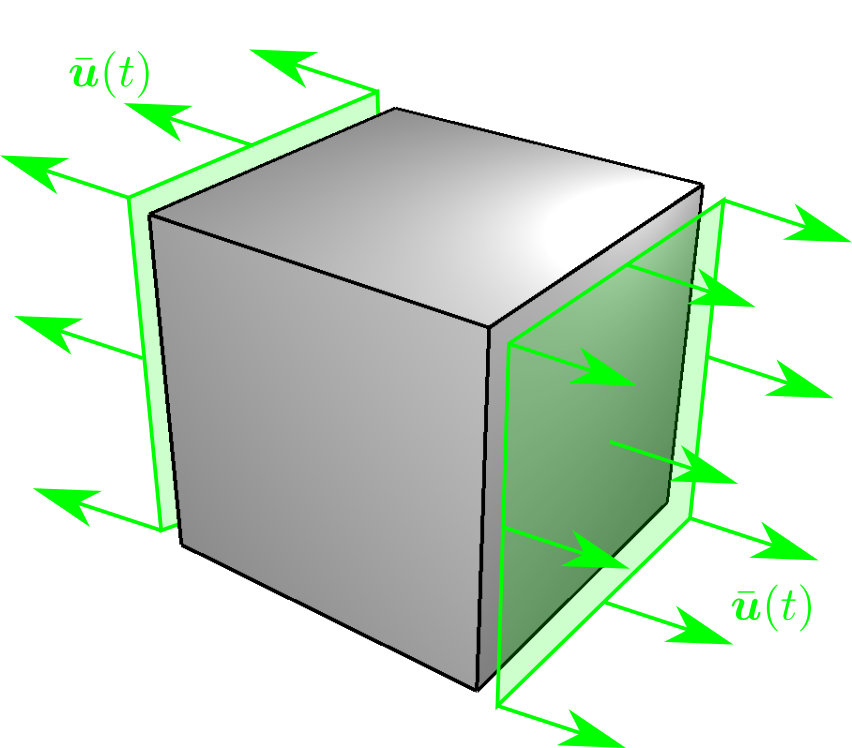}}
\put(230 ,00){\includegraphics[width=0.45\textwidth]{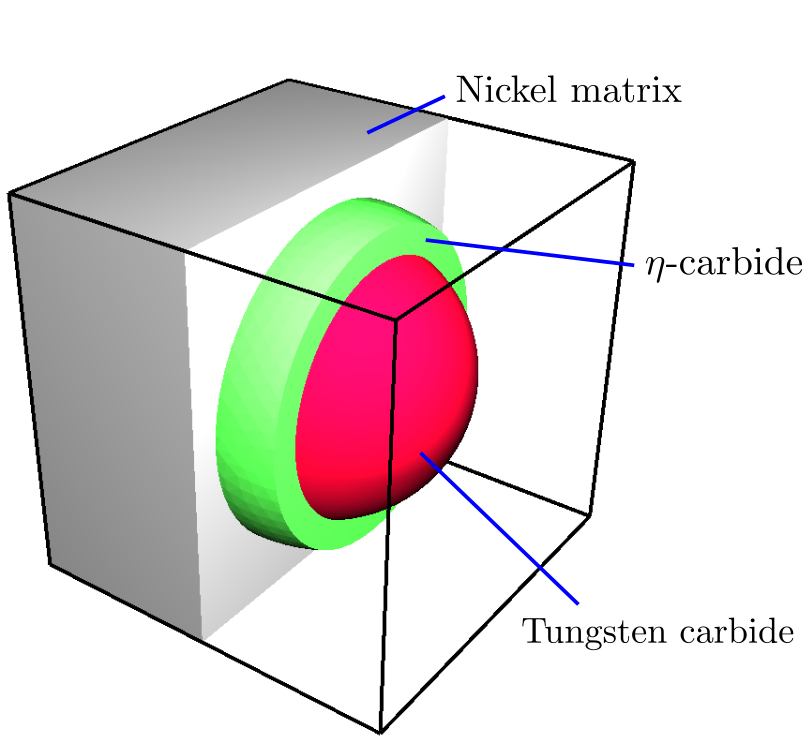}}
\put(0 ,0){(a)}
\put(220 ,0){(b)}
\end{picture}
\caption{(a) Boundary conditions of the simulation setup of the benchmark experiment and (b) the considered microstructure}
\label{fig:benchmark_setup}
\end{figure}

In order to proof the general feasibility of the proposed algorithm, a benchmark experiment on an artificial metal matrix composite microstructure as seen in figure~\ref{fig:benchmark_setup} is investigated.
Although more realistic, simplified microstructures could in principle be constructed by applying the concept of statistically similar representative volume elements as proposed in \citet{SchBalBraSch:2015:do3}, the microstructure here may already be considered sufficiently realistic for the purposes in this section. 
The microstructure has an edge length of $70$ $\mu$m and consists of a spherical tungsten carbide inclusion with a diameter of $50$ $\mu$m surrounded by a ductile nickel matrix.  
A brittle $\eta$-carbide layer with a thickness of $5$ $\mu$m lies in-between these phases. 
For comparative purposes, crack propagation through this microstructure is simulated based on eigenerosion and the following different types of approaches: 
\begin{itemize}
\item[(i)] As a reference representing the rather classical approach, an unstructured mesh using quadratic 10-node tetrahedral elements is considered, where the element faces conform with the assumed spherical morphology of the inclusion. 
\item[(ii)] As further reference, reflecting rather the voxel-based setting, a structured, regular FE mesh is analyzed, where every voxel is discretized with one $8$-node hexahedral element. 
\item[(iii)] As a more efficient, classical alternative, a semi-regular FE mesh is investigated, where $27$-node hexahedral elements considering local mesh refinement at the material interfaces using hanging nodes are taken into account. 
\item[(iv)] The results based on the proposed algorithm combining the eigenerosion and FCM is compared with the other calculations. 
\end{itemize}
The latter three types are based on artificially generated voxel data consisting of $56\times56\times56$ voxels which represents a virtual measurement of the morphology of the spherical inclusion in terms of voxels. 
For the proposed approach (iv), various voxel decomposition schemes are applied, all based on the discretization using $7\times 7\times 7$ finite cells such that each finite cell contains $8\times 8\times 8$ voxels. 
The resulting subcell decompositions are used for approach type (iii) to define the locally refined finite element discretization. 
The decomposition schemes are chosen such that the material boundaries are represented accurately in order to avoid an inaccurate representation as source for potentially erroneous results. 
\begin{figure}[t]
\begin{picture}(26.60,370.0)
\put(-30,190){\includegraphics[width=1.1\textwidth]{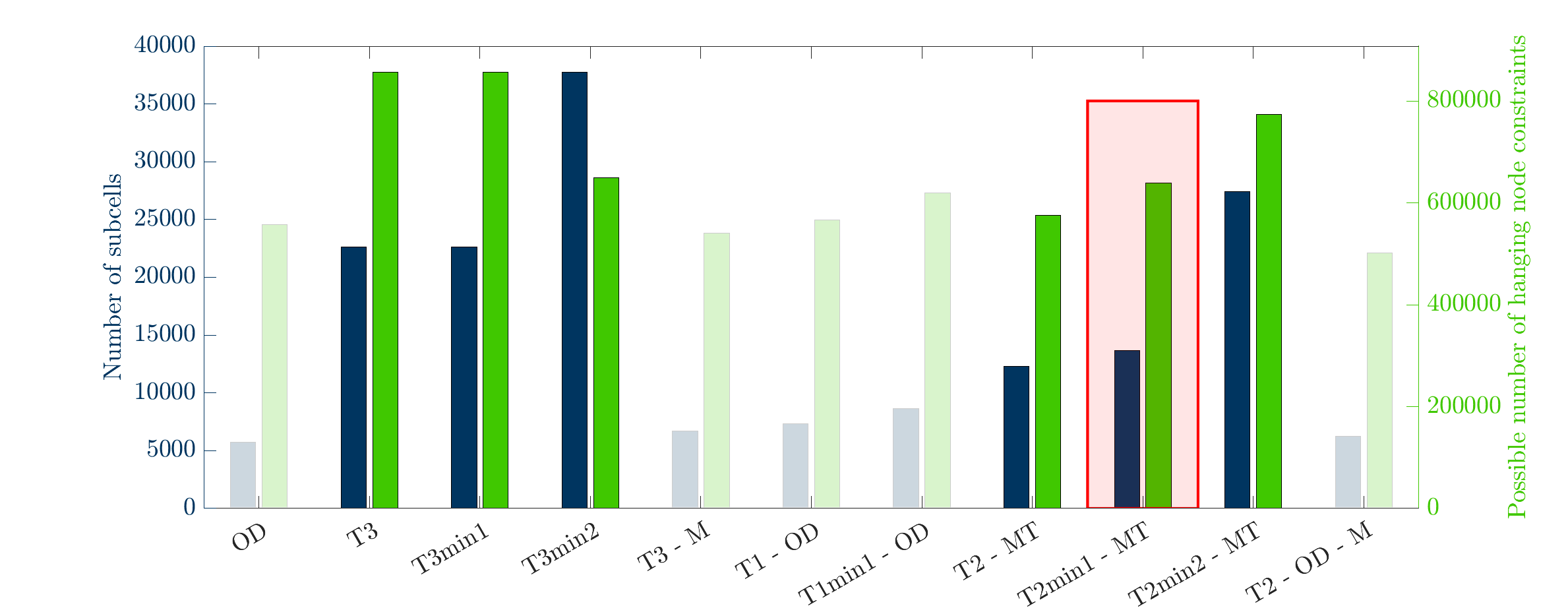}}
\put(-30,00){\includegraphics[width=1.1\textwidth]{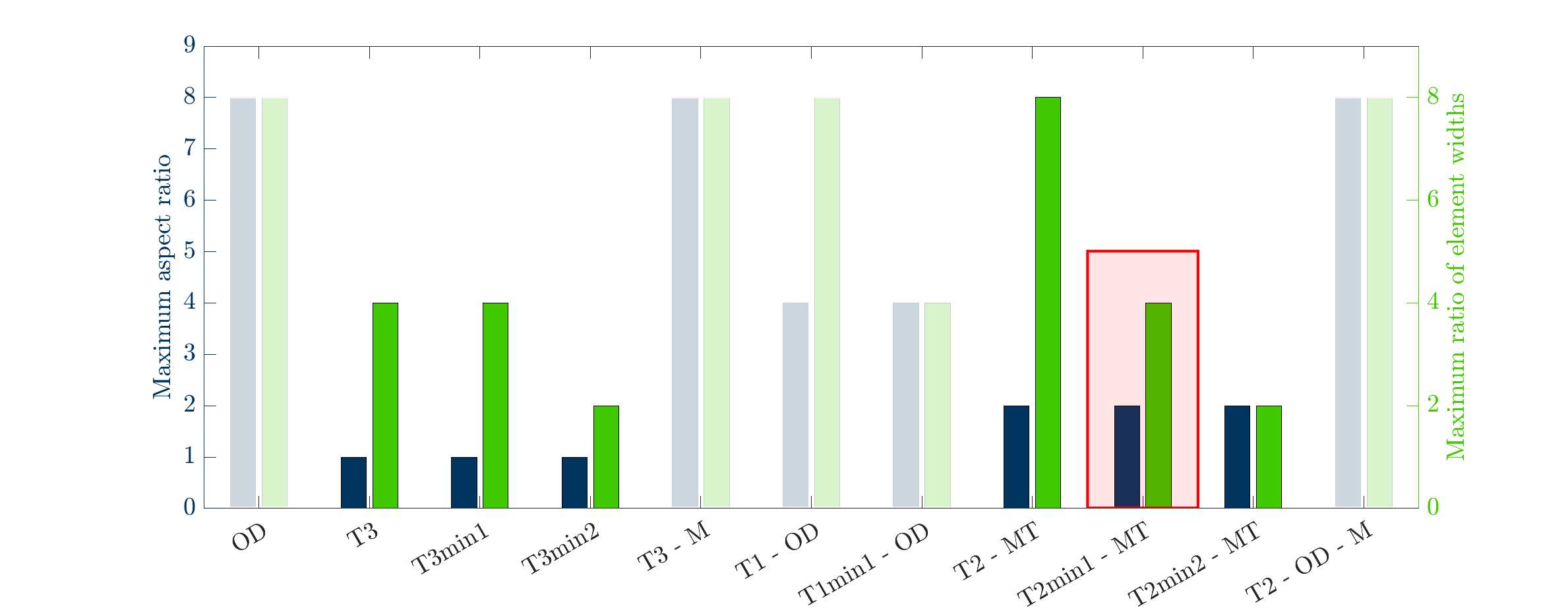}}
\put(0,0){(b)}
\put(0,190){(a)}
\end{picture}
\caption{Comparison of different decomposition schemes based on the voxel data of the artificial microstructure shown in figure~\ref{fig:benchmark_setup}: (a) total number of subcells and number of hanging node constraints, and (b) maximum value of aspect ratios of the subcells' edge lengths and ratio of maximum edge length (over all subcells) per minimum edge length (over all subcells). The transparency indicates that the corresponding discretization contains inconsistencies if the discretization type iii or iv are applied. If those discretizations are used, numerical problems may occur.  T2min1-MT = T2min1-OD is chosen for the numerical calculations.}
\label{fig:benchmark_discretization}
\end{figure}
The properties of the resulting discretizations considering quadratic shape functions are shown in figure~\ref{fig:benchmark_discretization}. 
As expected, the discretizations generated with the decomposition methods OD, T2-OD-M and T3-M result in the least numbers of subcells, while containing subcells with the highest possible aspect ratio. 
This may cause numerical instabilities in the case when the proposed algorithm is used as soon as the subcells are transformed to finite elements at the crack tip, or in the case of using discretization type iii. 
Furthermore, these discretizations contain inconsistencies if the subcells are transformed into finite elements.
These may lead to numerical problems if discretization types iii or iv are considered.
Hence, these discretizations should not be considered.
In contrast to that, the pure octree decompositions T3, T3min1 and T3min2 lead to a high number of subcells and thus, hanging nodes, but the element aspect ratio is kept $1$.
Thereby, the ratio of the maximum edge length over the whole microstructure problem per minimum edge length of all subcells remains low.
Additionally, the equidistant octree decomposition always leads to consistent discretizations.
The combination T2min1-MT enables a subcell decomposition with a small number of subcells and hanging node constraints on the one hand, while keeping the aspect ratios moderate on the other hand.
Furthermore, this approach supplies a consistent mesh at every time in the simulation with the proposed algorithm.
Therefore, this decomposition scheme is considered for the subsequent mechanical simulations. 

For the numerical calculation of the benchmark problem, the material parameters of the constituents as shown in table \ref{tab:benchmark_material_parameters} are considered. 
These have been chosen in line with experimental findings on the single components. 
\begin{table}[t]
\resizebox{\textwidth}{!}{%
\begin{tabular}{|gccccccc|}
 \hline
 \rowcolor{rub_blue40}
& $K\,[\mathrm{GPa}]$               & $\mu\,[\mathrm{GPa}]$ & $y_0\,[\mathrm{GPa}]$ & $y_{\infty}\,[\mathrm{GPa}]$ & $h^{\mathrm{exp}}$ & $h^{\mathrm{lin}}\,[\mathrm{GPa}]$ & $G_c\,[\mathrm{N/mm^2}]$ \\
Tungsten carbide &$308.12$ & $288.71$ & $10^{12}$ & $10^{12}$ & $0$ & $0$ & $0.0371$ \\
$\eta$-carbide & $394.38$ & $228.72$ & $10^{12}$ & $10^{12}$ & $0$ & $0$ & $0.0065$\\
Nickel & $225.6$ & $75.19$ & $260$ & $580$ & $9$ & $70 $ & $1.730$\\
\hline
\end{tabular}
}
\caption{Material parameters of metal matrix composite microstructure of the benchmark experiment.}
\label{tab:benchmark_material_parameters}
\end{table}
\begin{figure}[t]
\begin{picture}(26.60,175.0)
\put(-10,0){\includegraphics[width=0.55\textwidth]{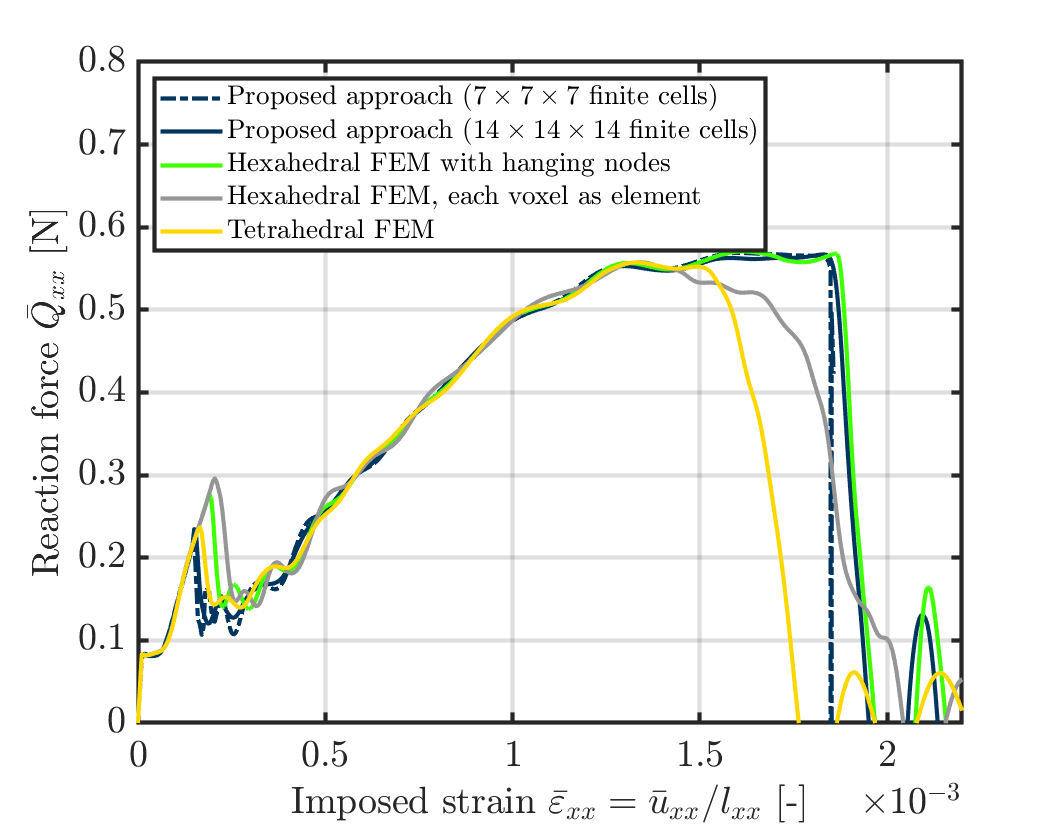}}
\put(220,0){\includegraphics[width=0.55\textwidth]{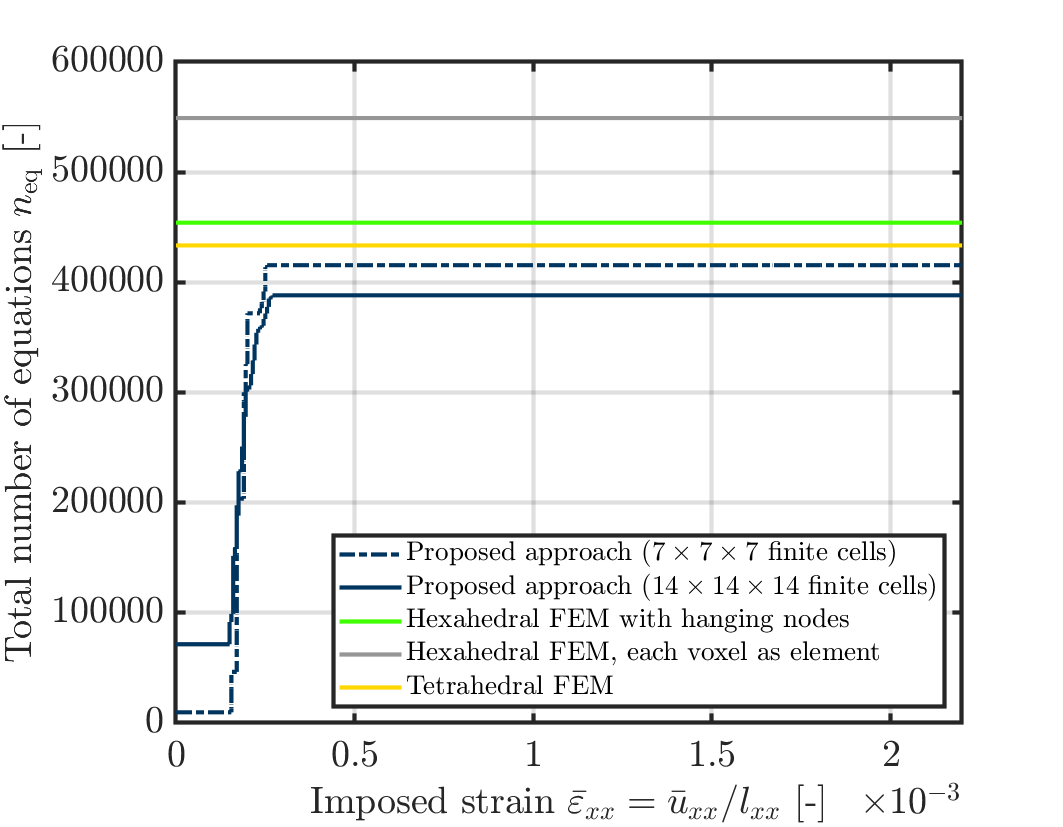}}
\put(0,0){(a)}
\put(240,0){(b)}
\end{picture}
\caption{(a) Resultant reaction force~$\bar{Q}_{xx}$ at the Dirichlet boundary versus imposed specimen strain and (b) number of equations within the linearized system of each Newton-Raphson step versus imposed specimen strain~$\bar{\varepsilon}_{xx}$. 
\label{fig:benchmark_reac_eq}}
\end{figure}
Displacements at the opposite surfaces with normal in $x$-direction, as shown in figure~\ref{fig:benchmark_setup}a, are linearly increased over time with a velocity of $\dot{\bar{u}}=350\,\mathrm{mm/s}$ such that predefined strains measured through the change of outer dimensions of the specimen $\bar{\varepsilon}_{xx}=\bar{u}_{xx}/l_{xx}$ are reached.
To further analyze the choice of polynomial degree in the finite cells, in addition to the tri-quadratic shape functions used in T2min1-OD, also tri-linear shape functions are considered. 
To this end, the benchmark microstructure is additionally decomposed using $14\times 14 \times 14 $ finite cells and approach T1-OD. 
Regarding discretization type ii, where all voxels are represented as single hexahedral finite elements, only elements with tri-linear shape functions are considered due to the high amount of degrees of freedom. 
Therefore, this calculation should not be misinterpreted as reference solution since the approximation quality is comparatively poor. 
For the eigenerosion regularization, the influence radius is chosen to $\varepsilon=0.5\,h$ in all simulations. 
Note, that for the brittle components in the microstructure, the initial yield stress~$y_0$ is set to a high value, which will not be reached in the simulation, to prohibit plastic deformations. \\
For the analysis of the numerical results, the resultant reaction force~$\bar{Q}_{xx}$ at the Dirichlet boundary versus imposed specimen strain~$\bar{\varepsilon}_{xx}$ is shown in figure~\ref{fig:benchmark_reac_eq}a. 
With increasing strain, the reaction force increases until the $\eta$-carbide layer cracks at the strain $\bar{\varepsilon}_{xx}=0.2 \cdot 10^{-3}$ and thus, the reaction forces drop rapidly. 
Afterwards, the structural response increases again until the nickel matrix breaks into two parts so that the reaction force decreases to zero. 
As expected, the results show a quite different response for the hexahedral FEM where each voxel is represented by one element, which is due to the lack in approximation accuracy of the tri-linear shape functions. 
In contrast to that, the proposed approach converges well with increasing number of finite cells to the hexahedral FE calculation using local refinement everywhere.
Whereas the reaction force using $7\times 7 \times 7$ finite cells still differs significantly from the hexahedral FE calculation with local mesh refinement in the region where the matrix breaks, it is almost identical when using $14\times 14\times 14 $ finite cells. 
Hence, this benchmark experiment shows the mesh independency of the proposed approach. 
\\
In addition to these calculations, also the results of the conforming, unstructured tetrahedral FE calculation are depicted and show a quite different quantitative response in the reaction force. 
This is not surprising since the considered microstructure morphology is different. 
Whereas here the inclusion shape is almost perfectly spherical, the shape is non-smooth for the FC calculation where voxel data is directly considered. 
Note that this does not represent a shortcoming of the proposed approach, it rather illustrates a general challenge when simulating heterogeneous structures which are solely given as voxel data. 
By constructing a somewhat interpolated interface morphology serving as reference for the generation of the conforming, unstructured FE mesh, the resulting morphology may not necessarily correspond to the ``real'' one. 
In fact, the FCM itself corresponds to a somewhat interpolated interface morphology, however differently interpolated, such that it can not be said which kind of interpolation works better in general. 
Here, the interface morphology is known because a spherical inclusion is considered and thus, no interpolation based on the voxel data is required. 
Therefore, the difference in the response of all voxel-based calculations compared to the conforming FE calculation just reflects the limitation of using voxel data of limited resolution. 
However, such limited resolutions are standard in real three-dimensional measurements of heterogeneous structures, e.g., based on micro-CT.\\ 
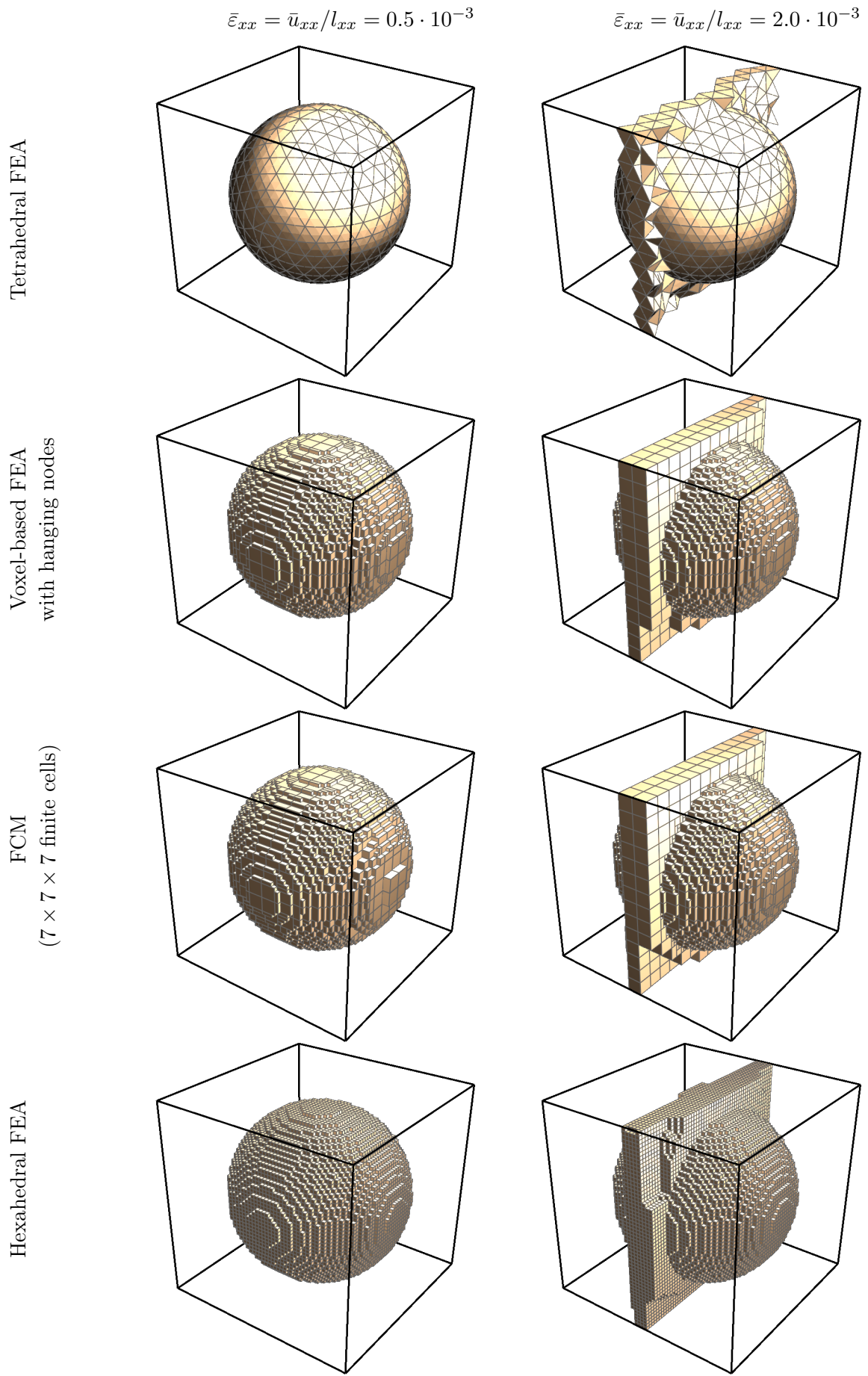
\begin{figure}[hbtp]
\unitlength1cm
\begin{picture}(16,24)
\put(0.5, 19.0){\rotatebox{90}{ Tetrahedral FEA}}
\put(0.5, 13.5){\rotatebox{90}{ Voxel-based FEA}}
\put(1.0, 13.3){\rotatebox{90}{ with hanging nodes}}
\put(0.5,  9.0){\rotatebox{90}{ FCM}}
\put(1.0,  7.5){\rotatebox{90}{ ($7\times 7 \times 7$ finite cells)}}
\put(0.5,  2.0){\rotatebox{90}{ Hexahedral FEA}}
\put(4.2, 24.0){ $\bar{\varepsilon}_{xx}=\bar{u}_{xx}/l_{xx}=0.5\cdot10^{-3}$}
\put(11.0,24.0){ $\bar{\varepsilon}_{xx}=\bar{u}_{xx}/l_{xx}=2.0\cdot10^{-3}$}
\put( 9.8,17.7){\includegraphics[width=0.38\textwidth]{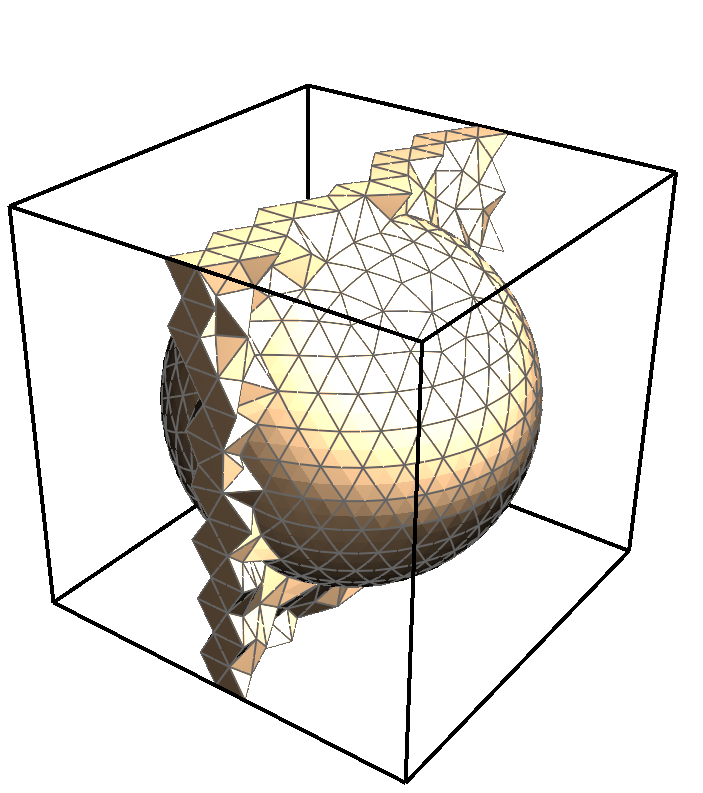}}
\put( 9.8,11.8){\includegraphics[width=0.38\textwidth]{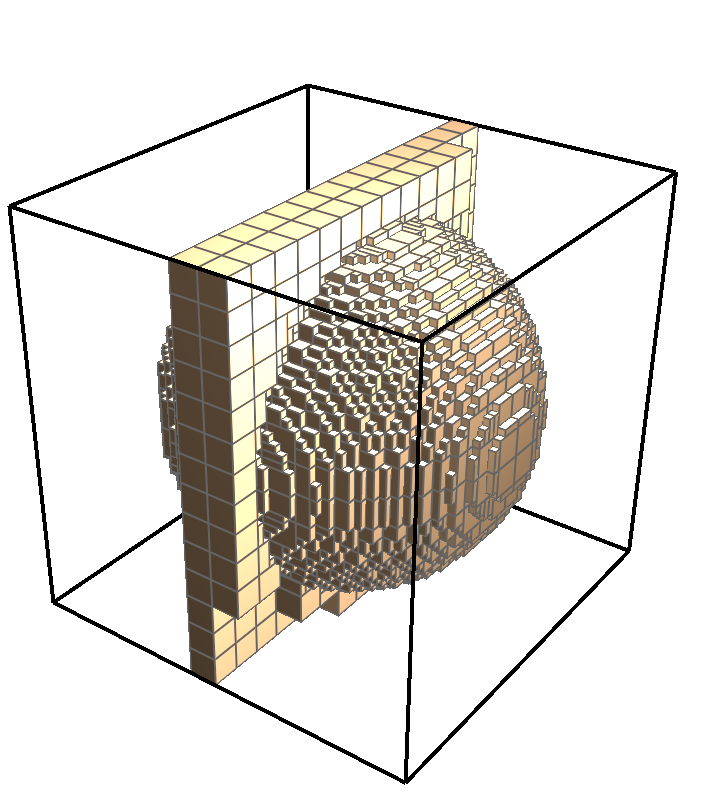}}
\put( 9.8, 5.9){\includegraphics[width=0.38\textwidth]{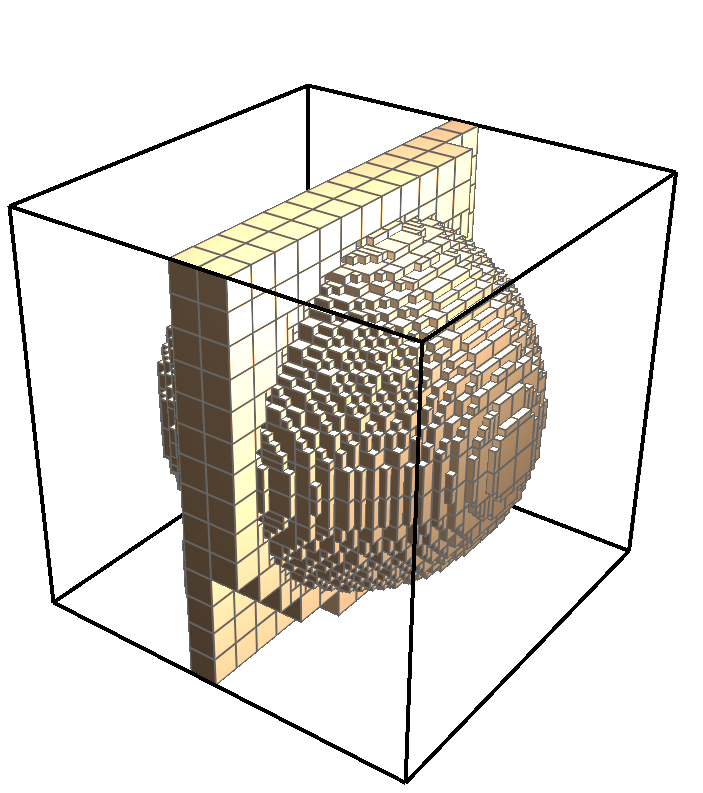}}
\put( 9.8, 0.0){\includegraphics[width=0.38\textwidth]{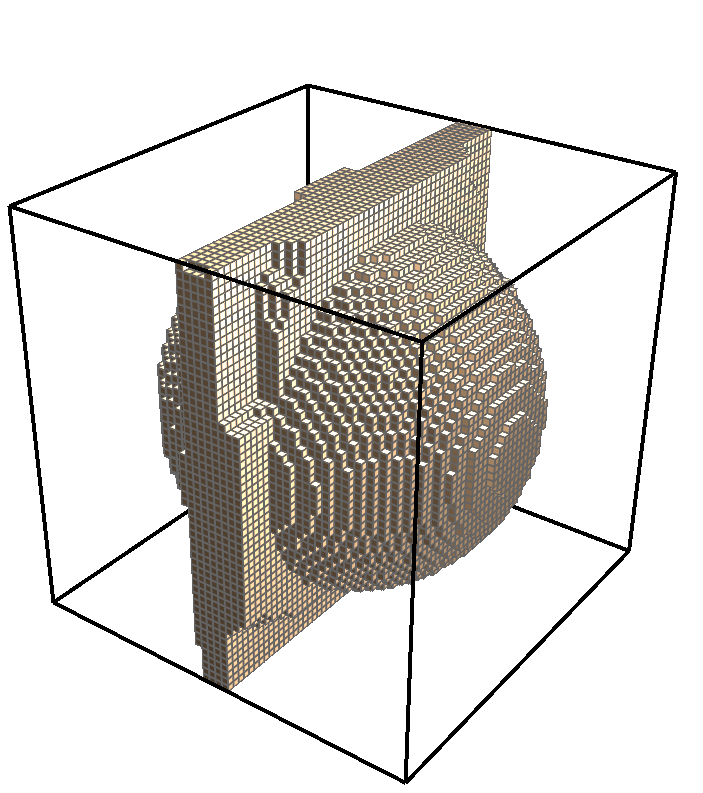}}
\put(3.0,17.7){\includegraphics[width=0.38\textwidth]{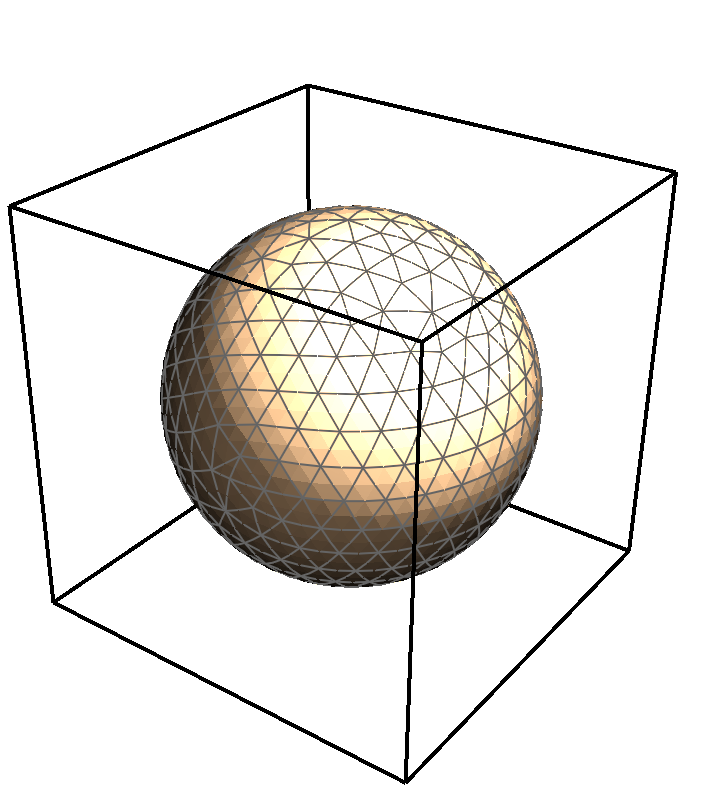}}
\put(3.0,11.8){\includegraphics[width=0.38\textwidth]{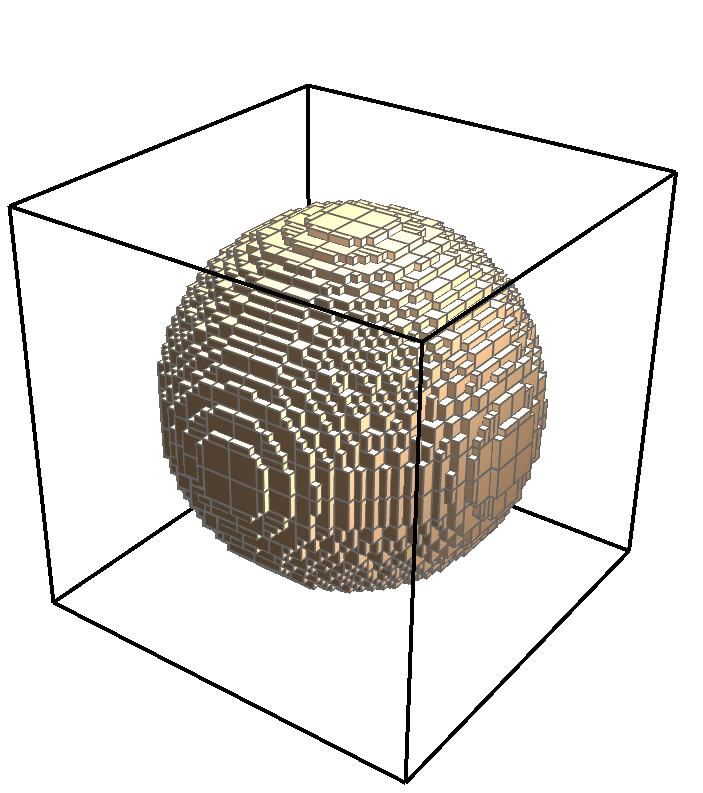}}
\put(3.0, 5.9){\includegraphics[width=0.38\textwidth]{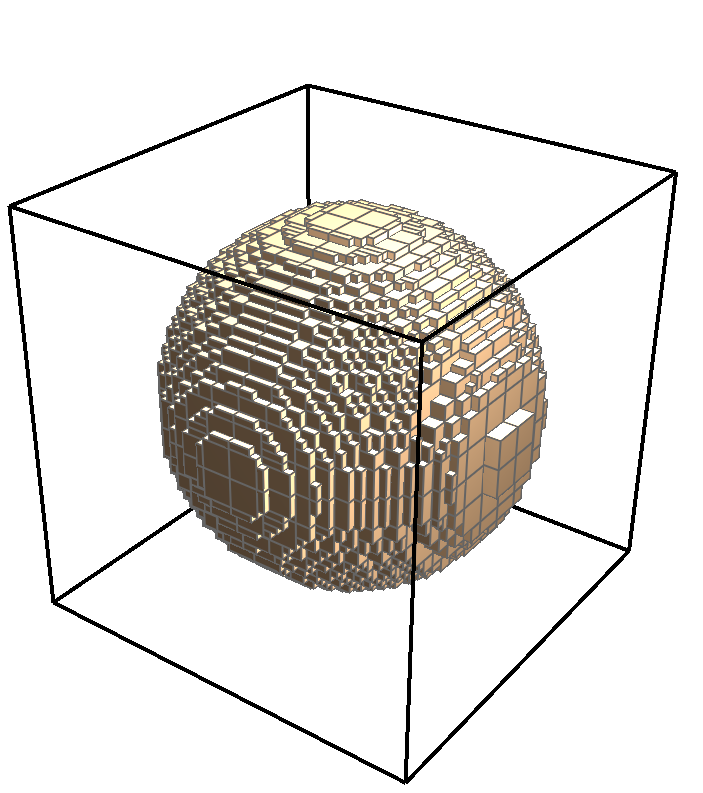}}
\put(3.0, 0.0){\includegraphics[width=0.38\textwidth]{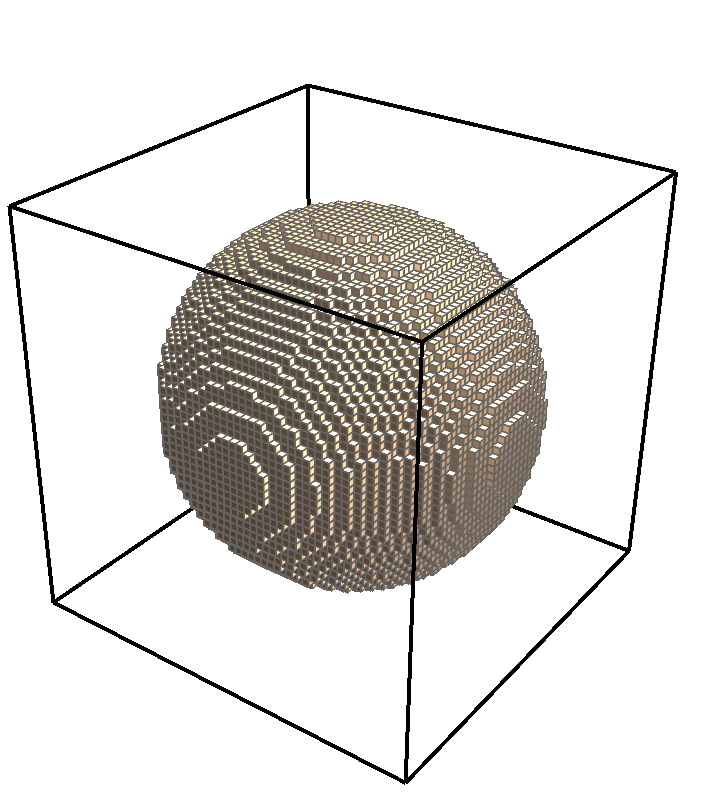}}
\end{picture}
\caption{Crack path through artificial metal matrix composite microstructures at a specimen strain of $\bar{\varepsilon}_{xx}=0.5\cdot10^{-3}$ (left) and $\bar{\varepsilon}_{xx}=2.0\cdot10^{-3}$ (right).}
\label{fig:benchmark_crack_path}
\end{figure}
Aside from the quantitative response, also qualitatively the proposed approach corresponds quite well with all other calculations. 
This can be seen by comparing the crack paths as shown in figure~\ref{fig:benchmark_crack_path}. 
There, the eroded elements are depicted at $\bar{\varepsilon}_{xx} = 0.5\cdot10^{-3}$ (after complete fracture of the $\eta$-carbide layer) and at $\bar{\varepsilon}_{xx} = 2.0\cdot10^{-3}$ (after complete failure of the specimen) and show significant similarities, again reflecting the mesh independence. 
Furthermore, it is shown that one constant $c$ for determining the influence radius $\varepsilon$ suffices even for different element and discretization types.\\ 
\begin{figure}[t]
\begin{picture}(26.60,190.0)
\put(20,0){\includegraphics[width=0.45\textwidth]{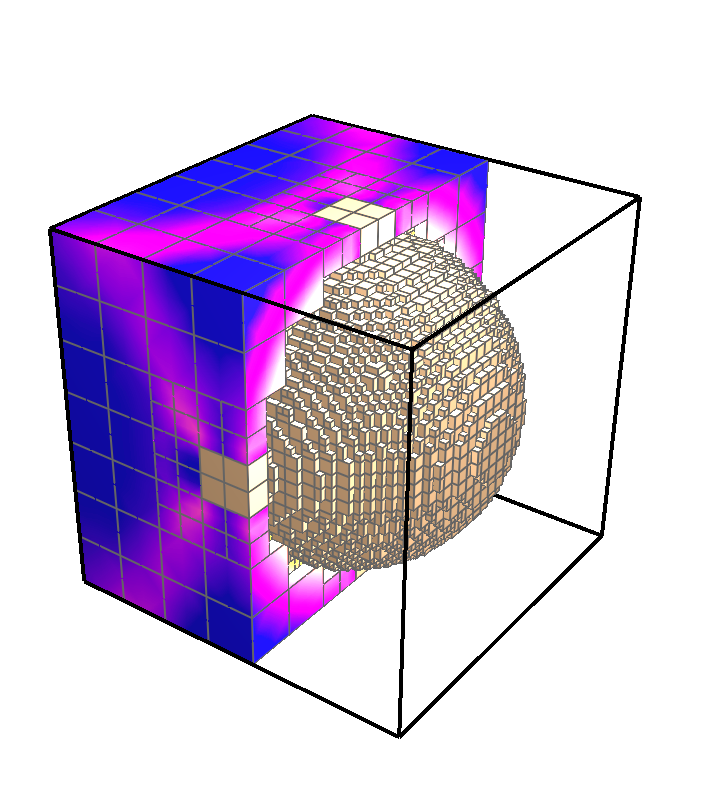}}
\put(0,0){(a)}
\put(00 ,25){\rotatebox{90}{ \includegraphics[width=150pt ,height=10 pt]{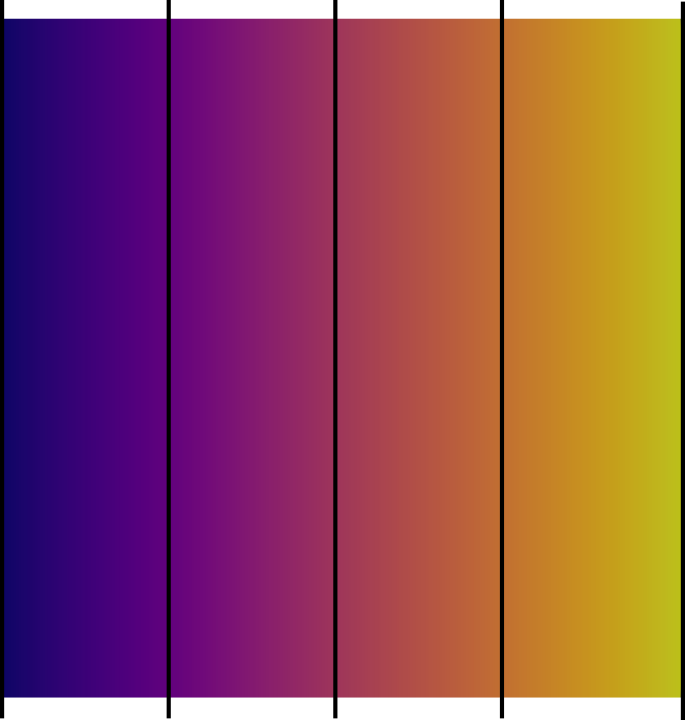}}}
\put(15,28){\footnotesize{$0$}}
\put(15,63){\footnotesize{$2.5$}}
\put(15,100){\footnotesize{$5.0$}}
\put(15,137){\footnotesize{$7.5$}}
\put(14,174){\footnotesize{$10.0$}}
\put(0,185){\footnotesize{$\times 10^{-5}$}}
\put(260,0){\includegraphics[width=0.45\textwidth]{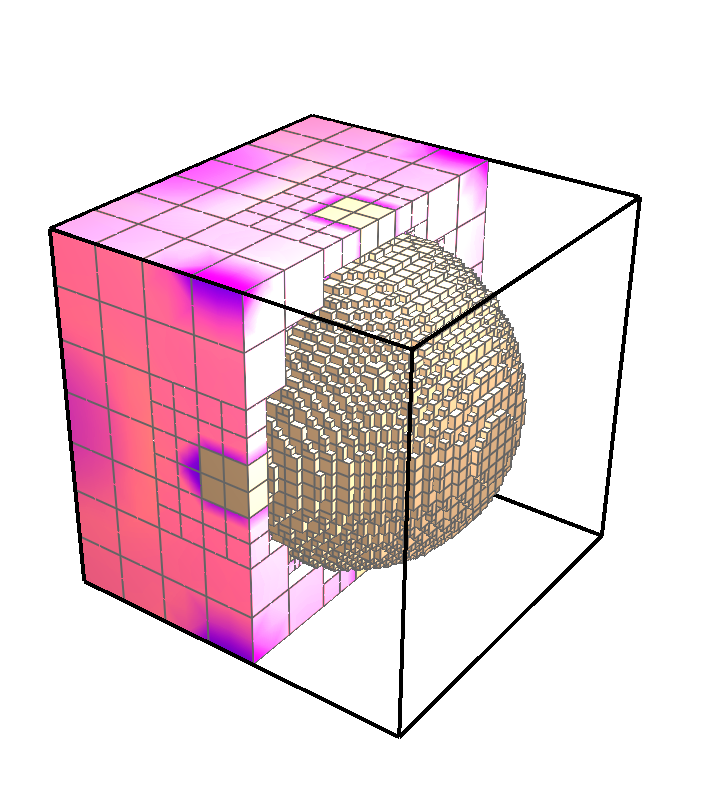}}
\put(240,0){(b)}
\put(240 ,25){\rotatebox{90}{ \includegraphics[width=150pt ,height=10 pt]{\figpath/colorbar_lines.png}}}
\put(255,28){\footnotesize{$0$}}
\put(255,63){\footnotesize{$100$}}
\put(255,100){\footnotesize{$200$}}
\put(255,137){\footnotesize{$300$}}
\put(255,174){\footnotesize{$400$}}
\put(240,185){\footnotesize{$\mathrm{MPa}$}}
\end{picture}
\caption{(a) Distribution of equivalent plastic strain~$\alpha$ and (b) von Mises Kirchhoff stress~$\tau^{\mathrm{vM}}$ over the specimen resulting from the simulation using the proposed approach at $\bar{\varepsilon}_{xx}=\bar{u}_{xx}/l_{xx}=0.18\%$ considering~$7\times 7 \times 7$ finite cells; the crack path is shown in gold color.
Here, the grey lines indicate either the faces of the active finite cells/subcells or finite elements in the case when the subcells are transformed to elements in the crack path. 
\label{fig:benchmark_vm_alpha}}
\end{figure}
The resulting equivalent plastic strain $\alpha$ and the von Mises Kirchhoff stress~$\tau^{\mathrm{vM}}$ resulting from the simulation using the proposed approach are shown in figure~\ref{fig:benchmark_vm_alpha}. 
Even if small specimen strains $\bar{\varepsilon}_{xx}$ are applied, plastic deformations occur due to microscopic strains larger than the macroscopic ones caused by the morphological heterogeneity. 

The main advantage of the proposed approach is its efficiency compared to the comparative calculations. 
This can be seen in the size of the resulting linearized system of equations to be solved in each Newton iteration, which is depicted in figure \ref{fig:benchmark_reac_eq}b. 
In the beginning of the deformation process, when no crack evolves and only finite cells and no elements including hanging nodes occur, the equation system is lower by a factor of $2.5$ compared to the final state in which all finite cells at the material boundaries are transformed into single elements.
In the regular hexahedral mesh $549,081$ equations occur.
If quadratic shape functions were considered, which would actually be needed for accuracy reasons, there would be $4,303,153$ equations and thus an increase of factor $60.1$. 
Furthermore, the number of elements for the assembling with $175,616$ elements is larger by a factor of $7.5$ compared to the $23,360$ subcells/elements. 
Additionally, the number of eroded elements increases the computational effort because for every eroded element, the system of nonlinear equations associated with mechanical equilibrium has to be solved again as part of the eigenerosion algorithm.
The final crack contains $4,046$ elements using the proposed approach and $18,406$ elements in the regular hexahedral mesh of the discretization type ii.
Especially with respect to the number of $600$ time steps, the large numbers of eroded elements strongly influence the number of solving steps and thus the computational effort.
This again demonstrates the gain in efficiency of the proposed approach compared to the discretization type ii.
Note that in this example the benefits of the extended FCM compared to the semi-regular hexahedral mesh with hanging nodes is even relatively small because almost all finite cells containing the $\eta$-carbide layer become part of the crack and all included subcells are thus transformed to finite elements. 
This will be much different and to the advantage of the proposed approach when considering more complex microstructures, where only a small fraction of finite cells become part of the cracks. 
Nevertheless, the benchmark experiment demonstrates the efficiency of the proposed approach compared to the alternative methods while being competitively accurate.

\subsection{Microstructure Based on Micro-CT Scan}
\begin{figure}[hbtp]
\begin{picture}(26.60,200.0)
\put(0,0){\includegraphics[width=0.50\textwidth]{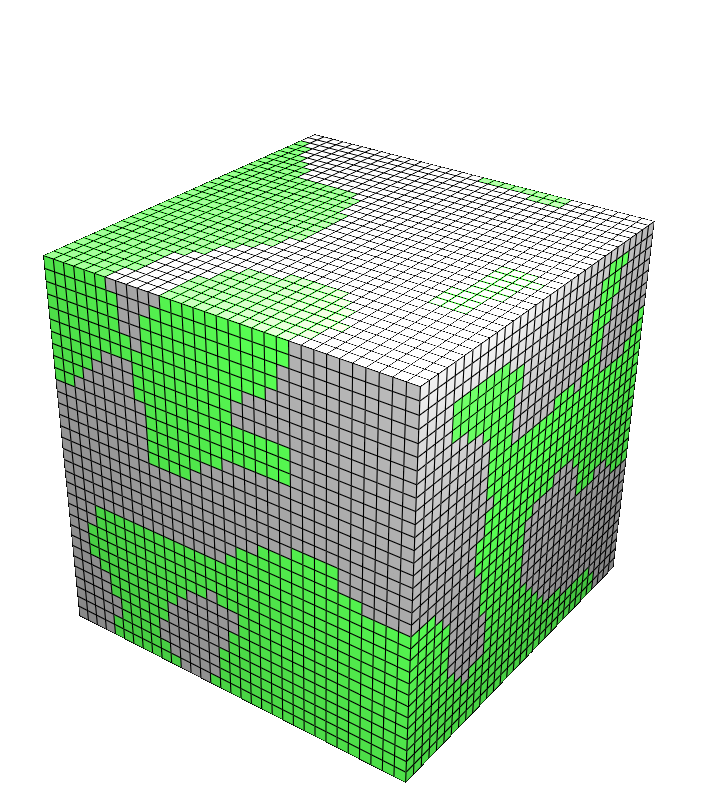}}
\put(0,0){(a)}
\put(230,0){\includegraphics[width=0.50\textwidth]{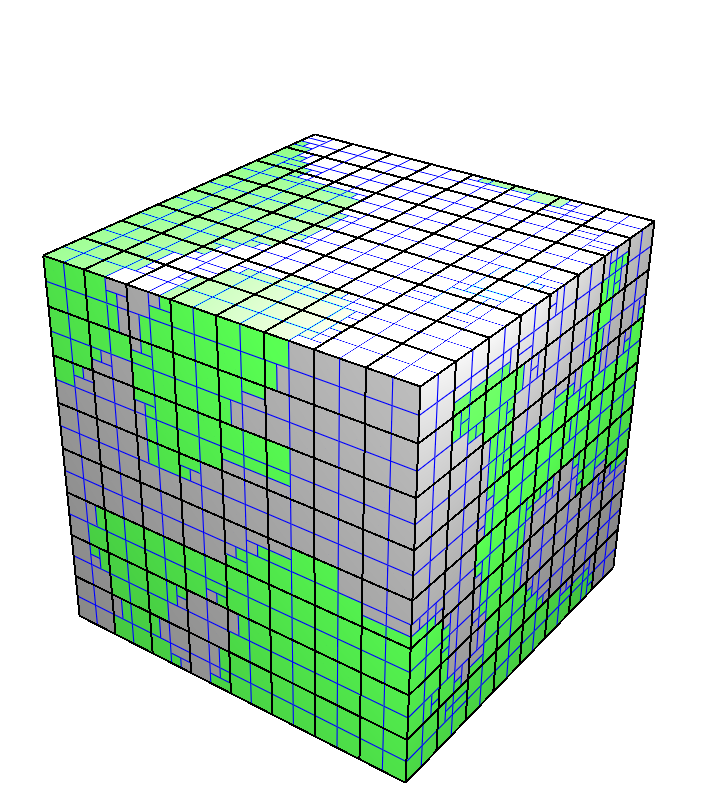}}
\put(230,0){(b)}
\end{picture}
\caption{(a) Binarized voxel data set of the considered metal-matrix composite (Ferrotitanite) containing $32\times 32\times 32$ voxels, each with the dimensions $2\, \mathrm{\mu m}\times 2\, \mathrm{\mu m} \times 2\, \mathrm{\mu m}$; the microstructure consists of brittle titanium carbide inclusions (green) and a ductile Nikro128 matrix (grey). (b) Discretization of microstructure with $8 \times 8 \times 8$ finite cells and the decomposition strategy T1min1-MT.}
\label{fig:FTN_voxel}
\end{figure}

\begin{figure}[!]
\unitlength1cm
\begin{picture}(16,24)
\put(0.0, 19.0){\rotatebox{90}{ Initial crack at $\bar{\varepsilon}_{xx}=0$}}
\put(0.0, 13.1){\rotatebox{90}{ von Mises stress~$\tau^{\mathrm{vM}}$}}
\put(0.5, 13.6){\rotatebox{90}{ at $\bar{\varepsilon}=0.017\,\%$}}
\put(0.0, 7.4){\rotatebox{90}{ Final morphology}}
\put(0.5, 7.6){\rotatebox{90}{ at $\bar{\varepsilon}_{xx}=2.5\,\%$}}
\put(0.0, 1.5){\rotatebox{90}{ Final crack path}}
\put(0.5, 1.7){\rotatebox{90}{ at $\bar{\varepsilon}_{xx}=2.5\,\%$}}
\put( 4.2, 24.0){ Initial crack 1}
\put(11.0,24.0){ Initial crack 2}
\put( 8.5,23.0){\includegraphics[width=0.1\textwidth]{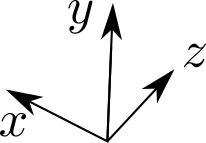}}
\put( 9.3,17.7){\includegraphics[width=0.4\textwidth]{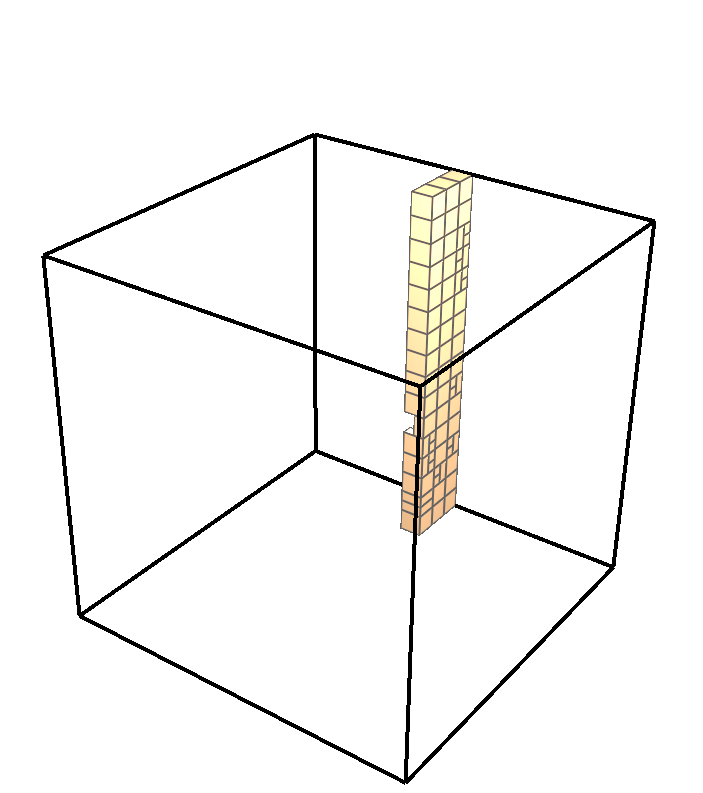}}
\put( 9.3,11.8){\includegraphics[width=0.4\textwidth]{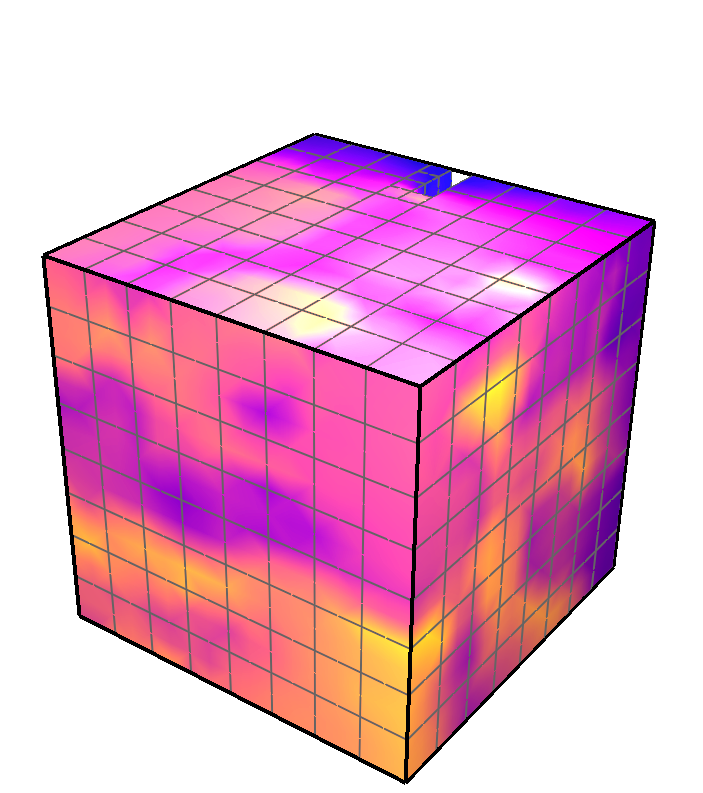}}
\put( 9.3, 5.9){\includegraphics[width=0.4\textwidth]{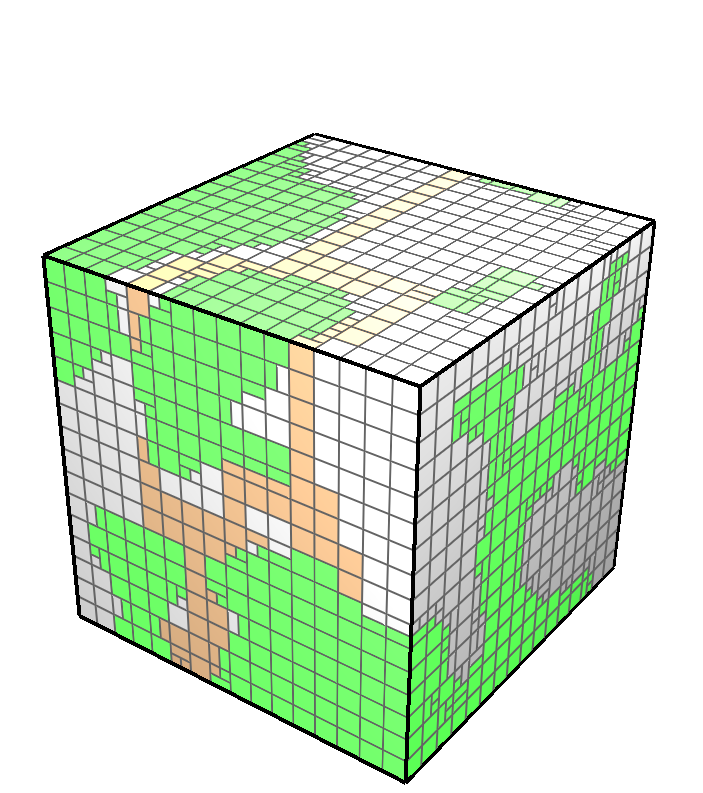}}
\put( 9.3, 0.0){\includegraphics[width=0.4\textwidth]{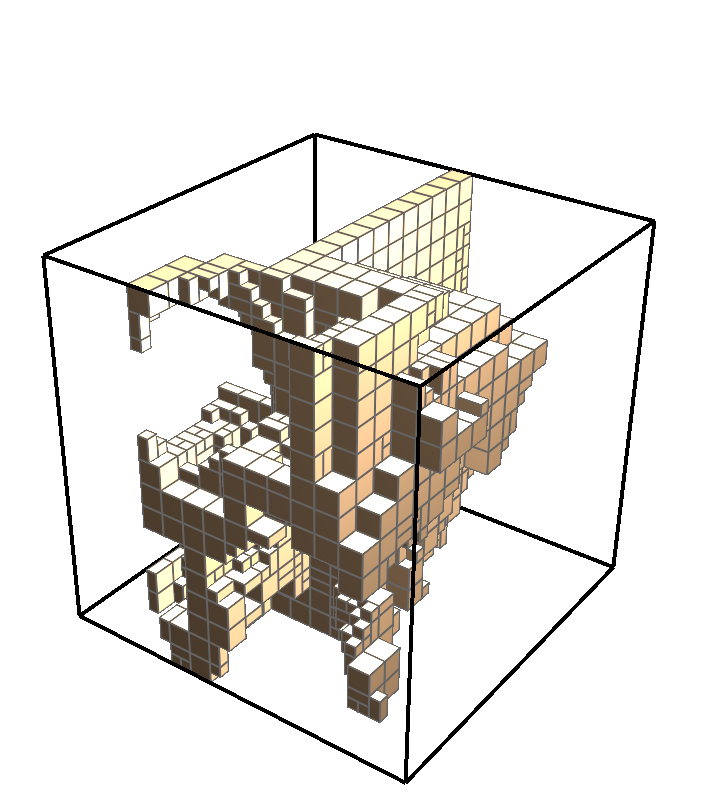}}
\put( 2.5,17.7){\includegraphics[width=0.4\textwidth]{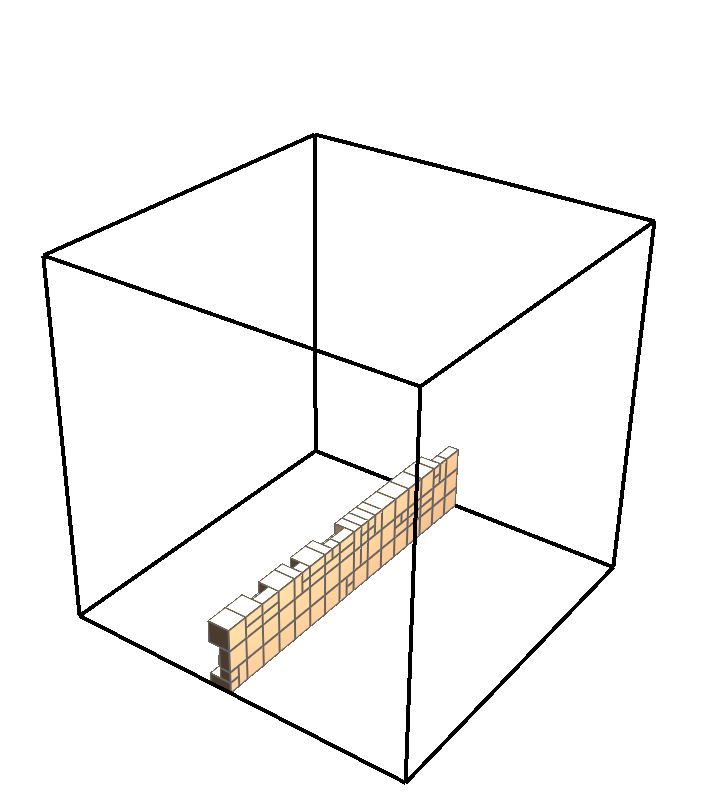}}
\put( 2.5,11.8){\includegraphics[width=0.4\textwidth]{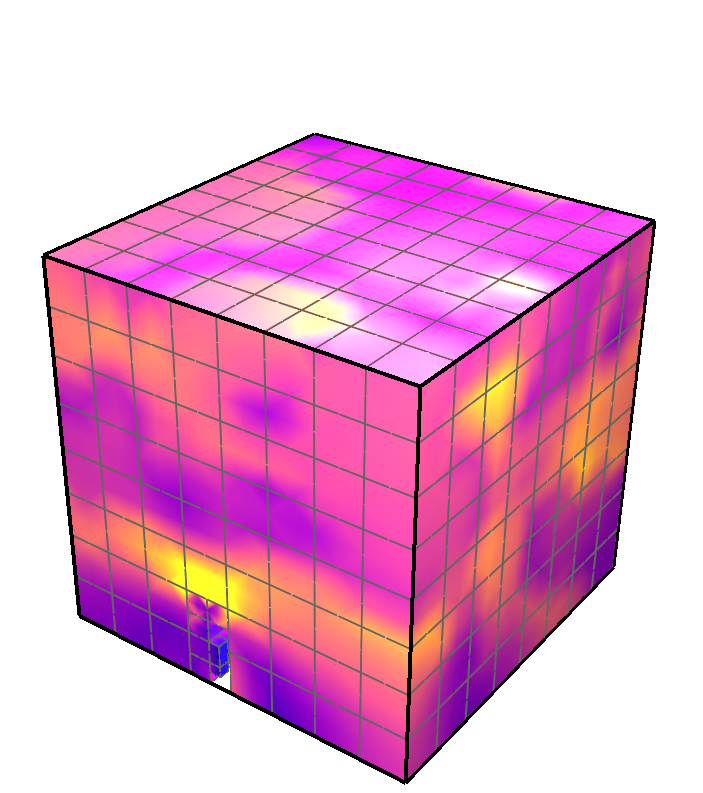}}
\put( 2.5, 5.9){\includegraphics[width=0.4\textwidth]{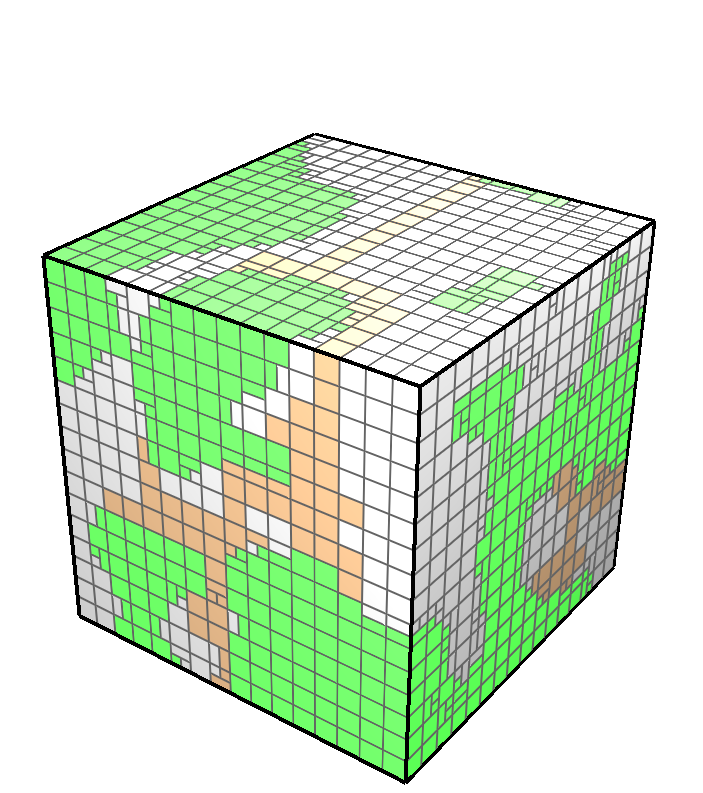}}
\put( 2.5, 0.0){\includegraphics[width=0.4\textwidth]{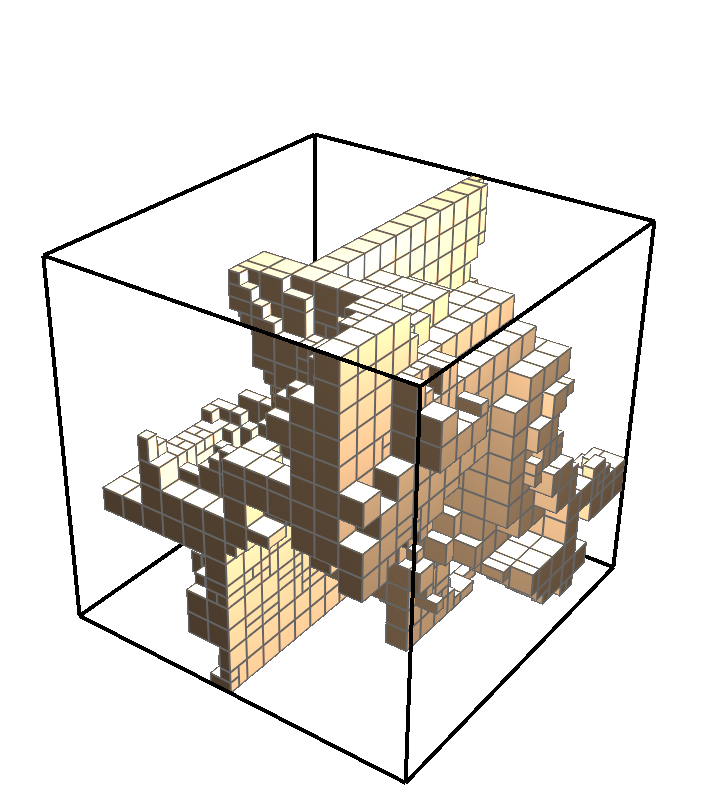}}
\put(2.1,12.6){\rotatebox{90}{ \includegraphics[width=125pt ,height=10 pt]{\figpath/colorbar_lines.png}}}
\put(2.0,17.3){\footnotesize{$\mathrm{MPa}$}}
\put(1.8,12.6){\footnotesize{$0$}}
\put(1.5,13.7){\footnotesize{$100$}}
\put(1.5,14.75){\footnotesize{$200$}}
\put(1.5,15.85){\footnotesize{$300$}}
\put(1.5,16.9){\footnotesize{$400$}}
\put(2.5,18.0){(a.1)}
\put(2.5,12.0){(b.1)}
\put(2.5, 6.2){(c.1)}
\put(2.5, 0.3){(d.1)}
\put(9.3,18.0){(a.2)}
\put(9.3,12.0){(b.2)}
\put(9.3, 6.2){(c.2)}
\put(9.3, 0.3){(d.2)}
\end{picture}
\caption{Simulated crack through metal-matrix microstructure at different specimen strains $\bar{\varepsilon}_{xx}=\bar{u}_{xx}/l_{xx}$ for initial crack (1) and (2); (a) shows the initial crack, (b) depicts the von Mises Kirchhoff stress, (c) and (d) show the microstructure and crack, respectively, for the fully fractured scenario. 
\label{fig:FTN_crack_path}}
\end{figure}

In order to demonstrate the capability of the proposed FCM/Eigenerosion approach to simulate real-world, heterogeneous structures, the real microstructure of Ferrotitanite is investigated here. 
This metal-matrix composite consists of brittle titanium carbide inclusions surrounded by a ductile Nikro128 matrix. 
The microstructure data is obtained from Micro-CT scanning and the considered data set consists of $32\times 32\times 32$ voxels, cf. figure~\ref{fig:FTN_voxel}a. 
This specimen has a size of $64\times 64 \times 64\,\mathrm{\mu m^3}$ and it discretized with $8\times 8 \times 8$ finite cells applying the scheme T1min1-OD assuming tri-quadratic shape functions cf. figure~\ref{fig:FTN_voxel}b. 
%
\begin{table}
\resizebox{\textwidth}{!}{%
\begin{tabular}{|gcccccccc|}
 \hline
\rowcolor{rub_blue40}
& $K\,[\mathrm{GPa}]$               & $\mu\,[\mathrm{GPa}]$ & $y_0\,[\mathrm{GPa}]$ & $y_{\infty}\,[\mathrm{GPa}]$ & $h^{\mathrm{exp}}\,[-]$ & $h^{\mathrm{lin}}\,[\mathrm{GPa}]$ & $\eta\,[\mathrm{GPa\,s}]$ & $G_c\,[\mathrm{N/mm}]$ \\
Titanium carbide &$235.42$ & $191.53$ & $10^{12}$ & $10^{12}$ & $0$ & $0$ & $-$ & $0.114$ \\
NiBSi & $167.84$ & $77.47$ & $1.3$ & $1.5$ & $300$ & $5.0$ & $1.0$ & $0.022$\\
\hline
\end{tabular}
}
\caption{Considered material parameters of Ferrotitanite microstructure.}
\label{tab:FTN_microstructure_material_parameters}
\end{table}
Simulations assuming an unregularized elasto-plastic material behavior have proven difficult. 
For this structure large changes in the deformation fields due to the erosion of elements occur which led to a failing Newton-Raphson iteration. 
This is not surprising as in~\cite{WinBal:2022:SoC} it has been shown that simulations including eigenerosion and an unregularized elasto-plastic material law suffered from localized plastic zones. 
To avoid these issues, the elasto-viscoplastic material law given in section~\ref{sec:matmod} is considered, where the evolution of the internal variables is delayed and thereby localization effects are decreased. 
The material parameters of the two constituents are presented in table~\ref{tab:FTN_microstructure_material_parameters}. 
In case of the titanium carbide, estimated values were given by material science expert collaborators.
In contrast to that, the parameters of the matrix material were fitted to experimental tensile tests in case of the ductile matrix.
In both cases, the Griffith-type energy release rate~$G_c$ has been fitted to tensile tests.
Here, the geometry of the tensile tests is decreased in such a way, that its axial lengths is of ~$50\,\mu$m because the Griffith-type energy release rate~$G_c$ scales with the geometry size.
This procedure ensures that this parameter lies in the correct magnitude.\\
As boundary conditions, the displacements at two opposite faces of the specimen are prescribed into the $x$-directions normal to the faces as in the previously shown benchmark experiment.
Additionally, at one side in $y$-direction and one side in $z$-direction, the displacements in normal directions are restricted, which corresponds to symmetry conditions in a uniaxial tension scenario.
The velocity of the deformation is controlled by the given strain rate $\dot{\bar{\varepsilon}}_{xx}=0.01/\mathrm{min}$ defined on specimen level. 
An initially fractured zone at the edge of the specimen is imposed to control the position of the crack to go rather through the center of the specimen in order to avoid boundary effects. 
Simulations with two different initial cracks as shown in figure~\ref{fig:FTN_crack_path}a are carried out. 
Due to the rather isotropically distributed inclusions, the considered metal-matrix composite behaves isotropically at the larger scale and thus, a quite similar response is to be expected from the simulations with varying initial cracks. 
Figure~\ref{fig:FTN_crack_path}b demonstrate that the resulting von Mises Kirchhoff stress~$\tau^{\mathrm{vM}}$ in the first time step before cracking occurs gets higher at the tip of the initial cracks and in the titanium carbide inclusions compared to the ones in the ductile matrix. 
However, the crack propagates primarily through the ductile matrix cf. figure~\ref{fig:FTN_crack_path}c because of its lower resistance against crack propagation. 
Hence, the crack path goes around the hard phase particles and propagates in all three dimensions until the structure is split into two parts, cf. figure~\ref{fig:FTN_crack_path}d. 
Additionally, bifurcation of the crack path occurs. 
For both initial crack scenarios, qualitatively similar crack paths develop. 
\begin{figure}[t]
\begin{picture}(26.60,165.0)
\put(-10,0){\includegraphics[width=0.52\textwidth]{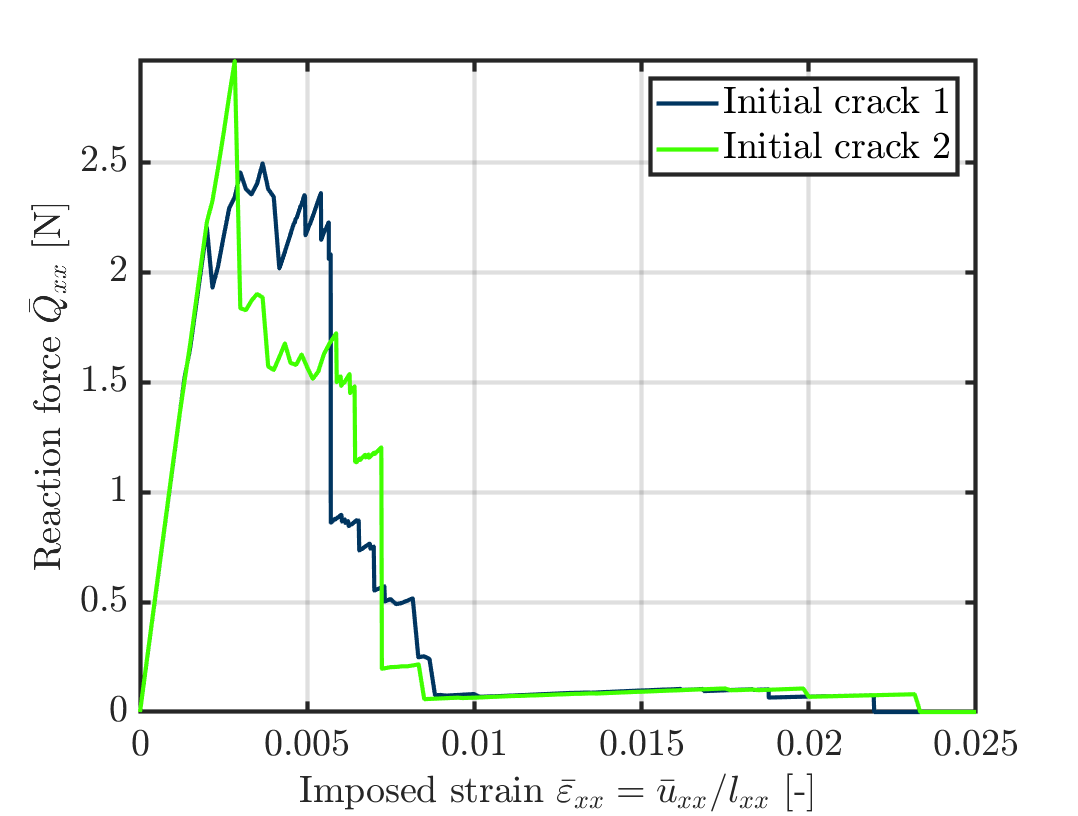}}
\put(220,0){\includegraphics[width=0.52\textwidth]{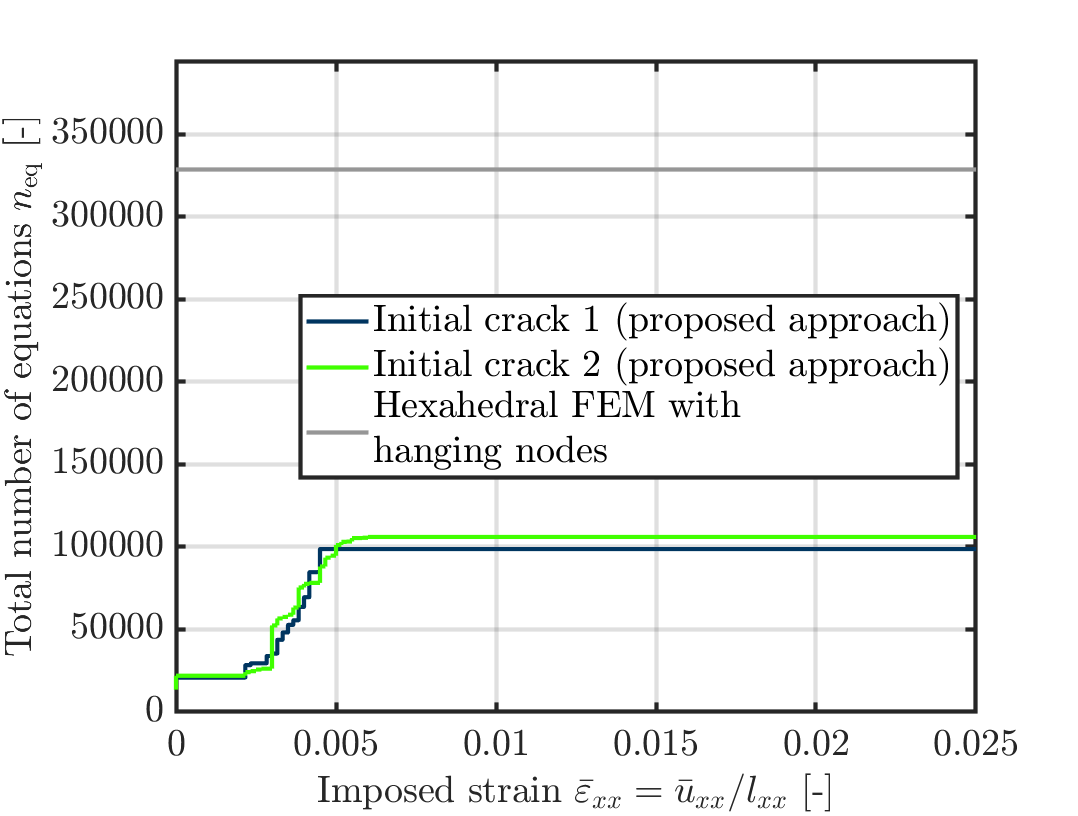}}
\put(0,0){(a)}
\put(240,0){(b)}
\end{picture}
\caption{(a) Resultant reaction force~$\bar{Q}_{xx}$ in the Dirichlet boundary in $x$-direction and (b) number of equations in the linearized system within each Newton iteration versus specimen strain $\bar{\varepsilon}_{xx}$. 
}
\label{fig:FTN_reac_eq}
\end{figure}
To analyze the quantitative response, figure~\ref{fig:FTN_reac_eq}a depicts the resultant reaction force at the Dirichlet boundary versus specimen strain. 
At first, the reaction force increases almost linearly until crack propagation starts at $\bar{\varepsilon}_{xx}=0.1\,\%$. 
Then, the reaction force decreases in some time steps due to erosion of elements and increases in some time steps due to the increased macroscopic strains until only a small connection between both sides remains at $\bar{\varepsilon}_{xx}=1\,\%$ until it breaks at $\bar{\varepsilon}_{xx}=2.2\,\%$. 
As can be seen, also quantitatively the response of the two initial crack scenarios is quite similar reflecting the overall isotropy of the material. \\
For the assessment of the efficiency, figure~\ref{fig:FTN_reac_eq}b depicts the number of linearized equations to be solved within each Newton-Raphson step versus specimen strain. 
In the final state, the number of equations $n_{\mathrm{eq}}$ is reduced by a factor of $3$ for the proposed FCM/Eigenerosion strategy compared to the hexahedral FEM with hanging nodes. 
This is not surprising since in the final state, $1,271$ (initial crack scenario 1) and $1,294$ (initial crack scenario 2) subcells have been eroded, resulting in approximately $1/6$ of the total number of subcells ($6,920$). 
In the beginning of the crack propagation, where only a few subcells are switched, the factor in the reduction of linearized equations is even $12$. 
Thus, by solely switching those subcells to finite elements where the crack propagates a significant efficiency gain is achieved. 
If each voxel was simulated as one single finite element, approximately $815,000$ linearized equations would occur which surpasses the one of the proposed strategy by a factor of $8$ in the final state and $24$ in the beginning. 
Additionally, the computational effort for the assembling would be highly increased because of the number of elements $32^3=32,768$ is high compared to the $6,920$ subcells of the enhanced FCM. 
Summarizing, this numerical example demonstrates the computational efficiency of the proposed FC/Eigenerosion strategy and its capability to simulate complex crack paths based on real voxel data.

\section{Conclusion}
The aim of this work was to develop a mesh-independent, robust and efficient framework for the simulation of ductile crack propagation through heterogeneous structures given as voxel data.
Therefore, the basic formulation of the eigenerosion for ductile crack propagation at finite strains has been chosen as a basis.
For an efficient simulation based on voxel data, the Finite Cell Method (FCM) has been incorporated additionally.
In our combined approach, the mesh is refined at the crack tip by transforming the subcells into single finite elements.
Thereby, although introducing an adaptive mesh-refinement in terms of the FE approximation, the integration scheme is kept and thus, additional projections of history variables are avoided.
For the application of this framework, the decomposition of the finite cells representing the voxel data as efficiently as possible has been a crucial ingredient.
In a numerical example of an artificial metal matrix composite, a microstructure consisting of a brittle sphere surrounded by a ductile matrix, it was shown that the octree decomposition with a merge on the lowest level led to the decomposition resulting in the most efficient simulation.
Furthermore, this exemplary microstructure was discretized with varying approaches.
In a virtual tensile test, a quite similar crack path occurred in all simulations.
By comparing the resultant forces, more pronounced quantitative differences between the calculations have been observed, specifically in the ductile softening phase of the simulation.
However, the results obtained from our proposed approach converged (with increasing numbers of finite cells) to the solution of a semi-regular hexahedral mesh using hanging nodes.
In this latter discretization, all material interfaces were locally adaptively meshed such that each of these elements was in line with the subcells considered as part of the proposed approach.
Therefore, the proposed approach turned out to be qualitatively and quantitatively as accurate as the semi-regular hexahedral discretization, which was not surprising.
However, as major advantage, our approach only required the full locally refined semi-regular mesh where the crack takes place.
In this example all material interphases, where the advantage of the FCM could be exploited, were eroded, and thus, the gain in computational costs could only be moderate.
Therefore, a second numerical example has been analyzed where a part of a real metal matrix composite microstructure obtained from micro-CT was considered.
In this example, the number of equations in the linearized system within the Newton-Raphson scheme could be reduced by a factor of approximately 3, compared to the semi-regular hexahedral discretization.
Comparing the calculations for two differently located initial cracks, a similar structural response and crack path was obtained reflecting the material's macroscopic isotropy.

\section*{Acknowledgement}

The authors gratefully appreciate financial funding by the Deutsche Forschungsgemeinschaft (DFG) in the framework of the Collaborative Research Center SFB 837 ``Interaction Modeling in Mechanized Tunneling'', project C6. 
Furthermore, the authors thank Arne R\"ottger (Bergische Universit\"at Wuppertal, Germany) and Lukas Brackmann (Ruhr-Universit\"at Bochum, Germany) for providing experimental data regarding the individual micro-components of the metallic microstructures. Furthermore,  we thank the Zentrum für Grenzfl\"achendominierte H\"ochstleistungswerkstoffe, ZGH (Ruhr-Universit\"at Bochum, Germany) for providing the micro-CT data.

\bibliographystyle{stylefiles/references_style}
\bibliography{references}

\begin{thebibliography}{52}
\providecommand{\natexlab}[1]{#1}
\providecommand{\url}[1]{\texttt{#1}}
\expandafter\ifx\csname urlstyle\endcsname\relax
  \providecommand{\doi}[1]{doi: #1}\else
  \providecommand{\doi}{doi: \begingroup \urlstyle{rm}\Url}\fi

\bibitem[Barenblatt(1962)]{Bar:1962:Tmt}
G.~I. Barenblatt.
\newblock The mathematical theory of equilibrium cracks in brittle fracture.
\newblock In \emph{Advances in applied mechanics}, volume~7, pages 55--129.
  Elsevier, 1962.
\newblock \doi{10.1016/S0065-2156(08)70121-2}.

\bibitem[Beese et~al.(2018)Beese, Loehnert, and Wriggers]{BeeLoeWri:2018:3dc}
S.~Beese, S.~Loehnert, and P.~Wriggers.
\newblock 3d ductile crack propagation within a polycrystalline microstructure
  using xfem.
\newblock \emph{Computational Mechanics}, 61\penalty0 (1):\penalty0 71--88,
  2018.
\newblock \doi{10.1007/s00466-017-1427-y}.

\bibitem[Belytschko and Black(1999)]{BelBla:1999:ecg}
T.~Belytschko and T.~Black.
\newblock Elastic crack growth in finite elements with minimal remeshing.
\newblock \emph{International journal for numerical methods in engineering},
  45\penalty0 (5):\penalty0 601--620, 1999.

\bibitem[Cheng et~al.(2020)Cheng, Tu, and Ghosh]{CheTuGho:2020:Wea}
J.~Cheng, X.~Tu, and S.~Ghosh.
\newblock Wavelet-enriched adaptive hierarchical fe model for coupled crystal
  plasticity-phase field modeling of crack propagation in polycrystalline
  microstructures.
\newblock \emph{Computer Methods in Applied Mechanics and Engineering},
  361:\penalty0 112757, 2020.
\newblock \doi{10.1016/j.cma.2019.112757}.

\bibitem[Demkowicz et~al.(1989)Demkowicz, Oden, Rachowicz, and
  Hardy]{DemOdeRacHar:1989:tau}
L.~Demkowicz, J.~T. Oden, W.~Rachowicz, and O.~Hardy.
\newblock Toward a universal hp adaptive finite element strategy, part 1.
  {C}onstrained approximation and data structure.
\newblock \emph{Computer Methods in Applied Mechanics and Engineering},
  77\penalty0 (1-2):\penalty0 79--112, 1989.
\newblock \doi{10.1016/0045-7825(89)90129-1}.

\bibitem[Dugdale(1960)]{Dug:1960:Yos}
D.~S. Dugdale.
\newblock Yielding of steel sheets containing slits.
\newblock \emph{Journal of the Mechanics and Physics of Solids}, 8\penalty0
  (2):\penalty0 100--104, 1960.
\newblock \doi{10.1016/0022-5096(60)90013-2}.

\bibitem[D{\"u}ster et~al.(2008)D{\"u}ster, Parvizian, Yang, and
  Rank]{DueParYanRan:2008:Tfc}
A.~D{\"u}ster, J.~Parvizian, Z.~Yang, and E.~Rank.
\newblock The finite cell method for three-dimensional problems of solid
  mechanics.
\newblock \emph{Computer methods in applied mechanics and engineering},
  197\penalty0 (45-48):\penalty0 3768--3782, 2008.
\newblock \doi{10.1016/j.cma.2008.02.036}.

\bibitem[Fangye et~al.(2020)Fangye, Miska, and Balzani]{FanMisBal:2020:aso}
Y.~F. Fangye, N.~Miska, and D.~Balzani.
\newblock Automated simulation of voxel-based microstructures based on enhanced
  finite cell approach.
\newblock \emph{Archive of Applied Mechanics}, pages 1--19, 2020.
\newblock \doi{10.1007/s00419-020-01719-x}.

\bibitem[Griffith(1921)]{Gri:1921:tpo}
A.~A. Griffith.
\newblock The phenomena of rupture and flow in solids.
\newblock \emph{Philosophical Transactions of the Royal Society A:
  Mathematical, Physical and Engineering Sciences}, 221\penalty0
  (582-593):\penalty0 163--198, 1921.
\newblock \doi{10.1098/rsta.1921.0006}.

\bibitem[Hillerborg et~al.(1976)Hillerborg, Mod{\'e}er, and
  Petersson]{HilModPet:1976:Aoc}
A.~Hillerborg, M.~Mod{\'e}er, and P.-E. Petersson.
\newblock Analysis of crack formation and crack growth in concrete by means of
  fracture mechanics and finite elements.
\newblock \emph{Cement and concrete research}, 6\penalty0 (6):\penalty0
  773--781, 1976.

\bibitem[Irwin(1957)]{Irw:1957:aos}
G.~Irwin.
\newblock Analysis of stresses and strains near the end of a crack traversing a
  plate.
\newblock \emph{Journal of Applied Mechanics}, 24:\penalty0 361--364, 1957.

\bibitem[Junker et~al.(2017)Junker, Schwarz, Makowski, and
  Hackl]{JunSchMakHac:2017:cmt}
P.~Junker, S.~Schwarz, J.~Makowski, and K.~Hackl.
\newblock A relaxation-based approach to damage modeling.
\newblock \emph{Continuum Mechanics and Therodynamics}, 29\penalty0
  (1):\penalty0 291--310, 2017.
\newblock \doi{10.1007/s00161-016-0528-8}.

\bibitem[Junker et~al.(2019)Junker, Schwarz, Jantos, and
  Hackl]{JunSchJanHac:2019:afa}
P.~Junker, S.~Schwarz, D.~Jantos, and K.~Hackl.
\newblock A fast and robust numerical treatment of a gradient-enhanced model
  for brittle damage.
\newblock \emph{International Journal for Multiscale Computational
  Engineering}, 17, 2019.
\newblock \doi{10.1615/IntJMultCompEng.2018027813}.

\bibitem[Junker et~al.(submitted)Junker, Riesselmann, and
  Balzani]{JunRieBal:2021:ear}
P.~Junker, J.~Riesselmann, and D.~Balzani.
\newblock Efficient and robust numerical treatment of a gradient-enhanced
  damage model at large deformations.
\newblock \emph{International Journal for Numerical Methods in Engineering},
  submitted.

\bibitem[Klinkel(2000)]{Kli:2000:tun}
S.~Klinkel.
\newblock \emph{Theorie und Numerik eines Volumen-Schalen-Elementes bei finiten
  elastischen und plastischen Verzerrungen}.
\newblock PhD thesis, Institut für Baustatik, Universit\"at Karlsruhe, 2000.

\bibitem[Kr{\"o}ner(1959)]{Kro:1959:akd}
E.~Kr{\"o}ner.
\newblock Allgemeine {K}ontinuumstheorie der {V}ersetzungen und
  {E}igenspannungen.
\newblock \emph{Archive for Rational Mechanics and Analysis}, 4\penalty0
  (1):\penalty0 273, 1959.

\bibitem[Lee(1969)]{Lee:1969:epd}
E.~H. Lee.
\newblock Elastic-plastic deformation at finite strains.
\newblock \emph{Journal of Applied Mechanics}, 36\penalty0 (1):\penalty0 1--6,
  1969.
\newblock \doi{10.1115/1.3564580}.

\bibitem[Li et~al.(2015)Li, Pandolfi, and Ortiz]{LiPanOrt:2015:mom}
B.~Li, A.~Pandolfi, and M.~Ortiz.
\newblock Material point erosion simulation of dynamic fragmentation of metals.
\newblock \emph{Mechanics of Materials}, 80:\penalty0 288--297, 2015.
\newblock \doi{10.1007/Fs10704-012-9788-x}.

\bibitem[Liang and Sofronis(2003)]{LiaSof:2003:Man}
Y.~Liang and P.~Sofronis.
\newblock Micromechanics and numerical modelling of the
  hydrogen--particle--matrix interactions in nickel-base alloys.
\newblock \emph{Modelling and Simulation in Materials Science and Engineering},
  11\penalty0 (4):\penalty0 523, 2003.
\newblock \doi{10.1088/0965-0393/11/4/308}.

\bibitem[Meng and Wang(2015)]{MenWan:2015:Poi}
Q.~Meng and Z.~Wang.
\newblock Prediction of interfacial strength and failure mechanisms in
  particle-reinforced metal-matrix composites based on a micromechanical model.
\newblock \emph{Engineering Fracture Mechanics}, 142:\penalty0 170--183, 2015.
\newblock \doi{10.1016/j.engfracmech.2015.06.001}.

\bibitem[Miehe et~al.(1994)Miehe, Stein, and Wagner]{MieSteWag:1994:ame}
C.~Miehe, E.~Stein, and W.~Wagner.
\newblock Associative multiplicative elasto-plasticity: formulation and aspects
  of the numerical implementation including stability analysis.
\newblock \emph{Computers \& structures}, 52\penalty0 (5):\penalty0 969--978,
  1994.
\newblock \doi{10.1016/0045-7949(94)0081-7}.

\bibitem[Miehe et~al.(2010)Miehe, Welschinger, and
  Hofacker]{MieWelHof:2010:ijnm}
C.~Miehe, F.~Welschinger, and M.~Hofacker.
\newblock Thermodynamically consistent phase-field models of fracture:
  Variational principles and multi-field fe implementations.
\newblock \emph{International Journal for Numerical Methods in Engineering},
  83:\penalty0 1273--1311, 2010.
\newblock \doi{10.1002/nme.2861}.

\bibitem[Mielke and Ortiz(2008)]{MieOrt:2008:aco}
A.~Mielke and M.~Ortiz.
\newblock A class of minimum principles for characterizing the trajectories and
  the relaxation of dissipative systems.
\newblock \emph{ESAIM: Control, Optimisation and Calculus of Variations},
  14\penalty0 (3):\penalty0 494--516, 2008.
\newblock \doi{10.1051/cocv:2007064}.

\bibitem[Miska and Balzani(2019)]{MisBal:2019:qou}
N.~Miska and D.~Balzani.
\newblock Quantification of uncertain macroscopic material properties resulting
  from variations of microstructure morphology based on statistically similar
  volume elements - application to dual-phase steels.
\newblock \emph{Computational Mechanics}, 64:\penalty0 1621--1637, 2019.
\newblock \doi{https://doi.org/10.1007/s00466-019-01738-8}.

\bibitem[Nagaraja et~al.(2019)Nagaraja, Elhaddad, Ambati, Kollmannsberger,
  De~Lorenzis, and Rank]{NagElhAmbKolDeLRan:2019:Pfm}
S.~Nagaraja, M.~Elhaddad, M.~Ambati, S.~Kollmannsberger, L.~De~Lorenzis, and
  E.~Rank.
\newblock Phase-field modeling of brittle fracture with multi-level hp-fem and
  the finite cell method.
\newblock \emph{Computational mechanics}, 63\penalty0 (6):\penalty0 1283--1300,
  2019.
\newblock \doi{10.1007/s00466-018-1649-7}.

\bibitem[Navas et~al.(2018)Navas, Rena, Li, and Ruiz]{NavRenYuLiRui:2018:mtd}
P.~Navas, C.~Y. Rena, B.~Li, and G.~Ruiz.
\newblock Modeling the dynamic fracture in concrete: an eigensoftening meshfree
  approach.
\newblock \emph{International Journal of impact engineering}, 113:\penalty0
  9--20, 2018.
\newblock \doi{10.1016/j.ijimpeng.2017.11.004}.

\bibitem[Newmark(1959)]{New:1959:amc}
N.~M. Newmark.
\newblock A method of computation for structural dynamics.
\newblock \emph{Journal of the engineering mechanics division}, 85\penalty0
  (3):\penalty0 67--94, 1959.

\bibitem[Nguyen et~al.(2016)Nguyen, Yvonnet, Zhu, Bornert, and
  Chateau]{NguYvoZhuBor:2016:Apf}
T.-T. Nguyen, J.~Yvonnet, Q.-Z. Zhu, M.~Bornert, and C.~Chateau.
\newblock A phase-field method for computational modeling of interfacial damage
  interacting with crack propagation in realistic microstructures obtained by
  microtomography.
\newblock \emph{Computer Methods in Applied Mechanics and Engineering},
  312:\penalty0 567--595, 2016.
\newblock \doi{10.1016/j.cma.2015.10.007}.

\bibitem[Oden et~al.(1989)Oden, Demkowicz, Rachowicz, and
  Westermann]{OdeDemRacWes:1989:tau}
J.~T. Oden, L.~Demkowicz, W.~Rachowicz, and T.~Westermann.
\newblock Toward a universal hp adaptive finite element strategy, part 2. {A}
  posteriori error estimation.
\newblock \emph{Computer methods in applied mechanics and engineering},
  77\penalty0 (1-2):\penalty0 113--180, 1989.
\newblock \doi{0.1016/0045-7825(89)90130-8}.

\bibitem[Pandolfi and Ortiz(2012)]{PanOrt:2012:aea}
A.~Pandolfi and M.~Ortiz.
\newblock An eigenerosion approach to brittle fracture.
\newblock \emph{International Journal for Numerical Methods in Engineering},
  92\penalty0 (8):\penalty0 694--714, 2012.
\newblock \doi{10.1002/nme.4352}.

\bibitem[Pandolfi et~al.(2013)Pandolfi, Li, and Ortiz]{PanLiOrt:2013:ijf}
A.~Pandolfi, B.~Li, and M.~Ortiz.
\newblock Modeling failure of brittle materials with eigenerosion.
\newblock \emph{Computational Modelling of Concrete Structures}, 1:\penalty0
  9--21, 2013.
\newblock \doi{10.1007/s10704-012-9788-x}.

\bibitem[Parvizian et~al.(2007)Parvizian, Düster, and
  Rank]{ParDueRan:2007:fcm}
J.~Parvizian, A.~Düster, and E.~Rank.
\newblock Finite cell method - h- and p-extension for embedded domain problems
  in solid mechanics.
\newblock \emph{Computational Mechanics}, 41:\penalty0 121 -- 133, 2007.
\newblock \doi{10.1007/s00466-007-0173-y}.

\bibitem[Perzyna(1966)]{Per:1966:Fpi}
P.~Perzyna.
\newblock Fundamental problems in viscoplasticity.
\newblock In \emph{Advances in applied mechanics}, volume~9, pages 243--377.
  Elsevier, 1966.
\newblock \doi{10.1016/S0065-2156(08)70009-7}.

\bibitem[Qinami et~al.(2019)Qinami, Pandolfi, and Kaliske]{QinPanKal:2019:vef}
A.~Qinami, A.~Pandolfi, and M.~Kaliske.
\newblock Variational eigenerosion for rate dependent plasticity in concrete
  modelling at small strain.
\newblock \emph{International Journal for Numerical Methods in Engineering},
  2019.
\newblock \doi{10.1002/nme.6271}.

\bibitem[Rachowicz et~al.(1989)Rachowicz, Oden, and
  Demkowicz]{RacOdeDem:1989:tau}
W.~Rachowicz, J.~T. Oden, and L.~Demkowicz.
\newblock Toward a universal hp adaptive finite element strategy part 3.
  {D}esign of hp meshes.
\newblock \emph{Computer Methods in Applied Mechanics and Engineering},
  77\penalty0 (1-2):\penalty0 181--212, 1989.
\newblock \doi{10.1016/0045-7825(89)90131-X}.

\bibitem[Ranjbar et~al.(2014)Ranjbar, Mashayekhi, Parvizian, D{\"u}ster, and
  Rank]{RanMasParDusRan:2014:Utf}
M.~Ranjbar, M.~Mashayekhi, J.~Parvizian, A.~D{\"u}ster, and E.~Rank.
\newblock Using the finite cell method to predict crack initiation in ductile
  materials.
\newblock \emph{Computational Materials Science}, 82:\penalty0 427--434, 2014.
\newblock \doi{10.1016/j.commatsci.2013.10.012}.

\bibitem[Schellekens and De~Borst(1993)]{SchDeB:1993:Otn}
J.~Schellekens and R.~De~Borst.
\newblock On the numerical integration of interface elements.
\newblock \emph{International Journal for Numerical Methods in Engineering},
  36\penalty0 (1):\penalty0 43--66, 1993.
\newblock \doi{10.1002/nme.1620360104}.

\bibitem[Scheunemann et~al.(2015)Scheunemann, Balzani, Brands, and
  Schr\"oder]{SchBalBraSch:2015:do3}
L.~Scheunemann, D.~Balzani, D.~Brands, and J.~Schr\"oder.
\newblock Design of {3D} statistically similar representative volume elements
  based on {M}inkowski functionals.
\newblock \emph{Mech. Mater.}, 90:\penalty0 185--201, 2015.

\bibitem[Schillinger et~al.(2012)Schillinger, Ruess, Zander, Bazilevs,
  D{\"u}ster, and Rank]{SchRueZanBazDueRan:2012:Sal}
D.~Schillinger, M.~Ruess, N.~Zander, Y.~Bazilevs, A.~D{\"u}ster, and E.~Rank.
\newblock Small and large deformation analysis with the p-and b-spline versions
  of the finite cell method.
\newblock \emph{Computational Mechanics}, 50\penalty0 (4):\penalty0 445--478,
  2012.
\newblock \doi{10.1007/s00466-012-0684-z}.

\bibitem[Schmidt et~al.(2009)Schmidt, Fraternali, and
  Ortiz]{SchFraOrt:2009:eae}
B.~Schmidt, F.~Fraternali, and M.~Ortiz.
\newblock Eigenfracture: an eigendeformation approach to variational fracture.
\newblock \emph{Multiscale Modeling \& Simulation}, 7\penalty0 (3):\penalty0
  1237--1266, 2009.
\newblock \doi{10.1137/080712568}.

\bibitem[Schneider et~al.(2016)Schneider, Klusemann, and
  Bargmann]{SchKluBarg:2016:atd}
K.~Schneider, B.~Klusemann, and S.~Bargmann.
\newblock Automatic three-dimensional geometry and mesh generation of periodic
  representative volume elements for matrix-inclusion composites.
\newblock \emph{Advances in Engineering Software}, 99:\penalty0 177--188, 2016.
\newblock \doi{10.1016/j.advengsoft.2016.06.001}.

\bibitem[Shahba and Ghosh(2019)]{ShaGho:2019:Cpf}
A.~Shahba and S.~Ghosh.
\newblock Coupled phase field finite element model for crack propagation in
  elastic polycrystalline microstructures.
\newblock \emph{International Journal of Fracture}, 219\penalty0 (1):\penalty0
  31--64, 2019.
\newblock \doi{10.1007/s10704-019-00378-6}.

\bibitem[Shakoor et~al.(2019)Shakoor, Navas, Mu{\~n}oz, Bernacki, and
  Bouchard]{ShaNavMunBerBou:2019:Cmf}
M.~Shakoor, V.~M.~T. Navas, D.~P. Mu{\~n}oz, M.~Bernacki, and P.-O. Bouchard.
\newblock Computational methods for ductile fracture modeling at the
  microscale.
\newblock \emph{Archives of Computational Methods in Engineering}, 26\penalty0
  (4):\penalty0 1153--1192, 2019.
\newblock \doi{10.1007/s11831-018-9276-1}.

\bibitem[Simo(1992)]{Sim:1992:afs}
J.~C. Simo.
\newblock Algorithms for static and dynamic multiplicative plasticity that
  preserve the classical return mapping schemes of the infinitesimal theory.
\newblock \emph{Computer Methods in Applied Mechanics and Engineering},
  99\penalty0 (1):\penalty0 61--112, 1992.
\newblock \doi{10.1016/0045-7825(92)90170-O}.

\bibitem[Simo and Miehe(1992)]{SimMie:1992:cmam}
J.~C. Simo and C.~Miehe.
\newblock Associative coupled thermoplasticity at finite strains: Formulation,
  numerical analysis and implementation.
\newblock \emph{Computer Methods in Applied Mechanics and Engineering},
  98\penalty0 (1):\penalty0 41--104, 1992.
\newblock \doi{10.1016/0045-7825(92)0170-O}.

\bibitem[Simo et~al.(1985)Simo, Taylor, and Pister]{SimTay:1984:cmam}
J.~C. Simo, R.~L. Taylor, and K.~S. Pister.
\newblock Variational and projection methods for the volume constraint in
  finite deformation elasto-plasticity.
\newblock \emph{Computer methods in applied mechanics and engineering},
  51:\penalty0 177 -- 208, 1985.

\bibitem[Sukumar et~al.(2003)Sukumar, Srolovitz, Baker, and
  Prevost]{SukSroBakPre:2003:Bfi}
N.~Sukumar, D.~Srolovitz, T.~Baker, and J.-H. Prevost.
\newblock Brittle fracture in polycrystalline microstructures with the extended
  finite element method.
\newblock \emph{International Journal for Numerical Methods in Engineering},
  56\penalty0 (14):\penalty0 2015--2037, 2003.
\newblock \doi{10.1002/nme.653}.

\bibitem[Voce(1955)]{Voc:1955:aps}
E.~Voce.
\newblock A practical strain hardening function.
\newblock \emph{Metallurgia}, 51:\penalty0 219--226, 1955.

\bibitem[Wingender and Balzani(2022)]{WinBal:2022:SoC}
D.~Wingender and D.~Balzani.
\newblock Simulation of crack propagation based on eigenerosion in brittle and
  ductile materials subject to finite strains.
\newblock \emph{Archive of Applied Mechanics}, pages 1--23, 2022.
\newblock \doi{10.1007/s00419-021-02101-1}.

\bibitem[Wulf et~al.(1996)Wulf, Steinkopff, and Fischmeister]{WulFis:1996:Fso}
J.~Wulf, T.~Steinkopff, and H.~Fischmeister.
\newblock Fe-simulation of crack paths in the real microstructure of an {Al}
  (6061)/{SiC} composite.
\newblock \emph{Acta materialia}, 44\penalty0 (5):\penalty0 1765--1779, 1996.
\newblock \doi{10.1016/1359-6454(95)00328-2}.

\bibitem[Xu and Needleman(1994)]{XuNee:1994:Nso}
X.-P. Xu and A.~Needleman.
\newblock Numerical simulations of fast crack growth in brittle solids.
\newblock \emph{Journal of the Mechanics and Physics of Solids}, 42\penalty0
  (9):\penalty0 1397--1434, 1994.
\newblock \doi{10.1016/0022-5096(94)90003-5}.

\bibitem[Yang et~al.(2012)Yang, Ruess, Kollmannsberger, D{\"u}ster, and
  Rank]{YanRuesKollDueRan:2012:Aei}
Z.~Yang, M.~Ruess, S.~Kollmannsberger, A.~D{\"u}ster, and E.~Rank.
\newblock An efficient integration technique for the voxel-based finite cell
  method.
\newblock \emph{International Journal for Numerical Methods in Engineering},
  91\penalty0 (5):\penalty0 457--471, 2012.
\newblock \doi{10.1002/nme.4269}.

\end{thebibliography}

\end{document}